\documentclass[12pt]{article}
\usepackage{graphicx}
\usepackage{amscd, amssymb}

\newenvironment{eq}{\begin{equation}}{\end{equation}}
\newenvironment{proof}{{\bf Proof}:}{\vskip 5mm }

\newtheorem{proposition}{Proposition}[section]
\newtheorem{lemma}[proposition]{Lemma}
\newtheorem{definition}[proposition]{Definition}
\newtheorem{theorem}[proposition]{Theorem}
\newtheorem{cor}[proposition]{Corollary}

\newtheorem{example}[proposition]{Example}
\newtheorem{remark}[proposition]{Remark}

\newcommand{\llabel}[1]{\label{#1}}
\newcommand{\comment}[1]{}
\newcommand{\sr}{\rightarrow}

\newcommand{\bdl}{\bar{\Delta}}
\newcommand{\zz}{{\bf Z\rm}}

\newcommand{\nn}{{\bf N\rm}}
\newcommand{\qq}{{\bf Q\rm}}

\newcommand{\LL}{{\bf L}}
\newcommand{\RR}{{\bf R}} 
\newcommand{\oo}{\otimes}

\newcommand{\uu}{\underline}

\newcommand{\af}{{\bf A}^1}

\newcommand{\dsr}{\stackrel{\sr}{\scriptstyle\sr}}

\newcommand{\BB}{_{\bullet}}

{\vskip
3mm}
\begin{document}
\begin{center}
{\bf\Large Motivic Eilenberg-MacLane spaces}
\vspace {3mm}

{\large\bf Vladimir Voevodsky}\footnote{School of Mathematics, Institute for Advanced Study,
Princeton NJ, USA. e-mail: vladimir@ias.edu}$^,$\footnote{Work on the earlier versions of this  paper was supported by NSF grant 0403367.}
\vspace {3mm}

{September 21, 2009}  
\end{center}
\begin{abstract}
In this paper we construct symmetric powers in the motivic homotopy categories of morphisms  and  finite correspondences associated with $f$-admissible subcategories in the categories of schemes of finite type over a field. Using this construction we provide a description of the motivic Eilenberg-MacLane spaces representing motivic cohomology on some $f$-admissible categories including  the category of semi-normal quasi-projective schemes and, over fields which admit resolution of singularities, on some admissible subcategories including the category of smooth schemes.  This description is then used to give a complete computation of the algebra of bistable motivic cohomological operations on smooth schemes over  fields of characteristic zero and to obtain partial results on unstable operations which are required for the proof of the Bloch-Kato conjecture. 
\end{abstract}

{\em \small MSC2010 14F42, 19E15}

\tableofcontents
\parskip = 0.2in

\setcounter{section}{-1}
\section{Introduction}
In this paper we analyze the structure of the motivic Eilenberg-MacLane spaces $K(A,p,q)_C$ for $p\ge q$.  From the perspective of the general motivic homotopy theory its most important result is Theorem \ref{mainst} which asserts that for a field $k$ of characteristic zero the algebra of all bistable operations in the motivic cohomology on $\Delta^{op}(Sm/k)_+^{\#}$ with coefficients in $\zz/l$ coincides with the motivic Steenrod algebra ${\cal A}^{*,*}(k,\zz/l)$ which was introduced in \cite{Redpub}. From the point of view of the proof of the Bloch-Kato conjecture which requires information about the unstable operations its most important results are Theorems \ref{main2}, \ref{2009.09.12.th1} and \ref{mth1}.   

The definition of motivic Eilenberg-MacLane spaces appears only in the third part of the paper (Section \ref{motEM}) after a number of techniques necessary for our analysis of these spaces has been developed. Following the general outline of the paper we start the introduction with the description of the results of the first two parts (Sections \ref{section2} and \ref{section3}). 

Let us fix a field $k$.  Let $c(k)=1$ if $char(k)=0$ and $c(k)=char(k)$ if $char(k)>0$. The number $c(k)$ is sometimes known as the "characteristic exponent" of $k$.  A full subcategory $C$ of $Sch/k$ such that
\begin{enumerate}
\item $Spec(k)$ and $\af$ are in $C$
\item for $X$ and $Y$ in $C$ the product $X\times Y$ is in $C$
\item if $X$ is in $C$ and $U\sr X$ is etale then $U$ is in $C$
\item for $X$ and $Y$ in $C$ the coproduct $X\amalg Y$ is in $C$
\end{enumerate}
will be called admissible. If in addition $C$ is closed under the formation of quotients with respect to actions of finite groups it will be called $f$-admissible. The category $Sm/k$ of smooth schemes over $k$ is essentially the smallest admissible $C$ since for any smooth $X$ and any admissible $C$ there exists a Zariski covering $\{U_i\sr X\}$ with $U_i\in C$. Unfortunately, $Sm/k$ is not $f$-admissible. There are two reasons for this. One is that the quotients may not exist for actions on smooth schemes which are not quasi-projective. This is easily resolved by considering the category of smooth quasi-projective schemes which is also admissible and whose category of sheaves in any topology which is at least as strong as the Zariski one, is equivalent to the category of sheaves on the whole $Sm/k$. Another one is that a quotient of a smooth scheme with respect to a finite group action need not be smooth. This is the main reason why we have to consider non-smooth schemes and one of the key sources of technical complexity of the paper.

For an admissible $C$ let $C_+$ be the full subcategory of the category of pointed objects in $C$ which consists of objects pointed by a disjoint base point, and $Cor(C,R)$ the category of finite correspondences over $C$. In our computations we will have to consider the $(Nis,\af)$-homotopy categories  $H_{Nis,\af}(C)$, $H_{Nis,\af}(C_+)$ and $H_{Nis,\af}(Cor(C,R))$ of $C$, $C_+$ and $Cor(C,R)$ respectively. To have a uniform treatment of these three cases we use the formalism of radditive functors developed in \cite{SRFsub} and summarized in Appendix \ref{radsum}. A radditive functor on a category $D$ with finite coproducts is a presheaf of sets on this category which takes finite coproducts to products. The category $Rad(D)$ of radditive functors on a small category $D$ is complete and cocomplete-complete with a set of compact generators and the category of simplicial radditive functors $\Delta^{op}Rad(D)$ on any $D$ provides a convenient environment for the homotopy theory. In our context, the category $Rad(C_+)$ is equivalent to the category $Rad(C)_{\BB}$ of pointed objects in $Rad(C)$ and $Rad(Cor(C,R))$ is equivalent to the category of presheaves with transfers of $R$-modules on $C$ which allows us to treat the cases of non-pointed, pointed and  "transfer enriched" homotopy theories as particular cases of the homotopy theory of simplicial radditive functors.    

The categories $H_{Nis,\af}(C)$, $H_{Nis,\af}(C_+)$ and $H_{Nis,\af}(Cor(C,R))$ are obtained by the application of the standard  localization constructions  to the simplicial objects in $Rad(C)$, $Rad(C_+)$ and $Rad(Cor(C,R))$ respectively. The category $H_{Nis,\af}(C)$ comes out to be  equivalent to the non-pointed homotopy category of the site with interval $(C_{Nis},\af)$ and the category $H_{Nis,\af}(C_+)$ to the pointed homotopy category of $(C_{Nis},\af)$ (see \cite{MoVo}).  The category $H_{Nis,\af}(Cor(C,R))$ is the $(Nis,\af)$-homotopy category of the category of finite correspondences on $C$ with coefficients in $R$ i.e. an unstable analog of the category $DM_{-}^{eff}(C,R)$. 

For each $C$ and $R$ we get a pair of adjoint functors 
$$\LL\Lambda^l_R:H_{Nis,\af}(C_+)\sr H_{Nis,\af}(Cor(C,R))$$
$$\Lambda^r_R:H_{Nis,\af}(Cor(C,R))\sr H_{Nis,\af}(C_+)$$
(see Theorem \ref{2009th1}), which are analogous to the pair which consists of the forgetting functor and the free $R$-module functor connecting the homotopy categories of pointed simplicial sets and simplicial $R$-modules. 

An inclusion of admissible categories $i:C\sr D$ defines pairs of adjoint functors $(\LL i^{rad}, i_{rad})$, $(\LL i^{rad}_+,i_{rad,+})$ and $(\LL i^{rad}_{tr},i_{rad,tr})$ between the homotopy categories of each type for $C$ and $D$. Corollary \ref{2009.07.24.4} asserts that the left adjoint of each pair is a full embedding i.e. $H_{Nis,\af}(C)$ is a full subcategory in $H_{Nis,\af}(D)$, etc. 

The left adjoints defined by $i$ commute with  $\LL\Lambda^l_R$ and the right adjoints with $\Lambda^r_R$. In an important addition to this elementary observation we prove in Theorem \ref{ressing} that if $k$ admits resolution of singularities and $C\subset Sm/k$ then $\LL\Lambda^l_R$ also commutes with the right adjoints defined by $i$. This result allows us to prove later Theorem \ref{mth1} which is a key step on the way to Theorem \ref{mainst}. 

In the second part of the paper we investigate the functors on the homotopy categories which are defined by the functors of generalized symmetric products on schemes. Let $\Phi$ be a pair of the form $(G,\phi:G\sr S_n)$ where $G$ is a group and $\phi$ an embedding of $G$ to the symmetric group on $n$ elements. Such a pair will be called a permutation group. If $C$ is an $f$-admissible category then a permutation group defines two generalized symmetric product functors $\tilde{S}^{\Phi}:X_+\mapsto (X^n/\phi(G))_+$ and $S^{\Phi}:X_+\mapsto (X_+)^n/\phi(G)$ from $C_+$ to itself. 

Using results of \cite{Delnotessub} we construct extensions of these functors to functors $\LL\tilde{S}^{\Phi}$ and $\LL S^{\Phi}$ from $H_{Nis,\af}(C_+)$ to itself. In order to do it we use the following construction. For a small category with finite coproducts $D$ let $D^{\#}$ be the full subcategory of directed colimits of representable functors in $Rad(D)$ i.e the formal completion of $D$ with respect to directed (and therefore filtered) colimits. It is also known as the category of ind-objects over $D$ (see e.g. \cite[I 8.2.4, p.70]{SGA41}). The homotopy theory of simplicial radditive functors provides a construction of a canonical resolution functor $L_*:\Delta^{op}Rad(D)\sr \Delta^{op}D^{\#}$ and an assertion that the projection functor from $\Delta^{op}D^{\#}$ to any homotopy category $H(D,E)$ associated with $D$ by the "standard construction" is a localization. In particular, in order to extend a functor $F:D_1\sr D_2$ to a functor $H(D_1,E_1)\sr H(D_2,E_2)$ it is sufficient to show that the obvious extension of $F$ to a functor $\Delta^{op}D_1^{\#}\sr \Delta^{op}D_2^{\#}$ takes $E_1$-local equivalences to $E_2$-local equivalences. Note that this applies to all functors $F$ including the ones which do not commute with finite coproducts. In the context of $D=C_+$ and $F=\tilde{S}^{\Phi}$ or $F=S^{\Phi}$ the required property of $F$ with respect to $(Nis,\af)$-equivalences was proved in \cite{Delnotessub}. Moreover, it was shown there that for these functors the equivalence $\LL F(X)=X$ holds not only for $X\in \Delta^{op}(C_+)^{\#}$ but also for a more general class of simplicial ind-solid sheaves, which includes the Nisnevich quotient sheaves $X/U$ for open embeddings $U\subset X$. In particular, up to a $Nis$-equivalence one has
$$\LL \tilde{S}^n(X/(X-Z))=\tilde{S}^n(X/(X-Z))=S^nX/(S^nX-S^nZ)$$
where $S^n=S^{(S_n,Id)}$ is the ordinary symmetric $n$-th power. 

In the next section we consider the extensions of $\tilde{S}^{\Phi}$ first to the categories $Cor(C,R)$ and then to $H_{Nis,\af}(Cor(C,R))$. We prove that such extensions $\LL S^{\Phi}_{tr}$, compatible via functors $\LL\Lambda^l_R$ with $\LL\tilde{S}^{\Phi}$, exist provided that $C$ is $f$-admissible and the "characteristic exponent" $c(k)$ of $k$  is invertible in the coefficient ring $R$. Along with the existence of $\LL S^{\Phi}_{tr}$ we prove a number of results which allow us in Theorem \ref{tonextnew} of the following section to give a complete computation of $\LL S^l_{tr}(L_n)$ where $l\ne char(k)$ is a prime and $L_n=\LL\Lambda^r_{{\bf F}_l}(T^n)$ is the "homology" of the motivic $(2n,n)$-sphere with coefficients in the finite field ${\bf F}_l$.

In Section \ref{PTO} we prove several relatively simple results about split proper Tate objects i.e. objects of $H_{Nis,\af}(Cor(C,R))$ which are isomorphic to (possibly infinite) direct sums of objects of the form $\Sigma^iL_j$. Using these results together with the results of the preceding sections we prove that in case when the ring of coefficients is a field $F$ of "allowed" characteristic the subcategory $\overline{SPT}$ of proper split Tate objects is closed  under the formation of all ordinary symmetric powers $\LL S^n_{tr}$. In fact, the result proved in Theorem \ref{main2} is more precise and also provides a lower bound on the range of "weights" which may appear in the $n$-th symmetric power of a proper split Tate object of weight $\ge q$ with coefficients in a field of characteristic $l$. This bound on the weight of summands was first observed, in a particular case, by C.~Weibel in \cite{patching} and plays an important role in one of the key arguments of \cite{zlsub}. 

The next ingredient which goes into the proof of Theorem \ref{2009.09.12.th1}  and Theorem \ref{mainst} are Theorem \ref{susvoe} and Proposition \ref{2009.09.07.pr1} which, when taken together, form a motivic analog of the Dold-Thom Theorem for connected $CW$-complexes. The role of connectedness assumption is played in our context by "condition (D1)" of Definition \ref{2009.09.07.def1} which holds for a wide class of motivic spaces including the motivic Moore spaces constructed in the following section. The combination of \ref{susvoe} and \ref{2009.09.07.pr1} implies that for spaces $X$ satisfying condition (D1) and an $f$-admissible $C$ which is contained in the category of semi-normal schemes, there is a canonical isomorphism in $H_{Nis,\af}(C_+)$ of the form
\begin{eq}
\llabel{2009.09.20.eq6}
\Lambda^r_S\LL\Lambda^l_S(X)=\LL S^{\infty}[1/c(k)](X)
\end{eq}
where $S=\zz[1/c(k)]$ and $\LL S^{\infty}[1/c(k)]$ is the homotopy colimit of the sequence of morphisms from $\LL S^{\infty}$ to itself defined by the multiplication by $c(k)$ in the obvious abelian monoid structure of $S^{\infty}$.

In Section \ref{moore.section} we finally define, in the context of a general admissible $C$, the motivic cohomology functors and the motivic Eilenberg-MacLane spaces. We could have done it as soon as the categories $H_{Nis,\af}(Cor(C,R))$ and functors $\LL\Lambda^l_{R}$, $\Lambda^r_R$ were defined but decided that it is better to give an informal definition here and delay the formal one until Section \ref{moore.section}.

Let $C$ be an admissible category. Consider the motivic spheres $S^1_t=(\af-\{0\},1)$ and $S^q_t=(S^1_t)^{\wedge q}$ as radditive functors on $C_+$.  Set
$$l_{q,R}=\LL\Lambda^l_{R}(S^q_t).$$
It is the object of $H_{Nis,\af}(Rad(Cor(C,R)))$ which represents the "homology" of the sphere $S^q_t$.  We will write $l_q$ for $l_{q,\zz}$. One defines the (reduced) motivic cohomology of $X\in \Delta^{op}Rad(C_+)$ with coefficients in an abelian group $A$ by the formula: 
\begin{eq}
\llabel{descr0intro}
\tilde{H}^{p,q}_{un,C}(X,A)=\left\{\begin{array}{lll}
Hom_{H_{Nis,\af}(Cor(C,\zz))}(\LL\Lambda^l_{\zz}X,\,\,\Sigma^{p-q}(A\oo_{\LL}l_q))&\rm for&p\ge q\\\\
Hom_{H_{Nis,\af}(Cor(C,\zz))}(\Sigma^{p-q}\LL\Lambda^l_{\zz}X,\,\,A\oo_{\LL}l_q)&\rm for&p\le q
\end{array}
\right.
\end{eq}
Corollary \ref{2009.07.24.4} implies that for $X\in C$ and $C\subset D$ one has
$$H^{p,q}_{un,C}(X,A)=H^{p,q}_{un,D}(X,A)$$
i.e. that the groups defined by (\ref{descr0intro}) depend only on $X$ and not on the ambient category used in the definition. Because of this property we will often write $H^{*,*}_{un}$ instead of $H^{*,*}_{un,C}$. 

If $C\subset Sm/k$ and $k$ is a perfect field then Theorem \ref{2009femb} shows that 
\begin{eq}
\llabel{motcoh1intro}
\tilde{H}^{p,q}_{un,C}(X,A)=Hom_{DM_{-}^{eff}(C,\zz)}(\tilde{M}(X), A(q)[p])
\end{eq}
i.e. for smooth schemes over a perfect field the definition given above agrees with the standard one which goes back to \cite{collection}. 

We use the subscript $un$ to emphasize the fact that {\em the motivic cohomology of general schemes as defined by (\ref{descr0intro}) are not known to have suspension isomorphisms with respect to either of the two motivic suspensions.} In particular, the equality  (\ref{motcoh1intro}) is not known to hold for non-smooth schemes since the right hand side of this equality automatically satisfies the suspension isomorphism for the simplicial suspension. Among many possible extensions of motivic cohomology functors to non-smooth schemes the one defined by (\ref{descr0intro}) is, to the best of my knowledge, the most basic and fundamental one in the sense that all other definitions can be obtained from this one by the imposition of additional properties such as stability or descent  for special coverings.  

The motivic Eilenberg-MacLane space $K(A,p,q)_C$ is the object of $H_{Nis,\af}(C_+)$ which represents the motivic cohomology functor $\tilde{H}_{un}^{p,q}(-,A)$ on this category.  The standard adjunctions show that for $p\ge q$ one has
\begin{eq}
\llabel{2009.09.20.eq1}
K(A,p,q)_C=\Lambda_{\zz}^r(\Sigma^{p-q}(A\oo_{\LL}l_q))
\end{eq}
In Proposition \ref{2009.09.06.th1} we prove that if $A$ is a finitely generated abelian group and $p\ge q$, $p>0$, then there exist well behaved (see below) spaces $X(A,p,q)$ such that
\begin{eq}
\llabel{2009.09.20.eq2}
\LL \Lambda^l_{\zz}(X(A,p,q))=\Sigma^{p-q}(A\oo_{\LL}l_q)
\end{eq}
Combining (\ref{2009.09.20.eq1}) and (\ref{2009.09.20.eq2}) we obtain the following expression for $K(A,p,q)_C$ in case when $A$ is finitely generated, $p\ge q$ and $p>0$:
\begin{eq}
\llabel{2009.09.20.eq3}
K(A,p,q)_C=\Lambda_{\zz}^r\LL\Lambda_{\zz}^l(X(A,p,q)).
\end{eq}
This expression is the key to all the computations with the motivic Eilenberg-MacLane spaces which are done in this paper. For $p<q$ no similar description exists and because of that the theory of spaces $K(A,p,q)_C$ for $p<q$ is very different from the theory for $p\ge q$. At the moment we know very little about the $p<q$ case and what we do know suggests that the structure of the spaces $K(A,p,q)$ and in particular of $K(\qq,p,q)$ for $p<q$ and $q\ge 2$ is very complex (for a related discussion see \cite{Sasha}).  

However, even in the case $p\ge q$ the spaces $K(A,p,q)_C$ themselves are very hard to study. 
For example, the standard adjunctions show that for an inclusion of admissible subcategories $i:C\sr D$ one has
$$i_{rad,+}(K(A,p,q)_D)=K(A,p,q)_C$$
but the adjoint equality $\LL i^{rad}_+(K(A,p,q)_C)=K(A,p,q)_D$ is not known to hold, except in trivial cases, even under the resolution of singularities assumption, which creates one of the key technical complications in our theory.

Fortunately, for any commutative ring $R$, all the information about the motivic cohomology of these spaces with coefficients in $R$-modules, i.e.  about the motivic cohomological operations from $\tilde{H}_{un}^{p,q}(-,A)$ to the motivic cohomology with coefficients in $R$-modules, is encoded in the objects
$$M(A,p,q;R)_C=\LL\Lambda^l_R(K(A,p,q)_C)$$
which are much more tractable. Combining (\ref{2009.09.20.eq3}), (\ref{2009.09.20.eq6}) and some technical arguments which allow us to consider $S$-coefficients instead of $\zz$-coefficients and $S^{\infty}$ instead of $S^{\infty}[1/c(k)]$ we prove our next main theorem (Theorem \ref{2009.09.12.th1}) which asserts that for an $f$-admissible $C$, finitely generated $S$-module $A$, $p\ge q$, $p>0$ and an $S$-algebra $R$ there are natural isomorphisms
$$M(A,p,q;R)_C\cong \oplus_{n\ge 0} S^n_{tr}(\Sigma^{p-q}((A\oo_{\LL,S}R)\oo_{R}l_{q,R}))$$
Together with the results of the first (Theorem \ref{ressing}) and the second (Theorem \ref{main2}) parts of the paper, Theorem \ref{2009.09.12.th1} leads to Theorem \ref{mth1} and then to Corollary \ref{2009.09.15.cor2}. 

This  corollary and a few elementary results about proper split Tate objects are the only results from the earlier parts of the paper which are needed for the proof of our main "stable" theorem (Theorem \ref{mainst}) which asserts that if $k$ is a field of characteristic zero then one has
$$lim_{n}H^{*+2n,*+n}(K(\zz/l,2n,n)_{Sm/k},\zz/l)={\cal A}^{*,*}(k,\zz/l)$$
where ${\cal A}^{*,*}(k,\zz/l)$ is the motivic Steenrod algebra defined in \cite{Redpub}.

An analog of Theorem \ref{mainst} should hold with $Sm/k$ replaced by an $f$-admissible $C$ which is contained in semi-normal schemes. The proof of such an analog would not require Theorem \ref{ressing} and therefore it would apply to all perfect fields. Writing such a proof up would require two developments - to construct etale realization functors with properties similar to the topological realization constructed in Section \ref{top.real} and to extend the results of \cite{Redpub} to $T$-stable version of motivic cohomology of non-smooth schemes. In order to directly extend Theorem \ref{mainst} to all perfect fields following the strategy used in the present paper it would be sufficient to construct the etale analog of the topological realization and to prove
that for a smooth scheme $X$ with an action of a finite group $G$ and a subgroup $H\subset G$, the map
$$\LL\Lambda^l_R i_{rad}(X/H)\sr \LL\Lambda^l_R i_{rad}(X/G)$$
where $i:Sm/k\sr SN/k$ and $[G:H]^{-1}\in R$, is a split epimorphism. In the present approach this statement follows via Theorem \ref{ressing} from the resolution of singularities but may be there is a more direct way of proving it. 

It took me thirteen years to write this paper and numerous people provided useful comments on its intermediate versions. I would especially like to thank Chuck Weibel who has helped me a lot as a source of both advice and willpower. I am also very grateful to one of the anonymous referees who not only suggested a great number of small improvements but also pointed out two places where the original arguments had to be substantially corrected.

\newpage \section{Motivic homology and homotopy}
\llabel{section2}
\subsection{Main categories and functors}
\llabel{MCF}
Let $k$ be a field and $C$ an admissible subcategory in $Sch/k$.  Since $C$ has finite coproducts the constructions and results of Appendix \ref{radsum} apply to $C$. The radditive functors on $C$ can be identified with presheaves on the category of connected objects in $C$ or alternatively with sheaves in the topology on $C$ which is generated by coverings by connected components. In particular the class of projective equivalences in $\Delta^{op}Rad(C)$ is closed under coproducts. Fix a set $C_0$ of representatives of isomorphism classes of objects in $C$. Consider the following two sets of morphisms in $\Delta^{op}C$:
\begin{enumerate}
\item The set $G_{Nis}$ of generating Nisnevich equivalences is the set of morphisms of the form $K_Q\sr X$ where $Q$ is an upper distinguished square in $C_0$ i.e. a pull-back square of the form
\begin{eq}
\llabel{uppersq}
\begin{CD}
B@>>> Y\\
@VVV @VVpV\\
A  @>j>> X
\end{CD}
\end{eq}
where $p$ is etale, $j$ is an open embedding, $B=p^{-1}(A)$, $p:Y-B\sr X-A$ an isomorphism and $K_{Q}$ is the object of $\Delta^{op}C$ given by the push-out square 
$$
\begin{CD}
B\amalg B @>>> B\oo\Delta^1\\
@VVV @VVpV\\
A\amalg Y @>j>> K_Q
\end{CD}
$$
\item The set $G_{\af}$ of generating $\af$-equivalences is the set of morphisms of the form $X\times\af\sr X$ for $X$ in $C_0$.
\end{enumerate}
Define the $(Nis,\af)$-homotopy category of $C$ setting:
$$H_{Nis,\af}(C)=H(C,G_{Nis}\cup G_{\af})$$
and similarly for $H_{Nis}(C)$ and $H_{\af}(C)$.

By  \cite[]{HH2sub} there is a natural equivalence
$$H_{Nis,\af}(C)\sr H(C_{Nis},\af)$$
where on the right hand side we have the homotopy category of the site with interval $(C_{Nis},\af)$ as defined in \cite{MoVo}.

Let $i_0,i_1:Spec(k)\sr \af$ be the morphisms corresponding to the points $0$ and $1$ of $\af$ respectively. An elementary $\af$-homotopy between morphisms $f,g:X\sr Y$ in $\Delta^{op}Rad(C)$ is a morphism $h:X\times\af\sr Y$ such that $h\circ (Id\times i_0)=f$ and $h\circ (Id\times i_1)=g$. Two morphisms are called elementary $\af$-homotopic if there exists an elementary $\af$-homotopy between them. Two morphisms are called $\af$-homotopic if they are equivalent with respect to the equivalence relation generated by the relation of being elementary $\af$-homotopic. A morphism $f:X\sr Y$ is called an $\af$-homotopy equivalence if there exists a morphism $g:Y\sr X$ such that the compositions $f\circ g$ and $g\circ f$ are $\af$-homotopic to the corresponding identity morphisms. 
\begin{proposition}
\llabel{heDelta1}
The $\af$-homotopy equivalences belong to $cl_l(G_{\af})$.
\end{proposition}
\begin{proof}
By \cite[2007satr]{SRFsub} it is enough to show $\af$-homotopy equivalences become isomorphisms after localization with respect to $cl_l(G_{\af})$. 
For any $X$ the morphisms $Id_X\times i_0$ and $Id_X\times i_1$ become equal in the localized category since they are both sections of the isomorphism $X\times\af\sr X$. Therefore, after the localization any two $\af$-homotopic morphisms become equal. 
We conclude that an $\af$-homotopy equivalence $f$ becomes an isomorphism since there exists $g$ such that both compositions $f\circ g$ and $g\circ f$ are equal to the corresponding identities.
\end{proof}
\begin{proposition}
\llabel{wlocasDelta}
A morphism $f$ in $\Delta^{op}Rad(C)$ belongs to $cl_l(G_{Nis})$ if and only if it is a local equivalence in the Nisnevich topology as a morphism of presheaves.
\end{proposition}
\begin{proof}
It follows easily from the Nisnevich variant of the Brown-Gersten theorem \cite[Lemma 3.1.18, p.101]{MoVo}. See also \cite{HH2sub}.
\end{proof}

\begin{lemma}
\llabel{clfree}
The classes $cl_l(G_{Nis})$, $cl_l(G_{\af})$, $cl_l(G_{Nis}\cup G_{\af})$ are closed under coproducts and finite direct products.
\end{lemma}
\begin{proof}
The case of coproducts follows immediately from the fact that projective equivalences in $\Delta^{op}Rad(C)$ are closed under coproducts. In the case of products let us consider for example the class  $cl_l(G_{Nis}\cup G_{\af})$. By Theorem \ref{maindelta} we have
$$cl_l(G_{Nis}\cup G_{\af})=cl_{\bdl}((G_{Nis}\amalg Id_C)\cup (G_{\af}\amalg Id_{C})\cup W_{proj})$$
In particular $cl_l(G_{Nis}\cup G_{\af})$ is $\bdl$-closed. Since it is closed under compositions it is sufficient to show that for $f:X\sr Y$ in $cl_l(G_{Nis}\cup G_{\af})$ and $Z\in \Delta^{op}Rad(C)$ one has $f\times Id_Z\in cl_l(G_{Nis}\cup G_{\af})$. Since $f\times Id_Z=\Delta(g)$ where $g$ is a morphism of bisimplicial objects whose rows are of the form $f\times Id_{F}$ for some $F$ in $Rad(C)$ it is sufficient to show that $f\times Id_F\in cl_l(G_{Nis}\cup G_{\af})$ for $F\in Rad(C)$. Applying the standard representable resolution functor $L_*$ of Proposition \ref{2009ref0}  to $F$ and using the fact that projective equivalences are closed under direct products we reduce the problem to the case when $F\in C^{\#}$. Using the fact that our class is closed under filtered colimits we further reduce it to the case $F=X\in C$. The functor $(-)\times X$ clearly takes $\bdl$-closures to $\bdl$-closures. It remains to verify that for $f\in  (G_{Nis}\amalg Id_C)\cup (G_{\af}\amalg Id_{C})\cup W_{proj}$ one has $f\times Id_X\in cl_l(G_{Nis}\cup G_{\af})$ which is straightforward. 
\end{proof}

Let $C_+$ the category of disjointly pointed objects of $C$.  As customary we write $\vee$ for the coproduct in the pointed case. Since $C_+$ has finite coproducts the constructions and results of Appendix \ref{radsum} apply to $C_+$. The radditive functors on $C_+$ can be identified with pointed radditive functors on $C$. In particular the class of projective equivalences in $\Delta^{op}Rad(C_+)$ is closed under coproducts. 

For $X$ in $C$ we write $X_+$ for $X\amalg Spec(k)$ considered as an object of $C_+$. The functor $(-)_+$ is left adjoint to the forgetting functor and in particular commutes with colimits.  Define the $(Nis,\af)$-homotopy category of $C_+$ setting:
$$H_{Nis,\af}(C_+)=H(C_+,(G_{Nis}\cup G_{\af})_+)$$
and similarly for $H_{Nis}(C)$ and $H_{\af}(C)$.

By  \cite[]{HH2sub} there is a natural equivalence
$$H_{Nis,\af}(C_+)\sr H_{\BB}(C_{Nis},\af)$$
where on the right hand side we have the pointed homotopy category of the site with interval $(C_{Nis},\af)$ as defined in \cite{MoVo}.

\begin{lemma}
\llabel{clpoint}
The classes $cl_l((G_{Nis})_+)$, $cl_l((G_{\af})_+)$ and $cl_l((G_{Nis}\cup G_{\af})_+)$ in $\Delta^{op}Rad(C_+)$ are closed under coproducts, direct products and smash products,
\end{lemma}
\begin{proof}
The case of coproducts follows immediately from the fact that projective equivalences in $\Delta^{op}Rad(C_+)$ are closed under coproducts. 

Smash products on $Rad(C_+)$ are defined as the radditive extension of the functor $\wedge: C_+\times C_+\sr C_+$ which takes $(U_+,V_+)$ to $(U\times V)_+$. One verifies easily that under the equivalence between $Rad(C_+)$ and the category of pointed presheaves on connected objects of $C$ the smash product of radditive functors is given by the usual smash product of pointed presheaves (it is sufficient to verify this property for representable presheaves). This interpretation implies that smash products respect projective equivalences. Therefore by the same reasoning as in the proof of Lemma \ref{clfree} it remains to verify that for $f\in  (G_{Nis}\amalg Id_C)_+$ (resp. $f\in (G_{\af}\amalg Id_{C})_+$ and $Z\in C_+$ one has $f\wedge Id_Z\in cl_l((G_{Nis})_+)$ (resp. $f\wedge Id_Z\in cl_l((G_{\af})_+)$) which is straightforward.
\end{proof}

\begin{proposition}
\llabel{reflect3}
The functor $(-)_+:C\sr C_+$ satisfies the conditions of Theorem \ref{2009ref4} and Corollary \ref{2009ref5}(2) with respect to the pair of classes of morphisms $(G_{Nis},(G_{Nis})_+)$, $(G_{\af},(G_{\af})_+)$ and $(G_{Nis}\cup G_{\af},(G_{Nis}\cup G_{\af})_+)$. Therefore $(-)_+$ defines adjoint pairs of functors between the corresponding homotopy categories and in each case the right adjoint functor reflects isomorphisms.
\end{proposition}
\begin{proof}
Straightforward verification using the coproduct parts of Lemmas \ref{clfree} and \ref{clpoint}.
\end{proof}
According to our general convention we should denote the functors between the homotopy categories which correspond to $(-)_+$ by $\LL (-)^{rad}_+$ and $\RR (-)_{+,rad}$. Using the interpretation of radditive functors on $C_+$ as pointed radditive functors on $C$ it is easy to see that the functor $(-)_+^{rad}$ takes a non-pointed radditive functor $F$ to the radditive functor $F\amalg Spec(k)$ pointed by the canonical morphism $pt\sr F\amalg Spec(k)$. In particular, it respects all projective equivalences and therefore can be used instead of $\LL (-)_{+}^{rad}$. It also shows that we may use the notation $(-)_+$ for $(-)_+^{rad}$ without any danger of confusion. 

Similarly, $(-)_{+,rad}$ in this interpretation is the forgetting functor from pointed radditive functors to radditive functors and we will denote it by $\phi$.

Let $Cor(C,R)$ be the category of finite correspondences with coefficients in a commutative ring $R$ between objects of $C$. This category was described in detail in \cite{cancellationsub} (see also \cite{Deglise1}). To distinguish objects of $C$ from the corresponding objects of $Cor(C,R)$ we let
$$[-]_R:C\sr Cor(C,R)$$
denote the functor which is the identity on objects and which takes morphisms to their graphs. The category $Cor(C,R)$ is additive and in particular has finite coproducts and $[-]$ commutes with finite coproducts. The radditive functors on $Cor(C,R)$ can be identified with $R$-linear functors from $Cor(C,R)$ to the category of $R$-modules  i.e. presheaves with transfers with coefficients in $R$ on $C$. In particular the class of projective equivalences in $\Delta^{op}Rad(Cor(C,R))$ is closed under coproducts.

Define $(Nis,\af)$-homotopy category of $Cor(C,R)$ setting:
$$H_{Nis,\af}(Cor(C,R))=H(Cor(C,R),[G_{Nis}\cup G_{\af}]_R)$$
When $C=Sm/k$ this category is the full subcategory of the category $DM^{eff}_{-}(k,R)$ which consists of complexes of sheaves with transfers with homotopy invariant cohomology sheaves such that $\uu{H}^i=0$ for $i>0$. See Theorem \ref{2009femb} below.

Category $Cor(C,R)$ carries a natural tensor structure which is given by the product of schemes on objects and by external products of relative cycles on morphisms. This structure is a functor $Cor(C,R)\times Cor(C,R)\sr Cor(C,R)$. Its radditive extension is a functor
$$Rad(Cor(C,R))\times Rad(Cor(C,R))=Rad(Cor(C,R)\times Cor(C,R))\sr Rad(Cor(C,R))$$
which we call the tensor product on $Rad(Cor(C,R))$ and denote by $\oo$. Alternatively, one can define $\oo$ on $Rad(Cor(C,R))$ by first extending $\oo$ to $Cor(C,R)^{\#}$ and then defining $F\oo G$ as $\pi_0(L_*(F)\oo L_*(G))$. 
\begin{lemma}
\llabel{cltr}
The classes $cl_l([G_{Nis}]_R)$, $cl_l([G_{\af}]_R)$ and $cl_l([G_{Nis}\cup G_{\af}]_R)$ are closed under direct sums and their intersections with $\Delta^{op}Cor(C,R)^{\#}$ are closed under tensor products.
\end{lemma}
\begin{proof}
The case of coproducts i.e. direct sums follows immediately from the fact that projective equivalences in $\Delta^{op}Rad(Cor(C,R))$ are closed under direct sums. The proof of the tensor product part is parallel to the proofs of Lemmas \ref{clfree} and \ref{clpoint}. 
\end{proof}

Since the category $Cor(C,R)$ is pointed the functor $[-]$ factors through the functor $(-)_+$. We denote the corresponding functor $C_+\sr Cor(C,R)$ by $\Lambda_R$. By definition $\Lambda_R(X_+)=[X]_R$. This functor commutes with finite coproducts and therefore defines a pair of adjoint functors 
$$\Lambda^l_R=(\Lambda_R)^{rad}:Rad(C_+)\sr Rad(Cor(C,R))$$
$$\Lambda^r_R=(\Lambda_R)_{rad}:Rad(Cor(C,R))\sr Rad(C_+)$$
where $\Lambda^r_R$ is the right adjoint. If we interpret $Rad(Cor(C,R))$ as the category of presheaves with transfers and $Rad(C_+)$ as the category of pointed radditive presheaves then $\Lambda^r_R$ is the "forgetting" functor which takes a presheaf with transfers to the same presheaf considered as a presheaf of pointed sets. By Proposition \ref{2009ref2} the functors $\Lambda^l_R$ and $\Lambda^r_R$ define a pair of adjoint functors $\LL\Lambda^l_R$ and $\RR\Lambda^r_R$ between the homotopy categories $H(C_+)$ and $H(Cor(C,R))$.

We will denote the classes of local equivalences in the context of $C$ by $W_{Nis}$, $W_{\af}$ and $W_{Nis,\af}$, in the context of $C_+$ by $W^{+}_{Nis}$, $W^{+}_{\af}$ and $W^{+}_{Nis,\af}$ and in the context of $Cor(C,R)$ by $W^{tr}_{Nis}$, $W^{tr}_{\af}$ and $W^{tr}_{Nis,\af}$.

\begin{theorem}
\llabel{2009th1}
For any $R$ the functor $\Lambda_R$ satisfies the conditions of Theorem \ref{2009ref4} and Corollary \ref{2009ref5}(2) with respect to the pairs of classes of morphisms $(G_{Nis},(G_{Nis})_+)$, $(G_{\af},(G_{\af})_+)$ and $(G_{Nis}\cup G_{\af},(G_{Nis}\cup G_{\af})_+)$. 

In particular, one has:
\begin{enumerate}
\item $\Lambda^r_R(W_{Nis}^{tr})\subset W^{+}_{Nis}$, $\Lambda^r_R(W_{\af}^{tr})\subset W^{+}_{\af}$ and $\Lambda^r_R(W_{Nis,\af}^{tr})\subset W^{+}_{Nis,\af}$
\item $\Lambda^l_R(W^{+}_{Nis}\cap \Delta^{op}C_+^{\#})\subset W_{Nis}^{tr}$, $\Lambda^l_R(W^{+}_{\af}\cap \Delta^{op}C_+^{\#})\subset W_{\af}^{tr}$ and $\Lambda^l_R(W^{+}_{Nis,\af}\cap \Delta^{op}C_+^{\#})\subset W_{Nis,\af}^{tr}$,
\end{enumerate}
the functors between the corresponding homotopy categories defined by $\LL\Lambda^l_R$ and $\RR\Lambda^r_R$ are adjoint and in each case $\RR\Lambda^r_R$ reflects isomorphisms. 
\end{theorem}
\begin{proof}
For simplicity of notation we will write $\Lambda$ instead of $\Lambda_R$. We have to prove four inclusions $\Lambda^l((G_{Nis})_+\vee Id_{C_+})\subset W^{tr}_{Nis}$, $\Lambda^l((G_{\af}])_{+}\vee Id_{C_+})\subset W^{tr}_{\af}$, $\Lambda^r([G_{Nis}]\oplus Id_{Cor(C,R)})\subset W^{+}_{Nis}$ and $\Lambda^r([G_{\af}]\oplus Id_{Cor(C,R)})\subset W^{+}_{\af}$. The functor $\Lambda^l$ takes $\vee$ to the direct sum $\oplus$. Therefore the first two inclusions follow from the definitions and the direct sums part of Lemma \ref{cltr}.

The functor $\Lambda^r$ takes $\oplus$ to direct product. Therefore the direct product part of Lemma \ref{clpoint}  implies that in order to prove the third and the fourth inclusion it is sufficient to show that $\Lambda^r([G_{Nis}])\subset W^{+}_{Nis}$ and $\Lambda^r([G_{\af}])\subset W^{+}_{\af}$. The fact that the functor $\phi$ which forgets the distinguished point reflects equivalences of all the considered types further implies that it is sufficient to show that  $\phi\Lambda^r([G_{Nis}])\subset W_{Nis}$ and $\phi\Lambda^r([G_{\af}])\subset W_{\af}$.

Propositions \ref{wlocasDelta} implies that in order to prove that  $\phi\Lambda^r([G_{Nis}])\subset W_{Nis}$ it is sufficient to show that for any upper distinguished square $Q$ of the form (\ref{uppersq}) the morphism $[K_Q]\sr [X]$ is a local equivalence on $C_{Nis}$ as a morphism of presheaves of sets or, equivalently, as a morphism of presheaves of abelian groups. It is further equivalent to the condition that the morphism of associated sheaves of abelian groups is a local equivalence on $C_{Nis}$. Consider the functor
$$\gamma:Rad(Cor(C,R))\sr ShvAb(C_{Nis})$$
which is the composition of the forgetting functor from $Rad(Cor(C,R))$ to presheaves of abelian groups on $C$ with the associated sheaf functor. Clearly, $\gamma$ respects finite coproducts and therefore it commutes with the $K_Q$ construction. Hence the morphism we are interested in can be written as
$$\gamma([K_Q]\sr [X])=(K_{\gamma([Q])}\sr \gamma([X]))$$
Consider a square $S$ of pre-sheaves of abelian groups on $C_{Nis}$ of the form
$$
\begin{CD}
S1 @>>> S2\\
@VVV @VVV\\
S3 @>>> S4.
\end{CD}
$$
Using the fact that $C_{Nis}$ has enough points one verifies easily that the associated morphism $K_S\sr S4$ is a local equivalence if and only if the sequence 
$$0\sr a_{Nis}S1\sr a_{Nis}S2\oplus a_{Nis}S3\sr a_{Nis}S4\sr 0$$
of associated Nisnevich sheaves is exact. We conclude that the morphism $K_{\gamma([Q])}\sr \gamma([X])$ is a local equivalence by \cite[Proposition 4.3.9]{SusVoe2new}.

To prove that $\phi\Lambda^r([G_{\af}])\subset W_{\af}$ it is sufficient by Lemma \ref{heDelta1}, to show that for $X\in C$ the map of presheaves of sets $p:[X\times\af]\sr [X]$ is an $\af$-homotopy equivalence. The inverse equivalence is given by $i:[X]\sr [X\times\af]$ corresponding to the point $0$ of $\af$. The $\af$-homotopy
$$[X\times\af]\times\af\sr [X\times\af]$$
between the identity and $p\circ i$  is provided by the composition
$$[X\times\af]\times\af\sr [X\times\af\times\af]\sr [X\times\af]$$
where the first map is a particular case of a general map of the form
$$[X]\times Y\sr [X\times Y]$$
and the second one is obtained from the multiplication map $\af\times\af\sr \af$.
\end{proof}

Since by definition $\RR\Lambda^r_R$ is the natural descent to the homotopy categories of the simplicial extension of the functor $\Lambda^r_R$ between the radditive functors we will often write $\Lambda^r_R$ instead of $\RR\Lambda^r_R$. 

\begin{remark}\rm{}\rm
The functor $\Lambda^l_R$ does not preserve projective equivalences between all objects of $\Delta^{op}Rad(C_+)$.  Consider for example the morphism $p:Spec(L)\sr Spec(k)$ where $L/k$ is a Galois extension with the Galois group $G$ and let $\check{C}(p)$ be the corresponding Cech simplicial object. Then the morphism $p':\check{C}(p)\sr Im(p)$ where $Im(p)$ is the image of $p$ in $Rad(C)$, is a projective equivalence and so is $p'_+$. On the other hand the sections of $\Lambda^l_R((\check{C}(p))_+)$ over $Spec(k)$ form a simplicial abelian group which computes homology of $G$ with coefficients in $R$ and therefore in general it is not equivalent to $\Lambda^l_R(Im(p))$ which is a single object in dimension zero. 
\end{remark}\rm{}

The standard cosimplicial object $\Delta^{\BB}_{\af}$ in $Sch/k$ (see \cite[p. 88]{MoVo} or \cite[p. 16]{MVW}) lies in any admissible subcategory $C$. For a radditive functor $F$ on $C$ we let $C_*(F)$ denote the simplicial radditive functor with terms 
$$C_n(F):U\mapsto F(U\times \Delta^n_{\af})$$
Similarly for $F$ in $Rad(C_+)$ we set
$$C^+_n(F):U_+\mapsto F((U\times \Delta^n_{\af})_+)$$
and for $F$ in $Rad(Cor(C,R))$ we set
$$C_n^{tr}(F):[U]\sr F([U\times \Delta^n_{\af}])$$
(in \cite{MoVo} the functor $C_*$ was denoted by $Sing_*$).  

For $F\in \Delta^{op}Rad(C)$ we may consider $C_*(F)$ as a bisimplicial object. Then the diagonal $\Delta C_*(F)$ is defined and belongs to $\Delta^{op}Rad(C)$. Similarly we get the functors $\Delta C_*^+$ and $\Delta C_*^{tr}$ on $\Delta^{op}Rad(C_+)$ and $\Delta^{op}Rad(C,R)$. The projection $\Delta^{\BB}_{\af}\sr Spec(k)$ defines natural transformations $Id_{\Delta^{op}Rad(C)}\sr \Delta C_*$, $Id_{\Delta^{op} Rad(C_+)}\sr \Delta C^+_*$ and $Id_{\Delta^{op} Rad(Cor(C,R))}\sr \Delta C^{tr}_*$.

\begin{proposition}
\llabel{2009.07.24.1}
For any $F\in \Delta^{op}Rad(C)$ the object $\Delta C_*(F)$ is $G_{\af}$-local and the morphism $F\sr \Delta C_*(F)$ belongs to $cl_l(G_{\af})$. Similarly for any  
$F\in \Delta^{op}Rad(C_+)$ the object $\Delta C^+_*(F)$ is $(G_{\af})_+$-local and the morphism $F\sr \Delta C^+_*(F)$ belongs to $cl_l((G_{\af})_+)$ and for any   $F\in \Delta^{op}Rad(Cor(C,R))$ the object $\Delta C^{tr}_*(F)$ is $[G_{\af}]$-local and the morphism $F\sr \Delta C^{tr}_*(F)$ belongs to $cl_l([G_{\af}])$.
\end{proposition}
\begin{proof}
The projection $p:\Delta^{\BB}_{\af}\times\af\sr \Delta^{\BB}_{\af}$ and the embedding $i:\Delta^{\BB}_{\af}\sr \Delta^{\BB}_{\af}\times\af$ which corresponds to the point $0$ of $\af$, are mutually inverse cosimplicial homotopy equivalences. Indeed, let $\phi$ be morphism $\af\sr \af$ given by $x\mapsto 0$. To construct a cosimplicial homotopy between $Id_{\Delta^{\BB}}\times Id_{\af}$ and $Id_{\Delta^{\BB}}\times \phi$ observe first that the multiplication morphism $(x,y)\mapsto xy$ defines an $\af$-homotopy between $\phi$ and $Id$. Therefore, in order to construct a required cosimplicial homotopy it is sufficient to construct a cosimplicial homotopy between the morphisms $Id_{\Delta^{\BB}}\times i_0$ and $Id_{\Delta^{\BB}}\times i_1$ where $i_0,i_1:Spec(k)\sr \af$ are the morphisms corresponding to the points $0$ and $1$ of $\af$. Such a homotopy is given by the morphisms $\theta_i$ defined in \cite[Def. 2.17, p.17]{MVW}.

Applying any term-wise functor to $p$ and $i$  we obtain a pair of mutually inverse simplicial homotopy equivalences. In particular, for any $F\in Rad(C)$ and $U\in C$ the map
$$C_*(F)(U)\sr C_*(F)(U\times\af)$$
defined by the projection is a homotopy equivalence of simplicial sets and the same applies to $C_*^+$ and $C_*^{tr}$. Since the class of weak equivalences of simplicial sets is $\bdl$-closed 
we conclude that the first half of the proposition holds. 

For any $F\in Rad(C)$ and any $n\ge 0$ the morphism $F\sr C_n(F)$ is easily seen to be an $\af$-homotopy equivalence (see \cite[Lemma 3.7]{MoVo} or \cite[Lemma 3.2.2]{H3new}). Therefore, for any $F\in \Delta^{op}Rad(C)$ the morphism $F\sr \Delta C_*(F)$ is in $cl_{\bdl}(G_{\af})$ and, by Theorem  \ref{maindelta} we conclude that this morphism is in $cl_l(G_{\af})$. Similar argument applies in the two other contexts.
\end{proof}

\subsection{The category $H_{Nis,\af}(Cor(Sm/k))$ and the category $DM_{-}^{eff}(k)$}\llabel{conDM}
Let $C$ be an admissible subcategory in $Sch/k$ and $R$ a commutative ring. A (Nisnevich) sheaf with transfers of $R$-modules on $C$ is an object $F$ of $Rad(Cor(C,R))$ such that $\Lambda_R^r(F)$ is a (pointed) sheaf  on $C_{Nis}$. We let $Shv_{Nis}(Cor(C,R))$ denote the full subcategory of  Nisnevich sheaves in $Rad(Cor(C,R))$. The proof of the following result in the context of any admissible subcategory $C$ in $Sch/k$ is strictly parallel to their proof in in the context of $C=Sm/k$ (see \cite[Th. 3.1.4, Lemma 3.1.2]{H3new}, \cite{Deglise1}). 
\begin{proposition}
\llabel{asstr}
The inclusion functor $\iota^{tr}_{Nis}:Shv_{Nis}(Cor(C,R))\sr Rad(Cor(C,R))$ has a left adjoint $a_{Nis}^{tr}$ such that $\Lambda^r\circ a_{Nis}^{tr}=a_{Nis}\circ \Lambda^r$ where $a_{Nis}$ is the usual associated sheaf functor.  
\end{proposition}
\begin{cor}
\llabel{abeliantr}
The category $Shv_{Nis}(Cor(C,R))$ is abelian.
\end{cor}
\begin{lemma}
\llabel{lshv}
For any $X\in C$ the representable radditive functor $[X]\in Rad(Cor(C,R))$ is a Nisnevich sheaf with transfers.
\end{lemma}

Let $D_{-}(Shv_{Nis}(Cor(C,R)))$ be the derived category of complexes bounded from the above over $Shv_{Nis}(Cor(C,R))$. Note that a morphism in $Cmpl_{-}(Shv_{Nis}(Cor(C,R)))$ is a quasi-isomorphism if and only if it is a quasi-isomorphism as a morphism of complexes of Nisnevich sheaves of abelian groups. 

Let  $N:\Delta^{op}Rad(Cor(C,R))\sr Cmpl_{-}(Rad(Cor(C,R)))$ be the normalization functor. A morphism $f:X\sr Y$ in $\Delta^{op}Rad(Cor(C,R))$ is a projective equivalence if and only if $N(f)$ is a quasi-isomorphism and a local equivalence in the Nisnevich topology if and only if $a_{Nis}^{tr}N(f)$ is a quasi-isomorphism. Together with Proposition \ref{wlocasDelta} this implies that $N$ defines a functor
$$N_{Nis}:H_{Nis}(Cor(C,R))\sr  D_{-}(Shv_{Nis}(Cor(C,R)))$$
Let 
$$K:Cmpl_{-}(Rad(Cor(C,R)))\sr \Delta^{op}Rad(Cor(C,R))$$
be the right adjoint to $N$. As for any abelian category, it takes a complex to the simplicial object corresponding to the canonical truncation of this complex at level zero. Therefore the same reasoning as above implies that if $f$ is such that $a_{Nis}^{tr}(f)$ is a quasi-isomorphism then $K(f)$ is a local equivalence in the Nisnevich topology and by  Proposition \ref{wlocasDelta} an element of $W_{Nis}^{tr}$. We conclude that $K$ defines a functor
$$K_{Nis}:D_{-}(Shv_{Nis}(Cor(C,R)))\sr H_{Nis}(Cor(C,R))$$
which is right adjoint to $N_{Nis}$. In addition, since $N$ is a full embedding the adjunction $Id\sr 
KN$ is an isomorphism and since $N_{Nis}$ and $K_{Nis}$ are direct descends of $N$ and $K$  the adjunction $Id\sr K_{Nis}N_{Nis}$ is an isomorphism. We proved the following result.
\begin{proposition}
\llabel{Nfullemb}
The adjoint functors 
$$N:\Delta^{op}Rad(Cor(C,R))\sr Cmpl_{-}(Rad(Cor(C,R)))$$
$$K:Cmpl_{-}(Rad(Cor(C,R)))\sr \Delta^{op}Rad(Cor(C,R))$$
descend to adjoint functors
$$N_{Nis}:H_{Nis}(Cor(C,R))\sr  D_{-}(Shv_{Nis}(Cor(C,R)))$$
$$K_{Nis}:D_{-}(Shv_{Nis}(Cor(C,R)))\sr H_{Nis}(Cor(C,R))$$
The functors $N$ and $N_{Nis}$ are full embeddings and the functors $K$ and $K_{Nis}$ are localizations.
\end{proposition}

For a class $E$ of morphisms in a triangulated category let $cl_{vl}(E)$ denote the (left) Verdier closure of $E$ i.e. the class of morphisms whose cones belong to the localizing subcategory generated by cones of morphisms from $E$.

Define $DM_{-}^{eff}(C,R)$ as the localization
$$DM_{-}^{eff}(C,R)=D_{-}(Shv_{Nis}(Cor(C,R)))[cl_{vl}([G_{\af}])^{-1}]$$
In the case when $C=Sm/k$ our definition agrees with the standard one by \cite[Prop. 3.2.3]{H3new} i.e.
$$DM_{-}^{eff}(k,R)=DM_{-}^{eff}(Sm/k,R)$$
Let $W^{DM}$ be the class of morphisms in $Cmpl_{-}(Rad(Cor(C,R)))$ which become isomorphisms in $DM_{-}^{eff}(C,R)$. The standard properties of Verdier localization imply that a morphism $f$ in $D_{-}(Shv_{Nis}(Cor(C,R)))$ becomes an isomorphism in $DM_{-}^{eff}(C,R)$ if and only if it belongs to $cl_{vl}([G_{\af}])$. Therefore $W^{DM}$ coincides with the class of morphisms whose image in $D_{-}(Shv_{Nis}(Cor(C,R)))$ lies in $cl_{vl}([G_{\af}])$.
\begin{proposition}
\llabel{functorm}
One has 
$$N(W_{Nis,\af}^{tr})\subset W^{DM}$$
and therefore $N$ descends to a functor
$$N_{Nis,\af}:H_{Nis,\af}(Cor(C,R))\sr DM_{-}^{eff}(C,R)$$
\end{proposition}
\begin{proof}
By Proposition \ref{deltaver} we have
$$N_{proj}(W_{Nis,\af}^{tr})\subset cl_{vl}(N_{proj}([G_{Nis}])\cup N_{proj}([G_{\af}]))$$
Since triangulated functors map Verdier closures to Verdier closures it is sufficient to check that $N_{Nis}([G_{Nis}])$ and $N_{Nis}([G_{\af}])$ are contained in $cl_{vl}([G_{\af}])$ in $D_{-}(Shv_{Nis}(Cor(C,R)))$ which is obvious.
\end{proof}

\begin{theorem}
\llabel{2009femb}
Let $k$ be a perfect field and $C$ an admissible subcategory contained in $Sm/k$. Then the functor $N_{Nis,\af}$ is a full embedding.
\end{theorem}
\begin{proof}
Denote by $CC_*$ the functor from $Cmpl_{-}(Rad(Cor(R,C)))$ to itself of the form
$$CC_*(X)=Tot(N(C^{cmpl}_*(X)))$$
where $C_*^{cmpl}$ is the functor from complexes to simplicial complexes obtained by applying $C_*^{tr}$ to a complex $X$ term by term. The normalization of a simplicial complex in the simplicial direction is a bicomplex of which we take the total complex. This operation only involves finite direct sums since the bicomplex in question lies in the second and third quadrants with only finitely many rows lying in the third one. 

For an individual  radditive functor $F$ the cohomology presheaves of $NC^{tr}_*(F)$ are homotopy invariant by Proposition \ref{2009.07.24.1} or by \cite[Prop. 3.6]{H2new}. On the other hand it is easy to see that if $B$ is a bicomplex such that the cohomology presheaves of its rows or columns are $\af$-homotopy invariant then the same holds for $Tot(B)$. Therefore the cohomology presheaves of $CC_*(X)$ are homotopy invariant for any $X$. Similarly, from Proposition \ref{2009.07.24.1} and from \cite[Lemma 3.2.5]{H3new} we conclude that for an individual  radditive functor $F$ the morphism $F\sr NC^{tr}_*(F)$ is in $W^{DM}$ which easily implies that the same holds for any complex $X$ of such functors. 

This reasoning holds for any $C$ and $R$. If $C=Sm/k$ or equivalently if $C$ is contained in $Sm/k$ and $k$ is a perfect field then we know from the second part of \cite[Prop. 3.2.3]{H3new} that a morphism $f\in W^{DM}$ whose source and target are complexes with homotopy invariant cohomology presheaves is a quasi-isomorphism in the Nisnevich topology. Therefore $CC_*$ takes elements of $W^{DM}$ to Nisnevich quasi-isomorphisms and the functor $K_{Nis}\circ CC_*$ descends to a functor
$${\bf R}K_{Nis,\af}:DM_{-}^{eff}(C,R)\sr H_{Nis,\af}(Cor(C,R))$$
The natural transformation $Id\sr CC_*$ provides a construction of a pair of natural transformations $Id\sr {\bf R}K_{Nis,\af}\circ N_{Nis,\af}$ and $N_{Nis,\af}\circ {\bf R}K_{Nis,\af}\sr Id$ which form an adjunction. Up to this point the construction would go through with any other $\af$-localization functor instead of $CC_*$.  To prove that $N_{Nis,\af}$ is a full embedding we need to show that the first of these transformations is an isomorphism. For $X\in\Delta^{op}Rad(Cor(C,R))$ it is represented by the composition
\begin{eq}\llabel{2009.07.24.eq1}
X\sr K(N(X))\sr K(CC_*(N(X))))=K(Tot(N(C^{cmpl}_*(N(X)))))
\end{eq}
The first morphism is an isomorphism of the Dold-Thom correspondence. We obviously have
$$N(C^{cmpl}_*(N(X)))=N_{rows}N_{columns}C_*^{tr}(X)$$
and by Eilenberg-Zilber Theorem there is a natural homotopy equivalence of complexes of the form $Tot(N_{rows}N_{columns}B)\cong N(\Delta B)$ for any bisimplicial object $B$. Therefore
$$K(Tot(N(C^{cmpl}_*(N(X)))))=K(Tot(N_{rows}N_{columns}C_*^{tr}(X)))\cong K(N(\Delta C^{tr}_*(X)))= \Delta C^{tr}_*(X)$$
and the morphism (\ref{2009.07.24.eq1}) is isomorphic to the morphism $X\sr \Delta C_*^{tr}(X)$ which is a $[G_{\af}]$-local equivalence by Proposition \ref{2009.07.24.1}. Theorem is proved.
\end{proof}
\begin{remark}\rm{}\rm
I do not know whether or not the functor $N_{Nis,\af}$ is a full embedding for a general admissible $C$. The problem is that the main theorem of \cite{H2new} is only known for smooth schemes. On the other hand it should be possible to prove using general arguments that $N_{Nis,\af}$ becomes a full embedding after $H_{Nis,\af}$ is stabilized with respect to the simplicial suspension. From this point of view the main theorem of \cite{H2new} may be stated by saying that $H_{Nis,\af}(Cor(Sm/k))$ is $\Sigma_s$-stable. 
\end{remark}\rm{} 
The following proposition describes the behavior of the functors $N$ relative to the cofiber sequences in $H(Cor(C,R))$. For a general definition of a cofiber sequence see \cite[Def. 6.2.7, pp.156]{Hovey}. Since $Cor(C,R)$ is additive the coaction part of the cofiber sequence is determined by the boundary map (cf. \cite[Rm. 7.1.3, p.178]{Hovey}) and we will write the cofiber sequences as $X\sr Y\sr Z\sr \Sigma^1X$. Proposition \ref{2009ref6} shows that a sequence of this form in $H(Cor(C,R))$ or in any other of the homotopy categories of $Cor(C,R)$ is a cofiber sequence if and only if it is isomorphic to the image of the sequence defines by a term-wise coprojection sequence $X\sr Y\sr Z$ in $\Delta^{op}Cor(C,R)^{\#}$. Since such coprojection sequences are exact we get the following result.
\begin{proposition}
\llabel{2009.08.07.pr1}
With respect to the natural isomorphisms $N(\Sigma^1(Z))\cong N(Z)[1]$ the functor $N_{Nis,\af}$ maps cofiber sequences in $H_{Nis,\af}(Cor(C,R))$ to distinguished triangles in $DM_{-}^{eff}(C,R)$.
\end{proposition}
We also mention without proof the following result which can be easily deduced from the proof of Theorem \ref{2009femb}.
\begin{proposition}
\llabel{2009.08.07.pr2}
If $k$ is a perfect field and $C$ is contained in $Sm/k$ then a sequence of the form $X\sr Y\sr Z\sr \Sigma^1(Z)$ in $H_{Nis,\af}(Cor(C,R))$ is a cofiber sequence if and only if its image in $DM_{-}^{eff}(C,R)$ is a distinguished triangle.
\end{proposition}

\subsection{Change of the underlying category $C$}\llabel{changeofC}

For the purpose of the following discussion let us denote the classes $G_{Nis}$ and $G_{\af}$ in $\Delta^{op}Rad(C)$ for a given admissible subcategory $C$ by $G_{Nis}^C$ and $G_{\af}^C$.

Let $C$, $D$ be two admissible subcategories such that $i:C\subset D$. Since $i$ commutes with coproducts it defines a pair of adjoint functors $i^{rad}:Rad(C)\sr Rad(D)$ and $i_{rad}:Rad(D)\sr Rad(C)$. If we interpret radditive functors as sheaves in the topology defined by coverings by connected components then $i^{rad}$ and $i_{rad}$ become the usual inverse and direct image functors for the corresponding continuous map of sites. Note that this map of sites is not in general a morphism of sites i.e. $i^{rad}$ does not commute with limits. 
\begin{proposition}
\llabel{2009.07.24.3}
The inclusion functor $i$ satisfy the conditions of Theorem \ref{2009ref4} and Proposition \ref{2009ref5}(1) relative to the pairs of classes $(G_{Nis}^C,G_{Nis}^D)$ and $(G_{\af}^C,G_{\af}^D)$. 

The same holds for the functors $i_+:C_+\sr D_+$ and $i_{tr}:Cor(C,R)\sr Cor(D,R)$ relative to the pairs of classes $((G_{Nis}^C)_+,(G_{Nis}^D)_+)$,  $((G_{\af}^C)_+,(G_{\af}^D)_+)$ and $([G_{Nis}^C],[G_{Nis}^D])$, $([G_{\af}^C],[G_{\af}^D])$ respectively. 

Therefore we get pairs of adjoint functors $(i_{rad},\LL i^{rad})$, $(i_{+,rad},\LL i_+^{rad})$ and $(i_{tr,rad},\LL i_{tr}^{rad})$ between the projective, $Nis$- and $\af$-homotopy categories and the left adjoints in these pairs are full embeddings.  
\end{proposition}
\begin{proof}
The condition on the left adjoint functor is obvious in each of the cases. The condition on the right adjoint in the case of $C$ follows easily from Proposition \ref{wlocasDelta} and Proposition \ref{heDelta1}.  In the case of $C_+$ and $Cor(-,R)$ it follows from the same propositions and the fact that the  functor which forgets transfers reflect isomorphisms (Proposition \ref{reflect3} and Theorem \ref{2009th1}). 
\end{proof}
\begin{cor}
\llabel{2009.07.24.4}
The inclusion functor $i$ satisfy the conditions of Theorem \ref{2009ref4} and Proposition \ref{2009ref5}(1)  relative to the pair of classes $(G_{Nis}^C\cup G_{\af}^C,G_{Nis}^D\cup G_{\af}^D)$. The same holds for he functors $i_+:C_+\sr D_+$ and $i_{tr}:Cor(C,R)\sr Cor(D,R)$ relative to the pairs of classes $((G_{Nis}^C)_+\cup (G_{\af}^C)_+,(G_{Nis}^D)_+\cup (G_{\af}^D)_+)$ and $([G_{Nis}^C]\cup [G_{\af}^C],[G_{Nis}^D]\cup [G_{\af}^D])$ respectively.

Therefore we get pairs of adjoint functors $(i_{rad},\LL i^{rad})$, $(i_{+,rad},\LL i_+^{rad})$ and $(i_{tr,rad},\LL i_{tr}^{rad})$ between the $(Nis,\af)$-homotopy categories and the left adjoints in these pairs are full embeddings.  
\end{cor}

Let us analyze now how the adjoint pairs $(\LL i_+^{rad}, i_{+,rad})$ and  $({\bf L}i_{tr}^{rad}, i_{tr,rad})$  are related to each other. Consider the following diagram
\begin{eq}
\llabel{foursq}
\begin{CD}
H_{Nis,\af}(C_+) @>{\bf L}i_+^{rad}>> H_{Nis,\af}(D_+) @>i_{+,rad}>> H_{Nis,\af}(C_+)\\
@V{\bf L}\Lambda^l VV  @V{\bf L}\Lambda^l VV  @V{\bf L}\Lambda^l VV\\
H_{Nis,\af}(Cor(C,R)) @>{\bf L}i_{tr}^{rad}>> H_{Nis,\af}(Cor(D,R)) @>i_{tr,rad}>> H_{Nis,\af}(Cor(C,R))\\
@V\Lambda^r VV @V\Lambda^r VV @V\Lambda^r VV\\
H_{Nis,\af}(C_+) @>{\bf L}i_{+}^{rad}>> H_{Nis,\af}(D_+) @>i_{+,rad}>> H_{Nis,\af}(C_+)
\end{CD}
\end{eq}
By Corollary \ref{2009.07.24.4} the compositions of the horizontal arrows are canonically isomorphic to the corresponding identities.  The commutative square
$$
\begin{CD}
C_+ @>i_+>> D_+\\
@V\Lambda^l VV @VV\Lambda^l V\\
Cor(C,R) @>i_{tr}>> Cor(D,R)
\end{CD}
$$
shows that the lower right square of (\ref{foursq}) commutes. The upper left square is left adjoint to the lower right one and therefore it commutes  as well. The lower left square is unlikely to commute. We do not know whether the upper right one commutes in general but there is the following important partial result. Note that while the previous discussion holds without change in the context of projective, $Nis$- and $\af$-equivalences the theorem below is only known to be valid for the $(Nis,\af)$-equivalences.  
\begin{theorem}
\llabel{ressing}
Let $k$ be a field with resolution of singularities and $C\subset Sm/k$. Then for any admissible $D$ which contains $C$ the upper right square of (\ref{foursq}) commutes i.e. one has
$$\LL\Lambda^l i_{+,rad}=i_{tr,rad}\LL\Lambda^l.$$
\end{theorem}
\begin{proof}
Consider first the case $C=Sm/k$ and $D=Sch/k$ i.e. the commutativity of the square
\begin{eq}\llabel{hardsq1}
\begin{CD}
H_{Nis,\af}((Sch/k)_+) @>i_{+,rad}>> H_{Nis,\af}((Sm/k)_+)\\
@V{\bf L}\Lambda^l VV @V{\bf L}\Lambda^l VV\\
H_{Nis,\af}(Cor(Sch/k)) @>i_{tr,rad}>> H_{Nis,\af}(Cor(Sm/k))
\end{CD}
\end{eq}
We have a natural transformation
\begin{eq}
\llabel{new1}
{\bf L}\Lambda^l i_{+,rad}\sr i_{tr,rad}{\bf L}\Lambda^l
\end{eq}
arising from the adjunctions and the commutativity of the lower right square of (\ref{foursq}) and we need to show that it is an isomorphism. Consider the square
\begin{eq}
\llabel{2009.08.2.eq1}
\begin{CD}
{\bf L} \Lambda^l\, i_{+,rad}\, {\bf L}i_+^{rad}\, i_{+,rad} @>>> {\bf L}\Lambda^l\, i_{+,rad}\\
@VVV @VVV\\
i_{tr,rad}\,{\bf L} \Lambda^l\, {\bf L} i_+^{rad}\, i_{+,rad} @>>> i_{tr,rad}\,{\bf L}\Lambda^l
\end{CD}
\end{eq}
where the vertical arrows come from (\ref{new1}) and the horizontal ones from the adjunction 
$${\bf L}i_+^{rad} i_{+,rad}\sr Id.$$
We need to prove that the right vertical arrow is an isomorphism. We will do it by showing that the other three arrows are isomorphisms. 

By Corollary \ref{2009.07.24.4} the functors ${\bf L}i_+^{rad}$ and $\LL i_{tr}^{rad}$ are full embeddings and therefore the canonical morphisms
\begin{eq}
\llabel{new2}
i_{+,rad}\sr i_{+,rad}\, {\bf L}i_+^{rad}\, i_{+,rad}
\end{eq}
and
\begin{eq}
\llabel{new2tr}
i_{tr,rad}\sr i_{tr,rad}\, {\bf L}i_{tr}^{rad} i_{tr,rad}
\end{eq}
are isomorphisms. We conclude that the upper horizontal arrow in  (\ref{2009.08.2.eq1}) is an isomorphism and exchanging ${\bf L} \Lambda^l$ and ${\bf L} i_+^{rad}$ by commutativity of the upper left square of the main diagram that the left vertical arrow is an isomorphism.  It remains to show that under our assumptions the lower horizontal arrow is an isomorphism. It follows from Lemmas \ref{cdh0}-\ref{cdh3}.
\begin{lemma}
\llabel{cdh0}
In the diagram of projective homotopy categories similar to (\ref{foursq}) which is defined by the commutative square of functors
$$\begin{CD}
C @>i>> D\\
@V(-)_+VV @VV(-)_+V\\
C_+ @>i_+>> D_+
\end{CD}
$$
all four squares commute. The same holds for the diagrams of homotopy categories of all other types which we consider.
\end{lemma}
\begin{proof}
Straightforward using the identification of radditive functors on $C_+$ and $D_+$ with pointed radditive functors on $C$ and $D$ respectively.
\end{proof}
\begin{lemma}
\llabel{cdh1}
Let $k$ be a field with resolution of singularities, then the morphism
$${\bf L}i_+^{rad}i_{+,rad}(F)=i_+^{rad}\,\Delta L_*\, i_{+,rad}(F)\sr i_+^{rad}i_{+,rad}(F)\sr F$$
is a local equivalence in the cdh-topology for any $F$ in $\Delta^{op}Rad((Sch/k)_+)$.
\end{lemma}
\begin{proof}
It is clearly sufficient to show that $i_+^{rad}$ takes projective equivalences to local equivalences in the $cdh$-topology.  Since the forgetting functor from the pointed to non-pointed context reflects local equivalences in the cdh-topology and in view of Lemma \ref{cdh0} it is sufficient to prove that $i^{rad}$ takes projective equivalences in $\Delta^{op}Rad(Sm/k)$ to $cdh$-local equivalences in $\Delta^{op}Rad(Sch/k)$. 

Recall from \cite{HH2sub} that $Sm/k$ can be equipped with scdh-topology such that the natural functor $Sm/k\sr Sch/k$ defines a continuous map of sites 
$$\pi:(Sch/k)_{cdh}\sr (Sm/k)_{scdh}$$
and that when $k$ admits resolution of singularities $\pi$ is a morphism of sites and therefore  $\pi^*$ respects the local equivalences of simplicial sheaves. On the other hand the radditive functors can be interpreted as sheaves on the sites $(Sch/k)_{con}$ and $(Sm/k)_{con}$ whose topology is generated by coverings by connected components. 
Let $a_{scdh}$ be the associated sheaf functor from $Rad(Sm/k)$ to $Shv((Sm/k)_{scdh}$ and $a_{cdh}$ be the associated sheaf functor from $Rad(Sch/k)$ to $Shv((Sch/k)_{cdh}$. These functors may be considered as the inverse image functors for the obvious morphisms of sites $(Sm/k)_{scdh}\sr (Sm/k)_{con}$ and $(Sch/k)_{cdh}\sr (Sch/k)_{con}$. Therefore
$$a_{cdh}i^{rad}=\pi^*a_{scdh}$$
and since $\pi$ is a morphism of sites we conclude that $a_{cdh}i^{rad}$ takes projective equivalences i.e. local equivalences in the con-topology to local equivalences in the cdh-topology. 
\end{proof}
\begin{lemma}
\llabel{cdh2}
Let $f:X\sr Y$ be a morphism in $\Delta^{op}(Sch/k)_{+}^{\#}$ which is a local equivalence in the cdh-topology. Then $\Lambda^l(f)$ is a local equivalence in the cdh-topology.
\end{lemma}
\begin{proof}
By \cite{HH2sub} the class of local equivalences in the cdh-topology on $\Delta^{op}(Sch/k)_{+}^{\#}$ is $cl_{\bdl}((W_{Nis})_+\cup (W_{lcd})_+)$ where $W_{lcd}$ is defined in the same way as $W_{Nis}$ but with respect to the lower distinguished squares (i.e. abstract blow-up squares). The functor $\Lambda^l$ obviously takes $\bdl$-closures to $\bdl$-closures and therefore it is sufficient to verify that both for upper and for lower distinguished $Q$ the morphism $[K_Q]\sr [X]$ is an equivalence in the cdh-topology. For upper distinguished ones we know it from the proof of Theorem \ref{2009th1}. By the same argument the condition that $[K_Q]\sr [X]$ is a local equivalence in the cdh-topology for a lower distinguished square $Q$ is equivalent to the condition that for such a square the sequence of cdh-sheaves of abelian groups
$$0\sr a_{cdh}[B]\sr a_{cdh}[A]\oplus a_{cdh}[Y]\sr a_{cdh}[X]\sr 0$$
is exact in the cdh-topology. This is the statement of \cite[Prop. 4.3.3]{SusVoe2new}.
\end{proof}
\begin{lemma}
\llabel{cdh3}
Let $k$ be a field with resolution of singularities and $f:X\sr Y$ be a local equivalence in the cdh-topology in $\Delta^{op}Rad(Cor(Sm/k))$. Then the image of $f$ in $H_{Nis,\af}(Cor(Sm/k))$ is an isomorphism.
\end{lemma}
\begin{proof}
By Theorem \ref{2009femb} it is sufficient to check that the corresponding morphism is an isomorphism in $DM^{eff}_{-}$. It follows from \cite[Theorem 5.5(2)]{F-Vnew}.
\end{proof}
This finishes the proof of Theorem \ref{ressing} for $C=Sm/k$ and $D=Sch/k$. Changing $Sm/k$ to an admissible subcategory $C$ does not change, up to an equivalence, the categories of Nisnevich sheaves and therefore does not change the categories involved in the statement of the theorem. For an admissible subcategory $j:D\subset Sch/k$ which contains $C$  consider the diagram
\begin{eq}\llabel{2009.08.03.eq1}
\begin{CD}
H'(D_+) @>\LL j_{+}^{rad}>> H'((Sch/k)_+) @>j_{+,rad}>> H'(D_+) @>i_{+,rad}>> H'(C_+)\\
@V{\bf L}\Lambda^l VV  @V{\bf L}\Lambda^l VV @V{\bf L}\Lambda^l VV @V{\bf L}\Lambda^l VV\\
H'(D_+) @>\LL j_{tr}^{rad}>> H'(Cor(Sch/k)) @>j_{tr,rad}>> H'(Cor(D,R)) @>i_{tr,rad}>> H'(Cor(C,R))
\end{CD}
\end{eq}
where we write $H'$ instead of $H_{Nis,\af}$ to shorten the notation.  Let $S_1$, $S_2$ and $S_3$ be the squares of the diagram. The rectangle $S_2S_3$ commutes by the previous remark. The square $S_1$ commutes since it consists of two left adjoints. Therefore the ambient rectangle $S_1S_2S_3$ commutes. On the other hand it is equivalent to $S_3$ since by Corollary \ref{2009.07.24.4} the compositions $j_{+,rad}\LL j_{+}^{rad}$ and $j_{tr,rad}\LL j_{+}^{rad}$ are naturally isomorphic to the corresponding identities\footnote{I would like to thank one of the referees for this argument.}.
\end{proof}

\newpage \section{Symmetric powers}
\llabel{section3}
\subsection{Generalized symmetric powers on $H_{Nis,\af}(C_+)$}
In this section we assume that the underlying category $C$ is $f$-admissible i.e. it is admissible and categorical quotients exist for all finite group actions. For examples of such categories see Appendix \ref{apadmissible}. We will often use the following property of finite group quotients.
\begin{lemma}
\llabel{2009.08.3.l1}
Let $X, X'$ be schemes of finite type over a field $k$ and let $G\sr Aut(X)$, $G'\sr Aut(G')$ be actions on $X$ and $X'$ by finite groups $G$, $G'$ such that the categorical quotients $X/G$ and $X'/G'$ exist. Then $X/G\times X'/G'$ is the categorical quotient of $X\times X'$ by the product action  of $G\times G'$.
\end{lemma}
\begin{proof}
The quotients with respect to finite group actions commute with flat base changes. Therefore one has: 
$$(X/G)\times (X'/G')=(X\times (X'/G'))/G\times\{e\}=(X\times X'/\{e\}\times G')/G\times \{e\}=(X\times X')/G\times G'$$
\end{proof}
Let $\Phi=(G,\phi:G\sr S_n)$ be a permutation group i.e. a group together with an embedding into the symmetric group. Consider the functor $S^{\Phi}$ (resp. $\tilde{S}^{\Phi}$) from $C_+$ to itself of the form
$$X_+\mapsto (X_+)^n/G$$
(resp. of the form 
$X_+\mapsto (X_+)^{\wedge n}/G=(X^n/G)_+$
). 
Let $\Phi_i=(G_i,\phi_i:G_i\sr S_{n_i})$, $i=1,2$ be two permutation groups. Define their wreath product $\Phi_1*\Phi_2$ as follows. Let $\{n\}=\{1,\dots,n\}$. The direct power $G_1^{n_2}$ acts on $\{n_1\}\times \{n_2\}$ in the obvious way. Consider the action of $G_2$ on the same set which is the product of the action defined by $\phi_2$ on $\{n_2\}$ and the trivial action on $\{n_1\}$. Let $\phi:G\sr S_{n_1n_2}$ be the subgroup generated by the images of $G_1^{n_2}$ and $G_2$. We set 
$$\Phi_1*\Phi_2=(G,\phi:G\sr S_{n_1n_2})$$
(One can see that  $G$ is the semi-direct product of $G_1^{n_2}$ and $G_2$ with respect to the obvious action of the later on the former.) Lemma \ref{2009.08.3.l1} easily implies the following result. 
\begin{proposition}
\llabel{wreath}
For any $\Phi_1$ and $\Phi_2$ as above there are isomorphisms  of functors
\begin{eq}
\llabel{below4}
S^{\Phi_1*\Phi_2}=S^{\Phi_2} \circ S^{\Phi_1}
\end{eq}
and 
\begin{eq}
\llabel{below4tilde}
\tilde{S}^{\Phi_1*\Phi_2}=\tilde{S}^{\Phi_2} \circ \tilde{S}^{\Phi_1}
\end{eq}
\end{proposition}
Define $\Phi_1\times\Phi_2$ by the formula
$$\Phi_1\times\Phi_2=(G_1\times G_2, \phi:G_1\times G_2\sr S_{n_1+n_2})$$
where $\phi$ is the composition of $\phi_1\times\phi_2$ with the obvious embedding $S_{n_1}\times S_{n_2}\sr S_{n_1+n_2}$. We have another straightforward result.
\begin{proposition}
\llabel{product}
For any $\Phi_1$ and $\Phi_2$ as above there are isomorphisms of functors
\begin{eq}
\llabel{below5}
S^{\Phi_1\times\Phi_2}=S^{\Phi_1}\times S^{\Phi_2}
\end{eq}
and 
\begin{eq}
\llabel{below6}
\tilde{S}^{\Phi_1\times\Phi_2}=\tilde{S}^{\Phi_1}\wedge \tilde{S}^{\Phi_2}.
\end{eq}
\end{proposition}
Consider now the radditive extensions 
$$S^{\Phi,rad}:Rad(C_+)\sr Rad(C_+)$$
$$\tilde{S}^{\Phi,rad}:Rad(C_+)\sr Rad(C_+)$$
of $S^{\Phi}$ and $\tilde{S}^{\Phi}$. The isomorphisms of Propositions \ref{wreath} and \ref{product} extend immediately to these functors. To simplify the notation we will usually write $S^{\Phi}$ and $\tilde{S}^{\Phi}$ instead of $S^{\Phi,rad}$ and $\tilde{S}^{\Phi,rad}$. As all radditive extensions these functors commute with filtering colimits and reflexive coequalizers. The behavior of these functors with respect to finite coproducts corresponds to their behavior with respect to finite coproducts of schemes. For the case of ordinary symmetric products see Proposition \ref{minus1} below.  Proposition \ref{forcompl} gives an important example of how $\tilde{S}^{\Phi}$ behaves with respect to colimits of another type.

The main result which allows us to extend the symmetric powers functors to different homotopy categories is Theorem \ref{mainpreserve} below. The machinery needed for the proof of this theorem was developed in \cite{Delnotessub}. The proofs there are given in the case when $C=QP/k$ is the category of quasi-projective schemes over $k$ but they are applicable without any modification to any $f$-admissible category $C$.  Let $a_{+,Nis}$ be the functor of associated sheaf from $Rad(C_+)$ to $Shv_{\BB}(C_{Nis})$.
\begin{definition}
\llabel{2009.08.4.def1}
An object $F$ of $Rad(C_+)$ is called solid if there exists a filtration $pt=F_0\subset F_1\subset\dots\subset F_m=F$ such that for each $i=0,\dots,m-1$ there exist an open embedding $U_i\subset X_i$ in $C$ and an isomorphism 
$$a_{+,Nis}(F_{i+1}/F_i)\cong a_{+,Nis}(X_i/U_i)$$
An object is called ind-solid if it is a filtered colimit of solid objects. A simplicial object is called (ind-)solid if its terms are (ind-)solid. 
\end{definition}
Note  in particular that according to this definition  objects of $\Delta^{op}C_+$ are solid and objects of $\Delta^{op}C_+^{\#}$ are ind-solid.  Note also that the object of $Rad(C_+)$ represented by a pointed scheme $(X,x)$ is not ind-solid unless the distinguished point is disjoint. In view of \cite[Prop. 35 (Prop. 4.1.7 in the preprint version)]{Delnotessub} a radditive functor $F$ on $C_+$ is (ind-)solid if and only if $a_{+,Nis}(F)$ is (ind-)solid in the sense of \cite[Def. 7(Def. 4.1.5 in the preprint version)]{Delnotessub}. 
\begin{theorem}
\llabel{mainpreserve}
Let $\Phi=(G,\phi:G\sr S_n)$ be a permutation group. Then one has:
\begin{enumerate}
\item if $X$ is an (ind-)solid object of $\Delta^{op}Rad(C_+)$ the $\tilde{S}^{\Phi}(X)$ and $S^{\Phi}(X)$ are (ind-)solid objects of $\Delta^{op}Rad(C_+)$,
\item if  $f:X\sr Y$ is a Nis- (resp. an  $(Nis,\af)$-) equivalence between ind-solid objects in $\Delta^{op}Rad(C_+)$ then $\tilde{S}^{\Phi}(f)$ and $S^{\Phi}(f)$ are Nis- (resp. $(Nis,\af)$-) equivalences. In particular, for an ind-solid $F$ the obvious morphisms 
$$\tilde{S}^{\Phi}(L_*(F))\sr \tilde{S}^{\Phi}(F)$$
$${S}^{\Phi}(L_*(F))\sr {S}^{\Phi}(F)$$
are Nis-equivalences.
\end{enumerate}
\end{theorem}
\begin{proof} 
Let $C/G$ be the category of $G$-objects in $C$. Denote by $P^{\Phi}:C_+\sr (C/G)_+$ the functor which maps $X_+$ to $(X_+)^n$ with the permutational action of $G$ defined by $\phi$ and by $\tilde{P}^{\Phi}$ the functor which maps $X_+$ to $(X_+)^{\wedge n}=(X^n)_+$. 

Let further $Quot_G:(C/G)_+\sr C_+$ be the functor of the form $X_+\mapsto (X/G)_+$. Then $\tilde{S}^{\Phi}={Quot}_G\circ \tilde{P}^{\Phi}$ and $S^{\Phi}={Quot}_G\circ S^{\Phi}$. 

Define classes $(G_{\af})_+$ and $(G_{Nis})_+$ in $\Delta^{op}Rad((C/G)_+)$ in exactly the same way as we have defined these classes for $G=\{e\}$ in Section \ref{MCF} (see also \cite[p.389]{Delnotessub}). Also define the notion of solid and ind-solid object in $\Delta^{op}Rad((C/G)_+)$ in exactly the same way as we have done for $G=\{e\}$. Theorem \ref{mainpreserve} now follows from Propositions \ref{2009.08.05.pr0}-\ref{2009.08.04.pr4}.
\end{proof}
\begin{proposition}
\llabel{2009.08.05.pr0}
A morphism $f$ in $\Delta^{op}Rad((C/G)_+)$ belongs to $cl_l((G_{Nis})_+)$ (resp. to $cl_l((G_{Nis})_+\amalg (G_{\af})_+)$) if and only if $a_{+,Nis}(f)$ is a local (resp. $\af$-) equivalence in $\Delta^{op}(Shv_{\BB}((C/G)_{Nis}))$ in the sense of \cite{Delnotessub}.
\end{proposition}
\begin{proof}
Modulo the obvious analog of Proposition \ref{heDelta1} for $(C/G)_+$ this proposition is essentially equivalent to \cite[Th. 5 (Th. 3.6.1 in the preprint version)]{Delnotessub}.
\end{proof}
\begin{proposition}
\llabel{2009.08.04.pr1}
The radditive extensions of the functors $P^{\Phi}$ and $\tilde{P}^{\Phi}$ take (ind-)solid objects to (ind-)solid objects.
\end{proposition}
\begin{proof}
Let $I_{\Phi}$ be the object of $C/G$ which is the union of $n$ copies of $Spec(k)$ with the permutational action of $G$ defined by $\phi$. In \cite{Delnotessub} we denoted by $(-)^I$ the functor which takes a pointed sheaf $F$ on $(C/G)_{Nis}$ to the pointed sheaf $\uu{Hom}(I,F)$ and by $(-)^{\wedge I}$ a reduced version of this functor.  One verifies immediately that the functor $a_{+,Nis}(P^{\Phi})^{rad}$ is naturally isomorphic to $(-)^I\circ a_{+,Nis}$ and $a_{+,Nis}(\tilde{P}^{\Phi})^{rad}$ is naturally isomorphic to $(-)^{\wedge I}\circ a_{+,Nis}$. The proposition now follows from \cite[Th. 7 (Th. 5.2.3 in the preprint version)]{Delnotessub} in the case of  $P^{\Phi}$ and from \cite[Prop. 47 (Prop. 5.2.9 in the preprint version)]{Delnotessub} in the case of $\tilde{P}^{\Phi}$.
\end{proof}
\begin{proposition}
\llabel{2009.08.04.pr2}
Let $f:X\sr Y$ be a $Nis$- (resp. $(Nis,\af)$-) equivalence between ind-solid objects of $\Delta^{op}Rad(C_+)$. Then $(P^{\Phi})^{rad}(f)$ and $(\tilde{P}^{\Phi})^{rad}(f)$ are $Nis$- (resp. $(Nis,\af)$-) equivalences in $\Delta^{op}Rad((C/G)_+)$.
\end{proposition}
\begin{proof}
The statement for $\tilde{P}^{\Phi}$ follows from \cite[Prop. 48 (Prop. 5.2.11 in the preprint version)]{Delnotessub}. The statement for $P^{\Phi}$ follows by the same argument.
\end{proof}
\begin{proposition}
\llabel{2009.08.04.pr3}
The radditive extensions of the functor $Quot_G$ takes (ind-)solid objects to (ind-)solid objects.
\end{proposition}
\begin{proof}
The functor $a_{+,Nis}Quot^{rad}_G$ is naturally isomorphic to the functor $\eta_{+,\#}a_{+,Nis}$ where $\eta_{+,\#}$ is the pointed analog of the functor $\eta_{\#}:Shv((C/G)_{Nis})\sr Shv(C_{Nis})$ defined in \cite{Delnotessub}. It is shows in \cite{Delnotessub} that $\eta_{\#}$ is a left adjoint and in particular that it commutes with colimits. Therefore the same the same is true for $a_{+,Nis}Quot^{rad}_G$ and the statement of the proposition follows from the definition of (ind-)solid objects.
\end{proof}
\begin{proposition}
\llabel{2009.08.04.pr4}
Let $f:X\sr Y$ be a $Nis$- (resp. $(Nis,\af)$-) equivalence between ind-solid objects of $\Delta^{op}Rad((C/G)_+)$. Then $Quot_G^{rad}(f)$ is a $Nis$- (resp. $(Nis,\af)$-) equivalence in $\Delta^{op}Rad(C_+)$.
\end{proposition}
\begin{proof}
It follows from \cite[Prop. 45 (Prop. 5.1.4 in the preprint version)]{Delnotessub} through the identification $a_{+,Nis}Quot^{rad}_G=\eta_{+,\#}a_{+,Nis}$.
\end{proof}
As an immediate corollary from Theorem \ref{mainpreserve} we get the following result. 
\begin{cor}
\llabel{2009.08.05.cor1}
For any $\Phi=(G,\phi:G\sr S_n)$ there are unique (up to a canonical isomorphism) functors $\LL \tilde{S}^{\Phi}$ and $\LL S^{\Phi}$ on the $Nis$- and $(Nis,\af)$-homotopy categories which are determined by the conditions that the squares 
$$
\begin{CD}
\Delta^{op}C_+^{\#} @>\tilde{S}^{\Phi}>> \Delta^{op}C_+^{\#}\\
@VVV @VVV\\
H_{Nis,\af}(C_+) @>\LL \tilde{S}^{\Phi}>> H_{Nis,\af}(C_+)
\end{CD}
\,\,\,\,\,\,\,\,\,\,\,\,\,\,\,\,\,\,\,\,\,\,\,\,
\begin{CD}
\Delta^{op}C_+^{\#} @>{S}^{\Phi}>> \Delta^{op}C_+^{\#}\\
@VVV @VVV\\
H_{Nis,\af}(C_+) @>\LL {S}^{\Phi}>> H_{Nis,\af}(C_+)
\end{CD}
$$
and their analogs for the $Nis$-local categories, commute. In addition, for any ind-solid object $F$ in $\Delta^{op}Rad(C_+)$ one has $\LL \tilde{S}^{\Phi}(F)=\tilde{S}^{\Phi}(F)$ and $\LL S^{\Phi}(F)=S^{\Phi}(F)$ in both $Nis$- and $(\af,Nis)$-homotopy categories. 
\end{cor}
For a closed subset $Z$ of $X$ we let $\tilde{S}^{\Phi}(X_+)-\tilde{S}^{\Phi}(Z_+)$ denote the open subscheme in $S^{\Phi}(X_+)=(X^n/G)_+$ whose complement is the image of $Z^n\subset X^n$ under the canonical map from $X^n$.  
\begin{proposition}
\llabel{forcompl}
Let $Z$ be a closed subset of $X$. Consider $X/(X-Z)$ as a radditive functor on $C_+$ and let $a_{Nis}$ be the functor of associated Nisnevich sheaf. Then one has
$$a_{+,Nis}\tilde{S}^{\Phi}(X/(X-Z))=a_{+,Nis}(\tilde{S}^{\Phi}(X_+)/(\tilde{S}^{\Phi}(X_+)-\tilde{S}^{\Phi}(Z_+))).$$
\end{proposition}
\begin{proof}
As in the proof of Proposition \ref{2009.08.04.pr1}  observe that $a_{+,Nis}\tilde{S}^{\Phi}=\eta_{+,\#}(-)^{\wedge I}a_{+,Nis}$. Therefore the statement of our proposition follows from \cite[Example 8 (Example 5.2.8 in the preprint version)]{Delnotessub} and the fact that $\eta_{+,\#}$ commutes with colimits.
\end{proof}
\begin{example}\rm Let us consider the statement of Proposition \ref{forcompl} for $\tilde{S}^{\Phi}=\tilde{S}^2$ being the usual symmetric square. The pointed object $X/(X-Z)$ is the coequalizer of the reflexive pair $X_+\vee (X-Z)_+\dsr X_+$ where the first arrow is identity on $X_+$ and the inclusion on $(X-Z)_+$ and the second arrow is the identity on $X_+$ and the projection of $(X-Z)_+$ to the distinguished point. Applying to this equalizer the functor $S^2$ we see that $S^2(X/(X-Z))$ is the coequalizer of the pair
$$(S^2X)_+\vee (X\times (X-Z))_+\vee (S^2(X-Z))_+\dsr S^2X_+$$
and looking at the maps we conclude that it is defined by the push-out square of the form
$$\begin{CD}
X\times (X-Z) \amalg S^2(X-Z) @>>> S^2X\\
@VVV @VVV\\
Spec(k) @>>> S^2(X/(X-Z))
\end{CD}
$$
The statement of the proposition is valid in this case because the map $X\times (X-Z) \amalg S^2(X-Z) \sr S^2X-S^2Z$ is a Nisnevich covering. To see this consider the subset of the source at which the map is etale. It is of the form $(X\times(X-Z)-\Delta(X-Z))\amalg S^2(X-Z)$ and the restriction of the map to this subset is easily seen as arising from the upper distinguished square
$$
\begin{CD}
(X-Z)\times (X-Z)-\Delta(X-Z) @>>> X\times (X-Z)-\Delta(X-Z)\\
@VVV @VVV\\
S^2(X-Z) @>>> S^2X-S^2Z
\end{CD}
$$
However this is not a covering in the con-topology or even in the Zariski topology which make it  necessary to use the Nisnevich associated sheaf functor in Proposition \ref{forcompl}.
\end{example}
As was noted above the radditive functor on $C_+$ represented by a pointed scheme $(X,x)$ is not solid unless $x$ is a disjoint base  point. To extend the computations of derived  symmetric powers to such objects we need the following results.
\begin{proposition}
\llabel{2009.08.05.pr1}
Let $\Phi=(G,\phi)$ be a permutation group and $(X,x)$ an object of $\Delta^{op}Rad(C_+)$ represented by a pointed object of $\Delta^{op}C^{\#}$. Then the morphism
$$S^{\Phi}(\Delta L_*((X,x)))\sr S^{\Phi}((X,x))$$
is a projective equivalence and $S^{\Phi}((X,x))=(X^n/G,x^n)$.
\end{proposition}
\begin{proof}
It is clearly sufficient to verify the proposition for $(X,x)$ being a pointed object of $C$.  Consider the pair of adjoint functors
$$(-)_+:C^{\#}\sr C^{\#}_+$$
$$\phi:C^{\#}_+\sr C^{\#}$$
where $\phi$ is the functor which forgets the distinguished point. Let $F=(-)_+\circ \phi$ be the corresponding cocomplete-triple. Then $F((X,x))=X_+$. Since $\phi$ reflects projective equivalences and $\phi\,F=(-)\amalg pt$ respects projective equivalences, $F$ respects projective equivalences and  and in particular 
$$F L_*((X,x))\sr F((X,x))$$
is a projective equivalence between objects of $\Delta^{op}C_+^{\#}$. 

The natural isomorphisms $(-)^n_+=((-)_+)^{\wedge n}$ define functor isomorphisms
\begin{eq}
\llabel{punp}
F(S^{\Phi}(-))=\tilde{S}^{\Phi}(F(-))
\end{eq}
Therefore we have a commutative square
$$
\begin{CD}
F S^{\Phi} L_*((X,x)) @>>> \tilde{S}^{\Phi} F L_*((X,x))\\
@VVV @VVV\\
F S^{\Phi} ((X,x)) @>>> \tilde{S}^{\Phi} F((X,x))
\end{CD}
$$
in which the horizontal arrows are isomorphisms and the right vertical arrow is a projective equivalence by Proposition \ref{2009ref1}(2). We conclude that $FS^{\Phi} L_*((X,x)) \sr F S^{\Phi} ((X,x))$ is a projective equivalence. 

Since $Rad(C_+)=Rad(C)_{\bullet}$, one verifies easily that for any $X\in Rad(C_+)$ the map from the simplicial object $F_*(X)=(F^{\circ (i+1)}(X))_{i\ge 0}$ defined by the cocomplete-triple $F$ to $X$ is a projective equivalence. Since $F$ maps  projective equivalences to projective equivalences this implies that a morphism $f:X\sr Y$ in $\Delta^{op}Rad(C_+)$ is a projective equivalence if and only if $F(f)=(\phi f)_+$ is a projective equivalence which finishes the proof of the proposition.

The fact that $S^{\Phi}((X,x))=(X^n/G,x^n)$ follows easily from the isomorphism (\ref{punp}).
\end{proof}

In the case of the ordinary symmetric products associated with the permutation group $(S_n, Id)$ we will use the simplified notations $S^n$ and $\tilde{S}^n$. For small values of $n$ one has
$$
\begin{array}{ll}
S^0(X_+)=pt&\tilde{S}^0(X_+)=S^0\\
S^1(X_+)=X_+&\tilde{S}^1(X_+)=X_+\\
S^2(X_+)=(X^2/S_2)_+\vee X_+&\tilde{S}^2(X_+)=(X^2/S_2)_+
\end{array}$$
\begin{proposition}
\llabel{minus1}
For any $n>0$ there is a family of natural in $X, Y\in \Delta^{op}Rad(C_+)$ isomorphisms
\begin{eq}
\llabel{2010.1.29.oldeq1}
\tilde{S}^n(X\vee Y)=\vee_{n\ge i\ge 0}(\tilde{S}^iX\wedge \tilde{S}^{n-i}Y).
\end{eq}
\end{proposition}
\begin{proof}
The left hand side of (\ref{2010.1.29.oldeq1}) is the value on $(X,Y)$ of the simplicial extension of the radditive extension of the functor $C_+\times C_+\sr C_+$ of the form $(U_+,V_+)\mapsto \tilde{S}^n(U_+\vee V_+)$, the right hand side is the value on (X,Y) of the simplicial extension of the radditive extension of the functor $(U_+,V_+)\mapsto \vee_{n\ge i\ge 0}(\tilde{S}^i(U_+)\wedge \tilde{S}^{n-i}(U_+))$.  These two functors from $C_+\times C_+$ to $C_+$ are isomorphic due to the formula
$$(U\amalg V)^n/S_n=\amalg_{n\ge i\ge 0} (U^i/S_i\times V^{n-i}/S_{n-i}).$$
Therefore their radditive extensions and then simplicial extensions of radditive extensions are isomorphic as well.   
\end{proof}
Since radditive extensions commute with filtered colimits   Proposition \ref{minus1} immediately implies the following result.
\begin{cor}
\llabel{2009.08.06.cor1}
Let $(X_{\alpha})_{\alpha\in A}$ be a family of objects in $\Delta^{op}Rad(C_+)$. Then one has
$$\tilde{S}^n(\vee_{\alpha\in A} X_{\alpha})=\vee_{k_1\alpha_1+\dots + k_m\alpha_m\in S^nA} (\tilde{S}^{k_1}X_{\alpha_1}\wedge\dots\wedge \tilde{S}^{k_m}X_{\alpha_m}).$$
where $S^nA$ is the $n$-th symmetric power of the set $A$, $k_1+\dots+k_m=n$ and $\alpha_1,\dots,\alpha_m$ are pairwise distinct. 
\end{cor}

Consider a coprojection sequence $X\stackrel{f}{\sr}Y\stackrel{p}{\sr}Z$ in $Rad(C_+)$. Let us say that a morphism $f:A\sr B$ in a category has a strict image if the image of the morphism of representable functors defined by $f$ is representable i.e. if there exists a factorization $A\sr Im(f)\sr B$ of $f$ where the first morphism is a split epimorphism and the second one is a monomorphism. 
\begin{lemma}
\llabel{2009.08.06.l2}
Let $f:X\sr Y$ be a coprojection in $Rad(C_+)$. Then the composition 
$$\vee_{n\ge a\ge i}(\tilde{S}^a(X)\wedge \tilde{S}^{n-a}(Y))\sr \tilde{S}^n(X\vee Y)\stackrel{\tilde{S}^n(f\vee Id)}{\longrightarrow} \tilde{S}^n(Y)$$
has a strict image which we denote by $S_{\ge i}^n(X,Y)\sr \tilde{S}^n(Y)$. A choice of $i:Z\sr Y$ such that $f\vee i:X\vee Z\sr Y$  is an isomorphism defines an isomorphism 
$$S_{\ge i}^n(X,Y)\cong  \vee_{n\ge a\ge i}(\tilde{S}^a(X)\wedge \tilde{S}^{n-a}(Z))$$
over $\tilde{S}^n(Y)$.
\end{lemma}
\begin{proof}
Let us choose $i:Z\sr Y$ such that $f\vee i:X\vee Z\sr Y$  is an isomorphism. We get a commutative square
$$
\begin{CD}
\vee_{n\ge a\ge i}(\tilde{S}^a(X)\wedge \tilde{S}^{n-a}(Z)) @>>> \tilde{S}^n(X\vee Z) @>\cong>> \tilde{S}^n(Y)\\
@VVV @VVV @VV=V\\
\vee_{n\ge a\ge i}(\tilde{S}^a(X)\wedge \tilde{S}^{n-a}(Y)) @>>> \tilde{S}^n(X\vee Y) @>>> \tilde{S}^n(Y)
\end{CD}
$$
To finish the proof it remains to construct for each $j\ge i$ a morphism 
$$\tilde{S}^j(X)\wedge \tilde{S}^{n-j}(Y)\sr \vee_{n\ge a\ge i}(\tilde{S}^a(X)\wedge \tilde{S}^{n-a}(Z))$$
such that the diagram
$$
\begin{CD}
\tilde{S}^j(X)\wedge \tilde{S}^{n-j}(Y) @>>> \tilde{S}^n(Y)\\
@VVV @VV=V\\
\vee_{n\ge a\ge i}(\tilde{S}^a(X)\wedge \tilde{S}^{n-a}(Z)) @>>> \tilde{S}^n(Y)\\
@VVV @VV=V\\
\vee_{n\ge a\ge i}(\tilde{S}^a(X)\wedge \tilde{S}^{n-a}(Y)) @>>> \tilde{S}^n(Y)
\end{CD}
$$
commutes. We define this morphism as the composition
$$\tilde{S}^j(X)\wedge \tilde{S}^{n-j}(Y)\cong \tilde{S}^j(X)\wedge \tilde{S}^{n-j}(X\vee Z)=\vee_{n-j\ge k\ge 0} (\tilde{S}^j(X)\wedge \tilde{S}^{k}(X)\wedge \tilde{S}^{n-k-j}(Z))\sr$$ 
$$\vee_{n\ge a\ge j} (\tilde{S}^{a}(X)\wedge \tilde{S}^{n-a}(Z))\sr \vee_{n\ge a\ge i} (\tilde{S}^{a}(X)\wedge \tilde{S}^{n-a}(Z))$$
where $a=k+j$. 
\end{proof}
By construction we obtain a sequence of morphisms 
$$\tilde{S}^n(X)=\tilde{S}^n_{\ge n}(X,Y)\sr \dots \tilde{S}^n_{\ge 1}(X,Y)\sr \tilde{S}^n_{\ge 0}(X,Y)=\tilde{S}^n(Y)$$
As a corollary of the second statement of Lemma \ref{2009.08.06.l2} we see that the morphisms
$$\tilde{S}^n_{\ge i+1}(X,Y)\sr \tilde{S}^n_{\ge i}(X,Y)$$ 
are coprojections. The following result identifies the cofibers of these coprojections.
\begin{lemma}
\llabel{2009.08.07.l1}
For any coprojection sequence $X\sr Y\sr Z$  there are natural morphisms $\tilde{S}^n_{\ge i}(X,Y)\sr \tilde{S}^i(X)\wedge \tilde{S}^{n-i}(Z)$ such that the sequences
$$\tilde{S}^n_{\ge i+1}(X,Y)\sr \tilde{S}^n_{\ge i}(X,Y)\sr \tilde{S}^i(X)\wedge \tilde{S}^{n-i}(Z)$$
are coprojection sequences.
\end{lemma}
\begin{proof}
Note first that it is clear from Lemma \ref{2009.08.06.l2} that a choice of $i:Z\sr Y$ defines a coprojection sequence of the required form. The point of the lemma is to show that this sequence is independent on the choice of $i$. 

To do so let us show that the obvious morphism 
\begin{eq}
\llabel{2009.08.07.eq2}
\vee_{n\ge a\ge i}(\tilde{S}^a(X)\wedge \tilde{S}^{n-a}(Y)) \sr \tilde{S}^i(X)\wedge \tilde{S}^{n-i}(Z)\end{eq} 
factors through the projection
\begin{eq}
\llabel{2009.08.07.eq1}
\vee_{n\ge a\ge i}(\tilde{S}^a(X)\wedge \tilde{S}^{n-a}(Y)) \sr \tilde{S}_{\ge i}^n(X,Y)
\end{eq}
The proof of Lemma \ref{2009.08.06.l2} implies easily that this projection identifies $\tilde{S}_{\ge i}^n(X,Y)$ with the image of the projector on the left hand side of (\ref{2009.08.07.eq1}) of the form
$$\vee_{n\ge a\ge i}(\tilde{S}^a(X)\wedge \tilde{S}^{n-a}(Y))\cong \vee_{n\ge a\ge i}\vee_{n-a\ge j\ge 0} \tilde{S}^a(X)\wedge \tilde{S}^{n-a-j}(X)\wedge \tilde{S}^j(Z)\sr$$
$$\sr \vee_{n\ge b\ge i} \tilde{S}^{b}(X)\wedge \tilde{S}^{n-b}(Y)$$
which maps $\tilde{S}^a(X)\wedge \tilde{S}^{n-a-j}(X)\wedge \tilde{S}^j(Z)$ to $\tilde{S}^{n-j}(X)\wedge \tilde{S}^j(Y)$. The compositions of both the identity and this projector with (\ref{2009.08.07.eq2}) are zero on the summands $\tilde{S}^a(X)\wedge \tilde{S}^{n-a-j}(X)\wedge \tilde{S}^j(Z)$ with $j\ne n-i$ and coincide with the canonical isomorphisms
$$\tilde{S}^i(X)\wedge \tilde{S}^{0}(X)\wedge \tilde{S}^{n-i}(Z)\sr \tilde{S}^i(X)\wedge \tilde{S}^{n-i}(Z)$$
on the only summand with $j=n-i$. Therefore (\ref{2009.08.07.eq2}) factors through (\ref{2009.08.07.eq1}). The fact that the resulting morphism extends the morphism $\tilde{S}^n_{\ge i+1}(X,Y)\sr\tilde{S}^n_{\ge i}(X,Y)$ to a coprojection sequence is straightforward.  
\end{proof}
For $0\le a\le b\le n$ define objects $\tilde{S}^n_{a,b}(X,Y)$ by the coprojection sequences
$$\tilde{S}_{\ge b}^n(X,Y)\sr \tilde{S}_{\ge a}^n(X,Y)\sr \tilde{S}_{a,b}^{n}(X,Y)$$
Summarizing the previous discussion and extending it to simplicial objects we get the following theorem which described the behavior of $\tilde{S}^n$ with respect to coprojection sequences.
\begin{theorem}
\llabel{2009.08.07.th1}
Any term-wise coprojection sequence 
$$X\sr Y\sr Z$$
in $\Delta^{op}Rad(C_+)$ defines in a natural way a collection of objects $\tilde{S}^n_{a,b}(X,Y)$ for $0\le a\le b\le n$, natural isomorphisms
$$\tilde{S}_{0,n}^n(X,Y)=\tilde{S}^n(Y)$$
and for any $a=0,\dots,n$
$$\tilde{S}^{n}_{a,a}(X,Y)=\tilde{S}^a(X)\wedge \tilde{S}^{n-a}(Z)$$
and for any $0\le a\le b\le c\le n$ coprojection sequences
\begin{eq}
\llabel{2009.08.07.eq3}
\tilde{S}^{n}_{b+1,c}(X,Y)\sr \tilde{S}^n_{a,c}(X,Y)\sr \tilde{S}^n_{a,b}(X,Y)
\end{eq}
\end{theorem}
The following corollary describes a particularly useful in applications tower of coprojection sequences of the form (\ref{2009.08.07.eq3}).
\begin{cor}
\llabel{2009.08.07.cor1}
Under the assumptions of the theorem there is a tower of coprojection sequences of the form
$$\begin{array}{l}
\tilde{S}^n(X)\sr \tilde{S}^n(Y)\sr \tilde{S}^n_{0,n-1}(X,Y)\\\\
\tilde{S}^i(X)\wedge \tilde{S}^{n-i}(Z) \sr \tilde{S}^n_{0,i}(X,Y)\sr \tilde{S}^n_{0,i-1}(X,Y)\,\,\,\,\,i=n-1,\dots,2\\\\
X\wedge \tilde{S}^{n-1}Z\sr \tilde{S}^n_{0,1}(X,Y)\sr \tilde{S}^n(Z)
\end{array}
$$
\end{cor}
\begin{proof}
These are coprojection sequences (\ref{2009.08.07.eq3}) for $a=0$ and $c=b+1$.
\end{proof}
Using Proposition \ref{2009ref6} one can easily reformulate an analog of Theorem \ref{2009.08.07.th1} and Corollary \ref{2009.08.07.cor1} for cofiber sequences in any of the homotopy categories of $C_+$ which we have considered.

The following lemma is straightforward.
\begin{lemma}
\llabel{ntilde}
For any $n>0$ there is a family of natural in $X\in C_+$  coprojection sequences 
\begin{eq}\llabel{nmom}
S^{n-1}(X)\sr S^n(X)\sr \tilde{S}^n(X).
\end{eq}
\end{lemma}
As a corollary of the fact that the isomorphisms and sequences of the lemmas are natural we conclude that they extends to objects of $C_+^{\#}$ and further to objects of $\Delta^{op}C_+^{\#}$. 

We let $S^{\infty}(X)$ denote the colimit of the sequence 
$$S^0(-)\sr S^1(-)\sr S^2(-)\sr\dots$$
in $C_+^{\#}$. Note that for $X\in C$ one has
$$S^{\infty}(X_+)=(\coprod_{n\ge 1} X^n/S_n)_+=\vee_{n\ge 1} \tilde{S}^n(X_+)$$
however his decomposition is natural only for morphisms of the form $f_+$ where $f$ is a morphism in $C$ and not for general morphisms in $C_+$. 

We will also consider the functors $S^{\infty}[1/d]$ for integers $d>0$ defined by the formula
$$S^{\infty}[1/d](X)=colim(S^{\infty}(X)\stackrel{\times d}{\longrightarrow}S^{\infty}(X)\stackrel{\times d}{\longrightarrow}\dots)$$
where $\times d$ is the multiplication by $d$ map with respect to the abelian monoid structure of $S^{\infty}$.  
\begin{lemma}
\llabel{wasnoth}
For any $X,Y\in C_+^{\#}$ there is a natural isomorphism
\begin{eq}
\llabel{wasnoth1}
S^{\infty}[1/d](X\vee Y)=S^{\infty}[1/d](X)\times S^{\infty}[1/d](Y).
\end{eq}
\end{lemma}
\begin{proof}
The maps $X\vee Y\sr X$ and $X\vee Y\sr Y$ define the map from the left to the right hand side of (\ref{wasnoth1}). To verify that it is an isomorphism we may assume that $X=U_+$ and $Y=V_+$ for $U,V\in C$. The case $d=1$ follows then immediately from Proposition \ref{minus1}. The case $d>1$ follows from the case $d=1$ and the fact that finite products commute with filtered colimits.
\end{proof}
Corollary \ref{2009.08.05.cor1} together with the fact that $(Nis,\af)$-equivalences are closed under filtered colimits implies the following result.
\begin{cor}
\llabel{symonhinfty}
For any $d>0$  there are unique functors $\LL S^{\infty}[1/d]$ on the $Nis$- and $(Nis,\af)$-homotopy categories which are determined by the conditions that the square 
$$
\begin{CD}
\Delta^{op}C_+^{\#} @>S^{\infty}[1/d]>> \Delta^{op}C_+^{\#}\\
@VVV @VVV\\
H_{Nis,\af}(C_+) @>\LL S^{\infty}[1/d]>> H_{Nis,\af}(C_+)
\end{CD}
$$
and its analog for the $Nis$-local category, commute. In addition, for any ind-solid object $F$ in $\Delta^{op}Rad(C_+)$ one has $\LL S^{\infty}[1/d](F)=S^{\infty}[1/d](F)$ in both $Nis$- and $(\af,Nis)$-homotopy categories. 
\end{cor}

\begin{remark}\rm
\llabel{rem3.24}
Proposition \ref{2009.08.05.pr1} shows that if $(X,x)$ is a radditive functor represented by a pointed scheme then $S^{\infty}((X,x))$ is ind-represented by the usual infinite symmetric power $colim_{i\ge 0} S^n((X,x))$ and the morphism $\LL S^{\infty}((X,x))\sr S^{\infty}((X,x))$ is in this case a projective equivalence.
\end{remark}

\subsection{Generalized symmetric powers on $H_{Nis,\af}(Cor(C,R))$}

We continue to assume that the underlying category $C$ is $f$-admissible.
Set
$$c=c(k)=\left\{\begin{array}{ll}
1&{\rm if}\,\,char(k)=0\\
char(k)&{\rm otherwise}
\end{array}
\right.
$$ 
The number $c(k)$ is sometimes called the characteristic exponent of $k$.  In what follows we assume that  $c$ is invertible in our ring of coefficients $R$. We start with the following result.
\begin{proposition}
\llabel{start}
Let $X,Y\in C$ and let $G$ be a finite group acting on $X$. Then 
$$Hom_{Cor(C,R)}([X/G],[Y])=Hom_{Cor(C,R)}([X],[Y])^G$$
i.e. $[X/G]$ is the categorical quotient for the action of $G$ on $[X]$.
\end{proposition}
\begin{proof}
It follows from the fact that the functor represented by $[Y]$ on $Sch/k$ is a qfh-sheaf by \cite[Proposition 4.2.7]{SusVoe2new} and that for qfh-sheaves $F$ one has $F(X/G)=F(X)^G$.
\end{proof}
\begin{proposition}
\llabel{ext}
For any permutation group $\Phi=(G,\phi:G\sr S_n)$ there exist a unique (up to a canonical isomorphism) functor $S^{\Phi}_{tr}:Cor(C,R)\sr Cor(C,R)$ such that the square
$$
\begin{CD}
C_+ @>\tilde{S}^{\Phi}>> C_+\\
@V\Lambda^l VV @VV\Lambda^l V\\
Cor(C,R) @>S^{\Phi}_{tr}>> Cor(C,R)
\end{CD}
$$
commute.  
\end{proposition}
\begin{proof}
In view of Proposition \ref{start}  we may define $S^{\Phi}_{tr}$ by the rule
$$S^{\Phi}(X_+)=[X]^{\oo n}/G.$$
One verifies easily that the required squares commute. 
\end{proof}
\begin{remark}\rm
\llabel{nnplustr}
One can also define functors on $Cor(C,R)$ corresponding to the un-reduced symmetric powers $S^{\Phi}$.  However, some of the important natural transformations between these functors on $C_+$ do not extend to natural transformation over $Cor$. For example, the embeddings $S^n(X)\sr S^{n+1}(X)$ are not natural with respect to morphisms in $Cor$ which can be seen by looking at the morphism $d\cdot Id:[X]\sr [X]$ for $d>1$.
\end{remark}
Let $\Phi$ be as above, $i:H\sr G$ a subgroup of $G$ and $\Psi$ the permutation group $(H,\psi=\phi\circ i)$. Assume for a moment that $H$ is normal in $G$ and consider finite correspondences with coefficients in a commutative ring $R$ such that $d=[G:H]$ is invertible in $R$. Then for any $X$ in $Rad(Cor(C,R))$ there is an action of $G/H$ on $S^{\Psi}_{tr}(X)$ and it is more or less obvious that $S^{\Phi}_{tr}(X)$ is the direct summand of $G/H$-invariants in $S^{\Psi}_{tr}(X)$. We will need an analog of this observation in the case when $H$ is not necessarily normal in $G$. 

For any $g\in G$ let $\Psi_g$ be the permutation group corresponding to the subgroup $H\cap gHg^{-1}$ of $G$. Then for any $X$ in $Cor(C,R)$ there are two morphisms 
$$p:S^{\Psi_g}_{tr}(X)\sr S^{\Psi}_{tr}(X)$$
$$p':S^{\Psi_g}_{tr}(X)\sr S^{\Psi}_{tr}(X)$$
where 
$$p:X^{\oo n}/(H\cap gHg^{-1})\sr X^{\oo n}/H$$
is the projection and $p'$ is the map whose composition with $X^{\oo n}\sr X^{\oo n}/(H\cap gHg^{-1})$ is $x\mapsto gx$ followed by the projection. 
\begin{theorem}
\llabel{nonnormal}
Let $d=[G:H]$ be invertible in the ring of coefficients $R$. Then for any $X$ in $Cor(C,R)$ there is a split cocomplete-equalizer sequence:
$$\oplus_{g\in G} S^{\Psi_g}_{tr}(X)\dsr S^{\Psi}_{tr}(X)\sr S^{\Phi}_{tr}(X)$$
where the two arrows are given by $p$ and $p'$ on each summand. 
\end{theorem}
\begin{proof}
In view of Proposition \ref{start} the theorem is a particular case of  Proposition \ref{splitc}.
\end{proof}
\begin{cor}
\llabel{nnormal}
Under the assumptions of the theorem assume in addition that $H$ is normal in $G$. Then 
$$S_{tr}^{\Phi}(X)=(S_{tr}^{\Psi}(X))^{G/H}$$
i.e. $S_{tr}^{\Phi}(X)$ is the image of the projector $d^{-1}\sum_{u\in G/H} u$ acting on $S_{tr}^{\Psi}(X)$.
\end{cor}

As in the case of the symmetric power functors on $C_+$ we will write $S^{\Phi}_{tr}$ instead of $(S^{\Phi}_{tr})^{rad}$ for the radditive extensions of these functors. Note that Theorem \ref{nonnormal} and Corollary \ref{nnormal} immediately extend to the functors $S_{tr}^{\Phi}$ on $Rad(Cor(C,R))$ and further on $\Delta^{op}Rad(Cor(C,R))$. 
\begin{remark}\rm
Note however, that the definition of $S^{\Phi}_{tr}(X)$ as the categorical quotient of $X^{\oo n}$ by the permutational action of $G$ extends to $Cor(C,R)^{\#}$ but not to $Rad(Cor(C,R))$. Even for $X\in Cor(C,R)$, the quotient $S^{\Phi}_{tr}(X)$ is not the quotient of the radditive functor $X^{\oo n}$ by the action of $G$ in $Rad(Cor(C,R))$. The later exists and maps to $S^{\Phi}_{tr}(X)$ but this map is almost never an isomorphism or even an $(Nis,\af)$-equivalence unless the projection $X^{\oo n}\sr S^{\Phi}_{tr}(X)$ splits.  
\end{remark}
\begin{theorem}
\llabel{prestr}
For any $R$ such that $c(k)$ is invertible in $R$ and any $\Phi=(G,\phi:G\sr S_n)$ the functors $S^{\Phi}_{tr}$ take $Nis$- (resp. $(Nis,\af)$-) equivalences between objects of $\Delta^{op}Cor(C,R)^{\#}$ to $Nis$- (resp. $(Nis,\af)$-) equivalences. Therefore, there are unique functors $\LL S_{tr}^{\Phi}$ on the $Nis$- and $(Nis,\af)$-homotopy categories which are determined by the conditions that the squares 
$$
\begin{CD}
\Delta^{op}Cor(C,R)^{\#} @>S_{tr}^{\Phi}>> \Delta^{op}Cor(C,R)^{\#}\\
@VVV @VVV\\
H_{Nis,\af}(Cor(C,R)) @>\LL S_{tr}^{\Phi}>> H_{Nis,\af}(Cor(C,R))
\end{CD}
$$
and their analogs for the $Nis$-local categories, commute.
\end{theorem}
\begin{proof}
Let $X,Y$ be objects  of $\Delta^{op}Cor(C,R)^{\#}$ and $f:X\sr Y$ be a $Nis$-equivalence. Then by Theorem \ref{maindelta} we have $f\in cl_{\bdl}([G_{Nis}])$. As any radditive extension, the functor $S^{\Phi}_{tr}$ takes $\bdl$-closures to $\bdl$-closures. On the other hand $[G_{Nis}]=\Lambda^l((G_{Nis})_+)$ and 
$$S_{tr}^{\Phi}(\Lambda^l((G_{Nis})_+))=\Lambda^l(\tilde{S}^{\Phi}((G_{Nis})_+))\subset \Lambda^l(cl_l((G_{Nis})_+)\cap \Delta^{op}C_+^{\#})\subset cl_l((G_{Nis})_+)$$
where the equality holds by Proposition \ref{ext}, the first inclusion by  Theorem \ref{mainpreserve} and the second inclusion by  Theorem \ref{2009th1}. The same argument applies to the case of $(Nis,\af)$-equivalences.
\end{proof}

Let us consider now in more detail the case of ordinary symmetric powers $S^n_{tr}$. Generalizing a particular case of Corollary \ref{nnormal} to radditive functors we get:
\begin{lemma}
\llabel{summand2}
Let $n$ be an integer and $R$ be a ring where $n!$ is invertible. Then for any $X$ in $\Delta^{op}Rad(Cor(C,R))$ the obvious morphism $X^{\oo n}\sr S_{tr}^n(X)$ defines an isomorphism between $S_{tr}^n(X)$ and the image of the projector $(1/n!)\sum_{\sigma\in S_n} \sigma$ on $X^{\oo n}$.
\end{lemma}
\begin{proposition}
\llabel{directsummand}
Let $n$ be an integer, $l$ a prime and $n=\sum n_i l^i$ the $l$-primary decomposition of $n$. Let further $R$ be an $l$-local ring (i.e. a ring where all primes but $l$ are invertible). Then for any $X$ in $\Delta^{op}Rad(Cor(C,R))$ there is a split epimorphism
\begin{eq}
\llabel{willsplit}
p_{n,l}:\oo_{i} ((S_{tr}^{l})^{\circ i}(X))^{\oo n_i}\sr S_{tr}^n(X)
\end{eq}
such that both $p_{n,l}$ and its section are natural in $X$. 
\end{proposition}
\begin{proof}
Denote temporarily by $\Phi_l$ the permutation group $(S_l,Id)$.  For a sequence of non-negative integers $\uu{q}=(q_1,\dots,q_k)$ consider the permutation group
$$\Phi_{l,\uu{q}}=\prod_i (\Phi_l^{* i})^{\times q_i}=(G_{l,\uu{q}}, \phi_{l,\uu{q}}:G_{l,\uu{q}}\sr S_n)$$
where $n=\sum q_i l^i$. 

By Propositions \ref{wreath} and \ref{product} the left hand side of (\ref{willsplit}) is canonically isomorphic to $\tilde{S}^{\Phi}(X)$ where $\Phi=\Phi_{l,\uu{n}}(X)$ for $\uu{n}=(n_0,\dots,n_i,\dots)$. The morphism $p_{n,l}$ is associated with the embedding $G_{l,\uu{n}}\sr S_n$. Using the fact that $\uu{n}$ is the $l$-primary decomposition of $n$ and computing how many times $l$ divides $n!$ one concludes that $[S_n:G_{l,\uu{n}}]$ is prime to $l$ and therefore invertible in $R$. Our result follows now from Theorem \ref{nonnormal}.
\end{proof} 

\begin{proposition}
\llabel{directsums}
For any $n\ge 0$ there is a natural in $X, Y\in \Delta^{op}Rad(Cor(C,R))$ family of  isomorphisms of the form 
$$S^n_{tr}(X\oplus Y)=\oplus_{i\ge 0} S^i_{tr}(X)\oo S^{n-i}_{tr}(Y)$$
\end{proposition}
\begin{proof}
The same reasoning as in the proof of Proposition \ref{minus1} shows that it is sufficient to consider the case when $X,Y\in Cor(C,R)$. Then we have morphisms
$$S^i_{tr}(X)\oo S^{n-i}_{tr}(Y)\sr S^n_{tr}(X\oplus Y)$$
which are obvious from the definition of $S^i_{tr}$ as $(-)^{\oo i}/S_i$. These morphisms are clearly natural in $X$ and $Y$. On the other hand they are compatible with the morphisms which define the isomorphism of Proposition \ref{minus1} and therefore their sum gives an isomorphism. 
\end{proof}
\begin{cor}
\llabel{infdirectsums}
Let $(X_{\alpha})_{\alpha\in A}$ be a family of objects in $\Delta^{op}Rad(Cor(C,R))$. Then one has
$$S^n_{tr}(\bigoplus_{\alpha\in A} X_{\alpha})=\bigoplus_{k_1\alpha_1+\dots + k_m\alpha_m\in S^nA} S_{tr}^{k_1}X_{\alpha_1}\oo\dots\oo S^{k_m}_{tr}X_{\alpha_m}.$$
where $S^nA$ is the $n$-th symmetric power of the set $A$, $k_1+\dots+k_m=n$ and $\alpha_1,\dots,\alpha_m$ are pairwise distinct. 
\end{cor}
By exactly the same reasoning as in the proof of Theorem \ref{2009.08.07.th1} we get the following results.
\begin{theorem}
\llabel{2009.08.07.th2}
Any term-wise coprojection sequence 
$$X\sr Y\sr Z$$
in $\Delta^{op}Rad(Cor(C,R))$ defines in a natural way a collection of objects $S^n_{tr,a,b}(X,Y)$ for $0\le a\le b\le n$, natural isomorphisms
$$S_{tr,0,n}^n(X,Y)=S^n_{tr}(Y)$$
and for any $a=0,\dots,n$
$$S^{n}_{tr,a,a}(X,Y)=S_{tr}^a(X)\oo S_{tr}^{n-a}(Z)$$
and for any $0\le a\le b\le c\le n$ coprojection sequences
\begin{eq}
\llabel{2009.08.07.eq5}
S^{n}_{tr,b+1,c}(X,Y)\sr S^n_{tr,a,c}(X,Y)\sr S^n_{tr,a,b}(X,Y)
\end{eq}
\end{theorem}
The following corollary describes a particularly useful in applications tower of coprojection sequences of the form (\ref{2009.08.07.eq3}).
\begin{cor}
\llabel{2009.08.07.cor2}
Under the assumptions of the theorem there is a tower of coprojection sequences of the form
$$\begin{array}{l}
S^n_{tr}(X)\sr S_{tr}^n(Y)\sr S^n_{tr,0,n-1}(X,Y)\\\\
S_{tr}^i(X)\oo S_{tr}^{n-i}(Z) \sr S^n_{tr,0,i}(X,Y)\sr S^n_{tr,0,i-1}(X,Y)\,\,\,\,\,i=n-1,\dots,2\\\\
X\oo S_{tr}^{n-1}Z\sr S^n_{tr,0,1}(X,Y)\sr S_{tr}^n(Z)
\end{array}
$$
\end{cor}
\begin{proof}
These are coprojection sequences (\ref{2009.08.07.eq5}) for $a=0$ and $c=b+1$.
\end{proof}

Using Proposition \ref{2009ref6} one can easily reformulate an analog of Theorem \ref{2009.08.07.th2} and Corollary \ref{2009.08.07.cor2} for cofiber sequences in any of the homotopy categories of $Cor(C,R)$ which we have considered. Below we do it in the case of $H_{Nis,\af}(Cor(C,R))$.

\begin{proposition}
\llabel{scase}
\llabel{exactseq}
For any $n>0$ and any cofiber sequence 
$$X\sr Y\sr Z\sr \Sigma^1X$$
in $H_{Nis,\af}(Cor(C,R))$ there are objects $\LL S^n_{tr,a,b}(X,Y)$ such that for $0\le a\le n$ one has 
$$\LL S^n_{tr,a,a}(X,Y)=\LL S^a_{tr}X\oo_{\LL} \LL S^{n-a}_{tr}Y,$$
and 
$$\LL S^n_{tr,0,n}(X,Y)=\LL S^n_{tr}Y,$$
and cofiber sequences
\begin{eq}
\llabel{obm2}
\LL S^n_{tr, k+1,j}(X,Y)\sr \LL S^n_{tr,m,j}(X,Y)\sr \LL S^n_{tr, m,k}(X,Y)\sr \Sigma^1\LL S^n_{tr, k+1,j}(X,Y).
\end{eq}
\end{proposition}
\begin{cor}
\llabel{towers}
Under the assumptions of the proposition there is a tower of cofiber sequences
$$\begin{array}{l}
\LL S^n_{tr}X\sr \LL S^n_{tr}Y\sr \LL S^n_{tr,0,n-1}(X,Y)\sr \Sigma^1\LL S^n_{tr}X\\\\
\LL S^i_{tr}X\oo_{\LL} \LL S^{n-i}_{tr}Z\sr \LL S^n_{tr,0,i}(X,Y)\sr \LL S^n_{tr,0,i-1}(X,Y)\sr \Sigma^1(\LL S^i_{tr}X\oo_{\LL} \LL S^{n-i}_{tr}Z)\\\\
\mbox{\rm for $i=n-1,\dots, 2$}\\\\
X\oo_{\LL} \LL S^{n-1}_{tr}Z\sr \LL S^n_{tr, 0,1}(X,Y)\sr \LL S^n_{tr}Z\sr \Sigma^1(X\oo_{\LL} \LL S^{n-1}_{tr}Z)
\end{array}
$$
\end{cor}

The usual problems with the lifting of morphisms of cofiber sequences in homotopy categories to morphisms in the model categories prevent us from asserting that the cofiber sequences of symmetric products of Theorem \ref{scase} and Corollary \ref{towers} are natural with respect to morphisms of the base cofiber sequence. Due to this issue and in view of Proposition \ref{2009ref6}  we will formulate the results below in the context of coprojection rather than cofiber sequences. 
\begin{definition}
\llabel{oddeven}
An object $X$ of $H_{Nis,\af}(Cor(C,R))$ is called even (resp. odd) if the permutation isomorphism $\sigma:X\oo X\sr X\oo X$ is the identity (resp. the multiplication by $-1$).
\end{definition}
We have the following obvious fact.
\begin{lemma}
\llabel{oddprop} Tensor product of two odd or two even objects is even. Tensor product of an odd and an even object is odd.
\end{lemma}
\begin{lemma}
\llabel{suspodd}
Let $X$ be an even (resp. odd) object. Then $\Sigma^1_sX$ is odd (resp. even).
\end{lemma}
\begin{proof}
It follows from the fact that $\Sigma^1_sX=X\oo S^1_s$ and that $S^1_s$ is  odd. 
\end{proof}
\begin{lemma}
\llabel{powerofodd}
Let $R$ be an $l$-local ring, $X$ an object of $H_{Nis,\af}(Cor(C,R))$ and $1<n<l$ an integer. Then one has:
\begin{enumerate}
\item if $X$ is  odd then $S^n_{tr}(X)=0$,
\item if $X$ is  even then the map $X^{\oo n}\sr S^n_{tr}(X)$ is an isomorphism.
\end{enumerate}
\end{lemma}
\begin{proof}
Since the projection $\Delta^{op}Cor(C,R)^{\#}\sr H_{Nis,\af}$ is an additive functor  Proposition \ref{summand2} implies hat $S^n_{tr}(X)$ as an object of $H_{Nis,\af}$ is the image of the averaging projector. Therefore for an even $X$ we get $X$. The number of elements in $S_n$ is even for $n>1$ and therefore for an odd $X$ and $n>1$ we get zero. 
\end{proof}
Let $\Psi$ be the permutation group $(G, i:G\sr S_l)$ where $G=\zz/l$ is embedded into $S_l$ as the subgroup generated by the cycle $(1\dots l)$. The symmetric power $S_{tr}^{\Psi}$ associated with $\Psi$ is  the $l$-th cyclic power. The quotient $N(G)/G$ where  $N(G)$ is the normalizer of $G$ in $S_l$ is canonically isomorphic to $Aut(\zz/l)=(\zz/l)^*$. The definition of $S^{\Psi}_{tr}$ shows that this quotient acts on $S^{\Psi}_{tr}(X)$ for any $X$ in a manner natural in $X$. If $l-1$ is invertible in $R$ then we may define the direct summand $S^{\Psi}_{tr}(X)/Aut(\zz/l)$ of coinvariants of this action and Corollary \ref{nnormal} implies that  it is naturally isomorphic to $S^{\Psi'}(X)$ where $\Psi'$ corresponds to $N(G)\subset S_l$. In particular, the we get a natural morphism $S^{\Psi}_{tr}(X)/Aut(\zz/l)\sr S^l_{tr}(X)$.
\begin{proposition}
\llabel{2009.08.10.prop3}
Let $R$ be an $l$-local ring and $X$ an object of $\Delta^{op}Cor(C,R)^{\#}$ which is even in $H_{Nis,\af}(Cor(C,R))$. Then the morphism $S^{\Psi}_{tr}(X)/Aut(\zz/l)\sr S^l_{tr}(X)$ is an $(Nis,\af)$-equivalence.
\end{proposition}
\begin{proof}
Let us apply Theorem \ref{nonnormal} to the natural transformation of functors $S_{tr}^{\Psi}\sr S_{tr}^l$ on $X$. Since $\zz/l$ has no non-trivial subgroups,  one has $gGg^{-1}=G$ or $G\cap gGg^{-1}=\{e\}$ for any $g\in S_l$. Therefore, we can re-write the first term of the split coequalizer diagram of Theorem \ref{nonnormal} as $(\oplus_{g\in N(G)}S^{\Psi}_{tr}(X))\oplus  (\oplus_{g\in S_l\backslash N(G)}X^{\oo l})$ with the two arrows differing by the action of $N(G)$ on $S^{\Psi}_{tr}(X)$ on the summands of the first type and by the action of $S_l$ on $X^{\oo l}$ on the summands of the second type. Since split equalizer diagrams are absolute i.e. preserved by all functors (see e.g. \cite[p.149]{McLane2}) the image of this diagram in $H_{Nis,\af}(Cor(C,R))$ is still a split coequalizer diagram. Since the action of $S_l$ on $X^{\oo l}$ in this category is trivial we conclude that $S^l_{tr}(X)=S^{\Psi}_{tr}(X)/Aut(G)$. 
\end{proof}
\begin{proposition}
\llabel{oe1}
Let $R$ be an $l$-local ring and $n\le l$. Let 
$$X\sr Y\sr Z\sr\Sigma^1X$$
be a coprojection sequence in $\Delta^{op}Cor(C,R)^{\#}$. Then one has:
\begin{enumerate}
\item If $X$ is odd in $H_{Nis,\af}(Cor(C,R))$ then there are natural coprojection sequences  of the form
\begin{eq}
\llabel{x22}
S^n_{tr}(X) \sr S^n_{tr}(Y)\sr S^n_{tr,0,n-1}(X,Y)
\end{eq}
\begin{eq}
\llabel{x21}
X\oo S^{n-1}_{tr}(Z) \sr S^n_{tr,0,1}(X,Y)\sr S^n_{tr}(Z)
\end{eq}
and a natural $(Nis,\af)$-equivalence $S^n_{tr,0,n-1}(X,Y)\sr S^n_{tr,0,1}(X,Y)$.
\item If $Z$ is odd then there are natural coprojection sequences in of the form
\begin{eq}
\llabel{z22}
S^n_{tr}(X) \sr S^n_{tr}(Y)\sr S^n_{tr,0,n-1}(X,Y)
\end{eq}
\begin{eq}
\llabel{z21}
S^{n-1}_{tr}(X)\oo Z \sr S^n_{tr,0,n-1}(X,Y)\sr S^n_{tr,0,n-2}(X,Y)
\end{eq}
and a natural $(Nis,\af)$-equivalence $S^n_{tr,0,n-2}(X,Y)\sr S^n_{tr}(Z)$.
\end{enumerate}
\end{proposition}
\begin{proof}
The sequence (\ref{x22}) is the first one and the sequence (\ref{x21}) the last one of the family of sequences in Corollary \ref{towers}.  The remaining sequences of the family give us morphisms 
$$S^n_{tr,0,n-1}(X,Y)\sr \dots\sr S^n_{tr,0,1}(X,Y)$$
which are $(Nis,\af)$-equivalences by our assumption on $X$, Lemma \ref{powerofodd} and Proposition \ref{2009ref6}(3). 

Similarly, the sequence (\ref{z22}) is the first and (\ref{z21}) the second of the sequences of  Corollary \ref{towers}. The rest of the sequences define morphisms
$$S^n_{tr,0,n-2}(X,Y)\sr \dots\sr S^n_{tr,0,0}(X,Y)=S^n_{tr}(Z)$$
which are $(Nis,\af)$-equivalences by our assumption on $Z$, Lemma \ref{powerofodd} and Proposition \ref{2009ref6}(3).
\end{proof}
\begin{cor}
\llabel{homexact}
Let $l$ be a prime, $R$ be an $l$-local ring and $X$ is an object of $\Delta^{op}Cor(C,R)^{\#}$. Then one has:
\begin{enumerate}
\item If $l=2$ or $X$ is  odd in $H_{Nis,\af}(Cor(C,R))$ then there is a cofiber sequence of the form
\begin{eq}
\llabel{seq2}
\Sigma^{l-1}_sX^{\oo l}\sr \Sigma^1S^l_{tr}(X)\sr S^l_{tr}(\Sigma^1_sX)\sr \Sigma^l_sX^{\oo l}.
\end{eq}
\item If $l=2$ or $X$ is even in $H_{Nis,\af}(Cor(C,R))$ then there is a cofiber sequence of the form
\begin{eq}
\llabel{seq1}
\Sigma^1_sX^{\oo l}\sr \Sigma^1S^l_{tr}(X)\sr S^l_{tr}(\Sigma^1_sX)\sr \Sigma^2_sX^{\oo l}.
\end{eq}
\end{enumerate}
\end{cor}
\begin{proof}
There is an obvious coprojection sequence
\begin{eq}
\llabel{seq0}
X \sr Cone(X)\sr \Sigma^1_sX
\end{eq}
where $Cone(X)$ is the simplicial cone of $X$. Applying to this sequence Proposition \ref{oe1} and using Lemma \ref{powerofodd}(2) and Lemma \ref{suspodd} we get the required sequences. 
\end{proof}

We will now consider the finite correspondence versions of the infinite symmetric power $S^{\infty}$. 
\begin{proposition}
\llabel{split}
For any $X$ in $C_+$, any $n\ge 0$ and any $R$ the sequences
$$\Lambda^l_{R}S^{n}(X)\sr \Lambda^l_{R}S^{n+1}(X)\sr \Lambda^l_{R}\tilde{S}^{n+1}(X)=S^{n+1}_{tr}(X)$$
are split exact in a manner natural in $X$. In particular, there are natural in $X\in C_+$ isomorphisms
$$\Lambda^l_{R}(S^{\infty}(X))\cong \oplus_{i\ge 1}S^i_{tr}(\Lambda^l_R(X)).$$
\end{proposition}
\begin{proof}
 Let $*:Spec(k)\sr X$ denote the distinguished point of $X$ and $S^{n-1}(X)\sr S^{n}(X)$ the the inclusions given on the level of the products by $(x_1,\dots,x_{n-1})\mapsto (x_1,\dots,x_{n-1},*)$. We need to construct maps
$$\Lambda^lS^{n}(X)\sr \Lambda^lS^{n-1}(X)$$
which split the morphisms in $Cor$ defined by these inclusions. Since 
$$\Lambda^lS^{n}(X)=\Lambda^l(X^{n}/S_{n}) =\Lambda^l(X^{n})/S_{n}$$
(by Proposition \ref{start}) it suffice to construct maps $s_{n}:\Lambda^lX^{n}\sr \Lambda^lX^{n-1}$ such that 
\begin{enumerate}
\item $\Lambda^lX^{n}\stackrel{s_{n}}{\sr} \Lambda^lX^{n}\sr (\Lambda^lX^{n-1})/S_{n-1}$ is invariant under the action of $S_{n}$ on $X^{n}$,
\item $\Lambda^lX^{n-1}\stackrel{Id\times *}{\sr}\Lambda^lX^{n}\stackrel{s_{n}}{\sr}\Lambda^lX^{n-1}$ is the identity.
\end{enumerate}
Let $\{n\}$ be the set $\{1,\dots,n\}$. Any map $f:\{m\}\sr \{n\}$ defines in the obvious way a map $f_X:X^{n}\sr X^{m}$. For $n>0$, define $s_n:\Lambda^lX^{n}\sr \Lambda^lX^{n-1}$ by the formula
$$s_n=\sum_{0\le m\le n-1}\sum_{i:\{m\}\sr \{n\}} (-1)^{n-m-1}[i_X\times (*)^{(n-m-1)}]$$
where $i$ runs through all order preserving monomorphisms $\{m\}\sr\{n\}$. For example, for $n=3$ we get 
$$s_3=[pr_{1,2}]+[pr_{2,3}]+[pr_{1,3}]-[pr_1\times *]-[pr_2\times *]-[pr_3\times *]+[*\times *].$$
where $pr_{i_1,\dots,i_m}$ corresponds to the map $i:\{m\}\sr \{n\}$ whose image is $\{i_1,\dots,i_m\}$.  One verifies easily that the maps $s_n$ defined in this way satisfy the two conditions stated above.
\end{proof}
Set $S^{\infty}_{tr}=\oplus_{n\ge 1}S_{tr}^n$. For $[X]\in Cor(C,R)$ the morphisms $\tilde{S}^i(X)\wedge \tilde{S}^j(X)\sr \tilde{S}^{i+j}(X)$ define a morphism
$$S^{\infty}_{tr}([X])\oo S^{\infty}_{tr}([X])\sr S^{\infty}_{tr}([X])$$
which makes $S^{\infty}_{tr}([X])$ into a commutative monoid relative to $\oo$. One verifies immediately that this monoid structure is natural for morphisms in $Cor(C,R)$. 

Combining Proposition \ref{ext} with Proposition \ref{split} and extending them to radditive functors we get the following result. 
\begin{proposition}
\llabel{sinfty}
There is a commutative square of functors
$$
\begin{CD}
Rad(C_+) @>S^{\infty}>> Rad(C_+)\\
@V\Lambda^l_RVV @VV\Lambda^l_R V\\
Rad(Cor(C,R)) @>S^{\infty}_{tr}>> Rad(Cor(C,R))
\end{CD}
$$
\end{proposition}
Note that the functor $S^{\infty}[1/d]$ does not extend to a functor on $Cor$ since the definition of the "multiplication by d" maps $S^n\sr S^{nd}$ involves the diagonals which are not functorial for morphisms in $Cor$. 
\begin{proposition}
\llabel{directsuminfty}
For any $X, Y$ in $Cor(C,R)$ there is an isomorphism
$$S^{\infty}_{tr}(X\oplus Y)=[S^{\infty}_{tr}(X)\oo S^{\infty}_{tr}(Y)]\oplus S^{\infty}_{tr}(X)\oplus S^{\infty}_{tr}(Y)$$
which is natural in $X$ and $Y$.
\end{proposition}
\begin{proof}
We have
$$[S^{\infty}_{tr}(X)\oo S^{\infty}_{tr}(Y)]\oplus S^{\infty}_{tr}(X)\oplus S^{\infty}_{tr}(Y)=\sum_{i,j\ge 0, i+j>0} S^i_{tr}(X)\oo S^j_{tr}(Y)=$$
$$=\sum_{n>0}\sum_{i\ge 0} S^i_{tr}(X)\oo S^{n-i}_{tr}(Y)=S^{\infty}_{tr}(X\oplus Y)$$
where the last equality holds by Proposition \ref{directsums}.  
\end{proof}

\begin{proposition}
\llabel{prestrinfty}
The functor $S^{\infty}_{tr}$ take $(Nis,\af)$-equivalences between objects of $\Delta^{op}Cor^{\#}$ to $(Nis,\af)$-equivalences. In particular there exist unique up to a canonical isomorphism functors 
$$S^{\infty}_{tr}:H_{Nis,\af}(Cor(C,R))\sr H_{Nis,\af}(Cor(C,R))$$
such that the squares
$$
\begin{CD}
\Delta^{op}Cor^{\#} @>S^{\infty}_{tr}>> \Delta^{op}Cor^{\#}\\
@VVV @VVV\\
H_{Nis,\af} @>S^{\infty}_{tr}>> H_{Nis,\af}
\end{CD}
$$
commute. 
\end{proposition}
\begin{proof}
Follows immediately from Theorem \ref{prestr} and the fact that $(Nis,\af)$-equivalences are closed under direct sums.
\end{proof}

\begin{remark}\rm{}\rm
All the results proved in this section for $H_{Nis,\af}(Cor(C,R))$ also hold for the intermediate homotopy categories $H(Cor(C,R))$, $H_{\af}(Cor(C,R))$ and $H_{Nis}(Cor(C,R))$.
\end{remark}\rm{}

\subsection{The  motive of $\tilde{S}^l(T^n)$}
\llabel{TMO}
Everywhere in this section $l$ is a prime not equal to $char(k)$. The underlying category $C$ will be the category of quasi-projective schemes. Let $T^n=a_{+,Nis}({\bf A}^n/({\bf A}^n-\{0\}))$ be the standard model of the motivic n-sphere and
$$L_n=\LL\Lambda^l(T^n)=\Lambda^l(L_*(T^n))$$
its image in $H_{Nis,\af}(Cor(C,R))$.  The goal of this section is to compute the isomorphism class of $\LL S^l_{tr}(L_n)$ in the case when $R={\bf F}_l$.

Recall that for a linear representation $\rho:G\sr Aut(V)$ of a finite group $V$ we let $Th(\rho)$ denote the quotient sheaf (in the Nisnevich topology) 
$$Th(\rho)=a_{+,Nis}(V/V-\{0\})$$
where $V$ and $V-\{0\}$ are considered as representable sheaves on quasi-projective $G$-schemes (see \cite{zslicepub}). Applying the functor $Quot_G$ on sheaves we get a pointed sheaf $Quot_G(Th(\rho))$ on $C_{Nis}$. Since $Quot_G$ commutes with colimits and coincides with the scheme-theoretic quotient on representable sheaves we have
$$Quot_G(Th(\rho))=(V/G)/((V-\{0\})/G).$$

We start by computing the isomorphism class of $\LL\Lambda^lQuot_G(Th(\rho))$ in $H_{Nis,\af}(Cor(C,R))$ in the case when $\rho$ is a representation of the cyclic group $\zz/l$ and $R$ is any commutative ring. Similar computations were done independently by Nie (see \cite{Nie2008}).  In what follows quotients, wedge products etc. are considered in the category of pointed Nisnevich sheaves i.e. after the application of the associated sheaf functor to the result of the corresponding construction in radditive functors. 

Let us say that a linear representation $\rho:\zz/l\sr Aut(V)$ is free if the corresponding action of $\zz/l$ on $V-\{0\}$ is free. Since $char(k)\ne l$, any nontrivial representation $\rho$ has a canonical decomposition into a direct sum
$\rho=\lambda\oplus\tau$ where $\lambda$ is free and $\tau$ is a trivial. As was shown in \cite{zslicepub} one has an isomorphism of sheaves
$$Quot_{\zz/l}(Th(\lambda\oplus\tau))=Quot_{\zz/l}(Th(\lambda))\wedge T^d$$
where $d=dim(\tau)$. Therefore it is sufficient to consider free representations $\rho$. For a linear representation of any $G$ we have
$$Quot_G(Th(\rho))=Quot_G(V)/Quot_G(V-\{0\})$$
because $Quot_G$ is a left adjoint and therefore commutes with colimits (see \cite[\S 5.1]{Delnotessub}). 
Set
$$X_{\rho}=Quot_G(V-\{0\})=(V-\{0\})/G.$$
Since $V$ is $G$-equivariantly $\af$-contractible, the pointed sheaf $Quot_G(V)$ is $\af$-contractible and therefore
$Quot_G(Th(\rho))$ is the unreduced suspension of $X_{\rho}$ i.e. there is a cofiber sequence in $H_{Nis,\af}(C_+)$ of the form
\begin{eq}
\llabel{willsplit2}
(X_{\rho})_+\sr S^0\sr Quot_{G}(Th(\rho))\sr \Sigma^1((X_{\rho})_+).
\end{eq}
\comment{If $dim(V)\ne 0$ there exists a rational point $*$ in $X_{\rho}$ and a choice of such a point defines a splitting of the sequence (\ref{willsplit2}) and therefore an isomorphism in $H_{Nis,\af}(C_+)$ of the form
$$Quot_G(Th(V))=\Sigma^1(X_{\rho},*).$$
Collecting these arguments together we get the following result.
\begin{proposition}
\llabel{tofree}
Let $\lambda:\zz/l\sr Aut(V)$ be a linear representation of the cyclic group and $\lambda={\rho}\oplus{\tau}$ be its decomposition into the direct sum of a free and a trivial representation.  Assume that $dim(\rho)>0$ and let $*$ be a rational point in $X_{\rho}$. Then there is a natural isomorphism in $H_{Nis,\af}(C_+)$ of the form
$$Quot_G(Th(\lambda))=\Sigma^1(X_{\rho},*)\wedge T^d$$
where $d=dim(\tau)$.
\end{proposition}
}
Since we consider representations of a cyclic group, the scheme  $X_{\rho}$ is smooth and the full embedding part of  Proposition \ref{2009.07.24.3} together with Theorem \ref{2009femb}  imply that we may do our computation in the more familiar context of the triangulated category  $DM^{eff}_{-}(k,R)$. To keep in concordance with our earlier notation we will write $\tilde{M}(X)$ instead of $N_{Nis,\af}(X)$. 

Shifting the image of the cofiber sequence (\ref{willsplit2}) in $DM_{-}^{eff}$ by one step to the  left and taking into account the canonical isomorphism $\tilde{M}(S^0)=R$ we get a distinguished triangle 
$$\tilde{M}(Quot_{G}(Th(\rho)))[-1]\sr \tilde{M}((X_{\rho})_+ \sr {\bf F}_l\sr \tilde{M}(Quot_{G}(Th(\rho)))$$
Since $dim(V)>0$, the map $\tilde{M}((X_{\rho})_+ \sr R$ is an epimorphism which is split by any rational point of $X_{\rho}$ and its kernel is canonically isomorphic to the reduced motive $\tilde{M}(X_{rho})$ of the non-pointed  scheme $X_{\rho}$ as defined in \cite[p.192]{H3new}. Therefore we get the following result. 
\begin{proposition}
\llabel{2009.08.10.prop1}
Let $\lambda:\zz/l\sr Aut(V)$ be a non-trivial linear representation of the cyclic group and $\lambda={\rho}\oplus{\tau}$ be its decomposition into the direct sum of a free and a trivial representation. Then $Quot_G(Th(\lambda))\in H_{Nis,\af}((Sm/k)_+)$ and there is a natural isomorphism
$$\tilde{M}(Quot_G(Th(\lambda)))=\tilde{M}(X_{\rho})(d)[2d+1]$$
where $d=dim(\tau)$.
\end{proposition}
Let us assume now that $R={\bf F}_l$. Since the action of $\zz/l$ on $V-\{0\}$ is free the projection $V-\{0\}\sr X_{\rho}$ is an etale Galois covering with the Galois group $\zz/l$ which defines a class $u_{\rho}\in H^1_{et}(X_{\rho},\zz/l)$. This construction is clearly natural i.e. the following lemma holds.
\begin{lemma}
\llabel{nat}
Let $\rho:\zz/l\sr Aut(V)$, $\alpha:\zz/l\sr Aut(W)$ be two free representations and $f:W-\{0\}\sr V-\{0\}$ be an equivariant morphism. Then one has
$$u_{\alpha}=Quot_{\zz/l}(f)^*u_{\rho}.$$
\end{lemma}
To obtain a motivic interpretation of $u_{\rho}$ let us consider the following construction. Let $m_l$ be the object of $DM_{-}^{eff}(k,{\bf F}_l)$ corresponding  the sheaf with transfers ${\bf \mu}_l=ker({\bf G}_m\stackrel{(-)^z}{\sr}{\bf G}_m)$. 
If $k$ contains an $l$-th root of unity then a choice of such a root defines an isomorphism ${\bf F}_l\sr m_l$. In general, $m_l$ is an Artin motive which is a direct summand of the motive of the zero dimensional smooth scheme $Spec\, k[t]/((x^l-1)/(x-1))$. In particular the dual $m_l^*$ to $m_l$ is well defined and there is a canonical isomorphism $m_l^*\oo m_l\sr {\bf F}_l$. The object $m_l^{\oo i}$ corresponds to the etale sheaf $\mu_l^{\oo i}$ and in particular $m_l^{\oo (l-1)}={\bf F}_l$ and $m_l^*=m_l^{\oo (l-2)}$. 
\begin{lemma}
\llabel{2009.08.10.l1}
For any smooth scheme $X$ there is a canonical isomorphism
\begin{eq}
\llabel{2009.08.10.eq1}
H^1_{et}(X,\zz/l)= Hom_{DM}(M(X), m_l^*(1)[1])
\end{eq}
\end{lemma}
\begin{proof}
In $DM_{-,et}^{eff}$ we have ${\bf F}_l(1)[1]=m_l$ since ${\bf G}_m\stackrel{(-)^z}{\sr}{\bf G}_m$ is a surjection in the etale topology. Therefore in $DM_{-,et}^{eff}$ we have a canonical isomorphism $m_l^*(1)[1]={\bf F}_l$. This defines a map from the right hand side of (\ref{2009.08.10.eq1}) to the left hand side. To verify that this map is an isomorphism it is sufficient by the usual transfer argument to do it over an extension of $k$ of degree prime to $l$ which contains an $l$-th root of unity. Then it becomes isomorphic to the well known isomorphism $H^1_{et}(X,{\bf \mu}_l)=H^{1,1}(X,\zz/l)$.
\end{proof}
Using the isomorphism of Lemma \ref{2009.08.10.l1} we may consider $u=u_{\rho}$ as a morphism $M(X_{\rho})\sr m_l^*(1)[1]$. Applying to it an obvious analog of the Bockstein homomorphism we get a morphism 
$$v=\beta(u):M(X_{\rho})\sr  m_l^*(1)[2].$$
Using the tensor structure of $DM$ and the diagonal of $X_{\rho}$ we further define the "product classes"
$$v^i:M(X_{\rho})\sr (m^*_l)^{\oo i}(i)[2i]$$
$$uv^i:M(X_{\rho})\sr (m^*_l)^{\oo (i+1)}(i+1)[2i+1].$$
where $v^0:M(X_{\rho})\sr {\bf F}_l$ is the canonical morphism. 
\begin{proposition}
\llabel{2009.08.10.prop2}
The morphism 
\begin{eq}
\llabel{firstiso}
I(\rho)=\oplus_{i=0}^{n-1}(v^i\oplus uv^i):M(X_{\rho})\sr \bigoplus_{i=0}^{n-1}((m^*_l)^{\oo i}(i)[2i]\oplus (m^*_l)^{\oo (i+1)}(i+1)[2i+1])
\end{eq}
where $n=dim(\rho)$, is an isomorphism which identifies $\tilde{M}(X_{\rho})$ with the sub-object of the right hand side of the form
$$\tilde{M}(X_{\rho})=m^*_l(1)[1]\oplus \bigoplus_{i=1}^{n-1}((m^*_l)^{\oo i}(i)[2i]\oplus (m^*_l)^{\oo (i+1)}(i+1)[2i+1])$$
\end{proposition}
\begin{proof}
By the usual transfer argument we may assume that $k$ contains an $l$-th root of unity $\xi$. 
Then our representation $V$ can be written canonically as a direct sum $\oplus V_m$ where the restriction of $\rho$ to $V_m$ takes $1$ to the multiplication by $\xi^m$. The condition that $\rho$ is free means that $n=dim(\rho)>0$ and $V_0=0$. 

Consider first the case when $V=V_1$. Then \cite[Lemma 6.3]{Redpub} implies that $X_{\rho}$ is canonically isomorphic to the complement to the zero section of the line bundle ${\cal O}(-l)$ on ${\bf P}^{n-1}$ and our result follows easily by computations similar to the one in \cite[pp.18-19]{Redpub}.

Consider now the general $V=\oplus V_m$. Denote by $W$ the same space as $V$ but considered as a $\zz/l$-scheme with respect to the representation $\alpha$ where $\alpha(1)=(w\mapsto \xi w)$. Let us choose a basis $(e_{ij})$ in each $V_j$. Let $p:W\sr V$ be the morphism which takes  $\sum x_{ij}e_{ij}$ to $\sum x_{ij}^{j}e_{ij}$. This morphism is clearly $\zz/l$-equivariant and maps $W-\{0\}$ to $V-\{0\}$. Since $m$ runs from $1$ to $l-1$, the resulting morphism $q:W-\{0\}\sr V-\{0\}$ is finite and surjective and of degree $\prod m^{dim V_m}$ which is prime to $l$. The same is then true for the morphism $p=Quot_{\zz/l}(q):X_{\alpha}\sr X_{\rho}$. Since we work with ${\bf F}_l$ coefficients the morphism of motives $M(p)$ is a split epimorphism. 

Consider now the class $u_{\rho}$. By Lemma \ref{nat} we have $p^*(u_{\rho})=u_{\alpha}$ and therefore we have $I(\rho)\circ M(p)=I(\alpha)$. Since $M(p)$ is a split epimorphism and $I(\alpha)$ is an isomorphism we conclude that both $M(p)$ and $I(\rho)$ are isomorphisms.

The part of the proposition describing $\tilde{M}(X_{\rho})$ follows immediately from the fact that $v^0$ is the canonical morphism $M(X_{\rho})\sr {\bf F}_l$.
\end{proof}
\begin{cor}
\llabel{motcomp}
Let $\lambda:\zz/l\sr Aut(V)$ be a non-trivial linear representation of the cyclic group and $\lambda={\rho}\oplus{\tau}$ be its decomposition into the direct sum of a free and a trivial representation.  Then there is a canonical isomorphism in $H_{Nis,\af}(Cor(C,{\bf F}_l))$ of the form
$$\tilde{M}(Quot_{\zz/l}(Th(\lambda)))=$$
\begin{eq}
\llabel{fifthiso}
=\left(\bigoplus_{i=1}^{n-1}((m_l^*)^{\oo i}(i+d)[2i+2d]\oplus (m_l^*)^{\oo i}(i+d)[2i+2d+1])\right)\oplus (m^*_l)^{\oo n}(d+n)[2d+2n]
\end{eq}
where $n=dim(\rho)$ and $d=dim(\tau)$.
\end{cor}

We will consider now a special case when $\lambda:\zz/l\sr Aut(V)$ is the direct sum of $n$ copies of the regular representation of $\zz/l$. The additional feature which appears in this case is the action of the automorphism group $U=(\zz/l)^*$ of $\zz/l$ on $V$. Let us denote this action by $s:(\zz/l)^*\sr Aut(V)$. It does not commute with the action defined by $\lambda$ i.e. the morphisms $s(m)$ are not $\zz/l$-equivariant but for any $a\in (\zz/l)^*$ the square
\begin{eq}
\llabel{2009.08.10.oldeq1}
\begin{CD}
V-\{0\} @>s(a)>> V-\{0\}\\
@V\lambda(1)VV @VV\lambda(a)V\\
V-\{0\} @>s(a)>> V-\{0\}
\end{CD}
\end{eq}
commutes. Therefore $s$ defines an action of $(\zz/l)^*$ on the related quotients and in particular on $Quot_{\zz/l}(Th(V))$ and we write
$$r(a):Quot_{\zz/l}(Th(V))\sr Quot_{\zz/l}(Th(V))$$
for the automorphism corresponding to $a\in (\zz/l)^*$.
\begin{proposition}
\llabel{onm}
Let $\rho$ be the direct sum of $n>0$ copies of the regular representation of $\zz/l$. Then the isomorphism of Corollary \ref{motcomp} is of the form
$$\tilde{M}(Quot_{\zz/l}(Th(\rho)))=$$
\begin{eq}
\llabel{sixthiso}
\left(\bigoplus_{i=1}^{n(l-1)-1}((m_l^*)^{\oo i}(i+n)[2i+2n]\oplus (m_l^*)^{\oo i}(i+n)[2i+2n+1])\right)\oplus (m^*_l)^{\oo n}(nl)[2nl]
\end{eq}
With respect to this isomorphism the morphism $\tilde{M}(r(a))$ is of the form
\begin{eq}
\llabel{eqact}
M(r(a))=\left(\bigoplus_{i=1}^{n(l-1)-1} (a^{-i}Id \oplus a^{-i}Id)\right)\oplus a^{-n(l-1)}Id 
\end{eq}
\end{proposition}
\begin{proof}
In the decomposition  
\begin{eq}
\llabel{decomp}
V_{\rho}=V_{\lambda}\oplus V_{\tau}
\end{eq}
of $V$ into the free and trivial parts we have $dim(V_{\lambda})=n(l-1)$ and $dim(V_{\tau})=n$ which implies the first part of the proposition. To prove the second part observe first that the decomposition (\ref{decomp}) is invariant under the action of $(\zz/l)^*$ and that the corresponding action on $V_{\tau}$ is trivial. Therefore, the action of $(\zz/l)^*$ on (\ref{sixthiso})
is determined by its action on $M(X_{\lambda})$. The action of any endomorphism of $X_{\lambda}$ on its motive is determined by its action on the motivic cohomology class $u_{\lambda}$.   

In view of (\ref{2009.08.10.oldeq1}) we may consider $s(a)$ as an equivariant morphism assuming that the action on the first copy of $V_{\lambda}-\{0\}$ is given by $\lambda$ and on the second copy by $\lambda_a$ where $\lambda_a(1)=\lambda(a)$. Applying Lemma \ref{nat} we conclude that $s(a)^*(u_{\lambda})=a^{-1}\cdot u_{\lambda}$ which implies our result.
\end{proof}
We are ready now to prove the main theorem of this section. Recall that we let $L_n$ denote the image $\LL\Lambda^l(T^n)$ of $T^n$ in $H_{Nis,\af}(Cor(C,R))$.
\begin{theorem}
\llabel{tonextnew}
Let $l$ be a prime and $k$ be a perfect field of characteristic $\ne l$. Let $C$ be the category of quasi-projective schemes over $k$. Then for any $n>0$ there is a natural  isomorphism in $H_{Nis,\af}(Cor(C, {\bf F}_l))$ of the form
$$\LL S^l_{tr}(L_n)=L_{ln}\oplus\bigoplus_{i=1}^{n-1}(L_{i(l-1)+n}\oplus \Sigma^1_sL_{i(l-1)+n}).$$
\end{theorem}
\begin{proof}
By our definition of $L_n$ we have $L_n\in \Delta^{op}Cor(C,R)$ and in particular $\LL S^l_{tr}(L_n)=S^l_{tr}(L_n)$. By \cite[Corollary 15.8]{MVW}, $L_n$ is an even object and therefore Proposition \ref{2009.08.10.prop3} applies i.e.
$$S^l_{tr}(L_n)=S^{\Psi}_{tr}(L_n)/Aut(\zz/l)=S^{\Psi}_{tr}(\Lambda^l(L_*(T^n)))/Aut(\zz/l)=$$$$=\Lambda^l\tilde{S}^{\Psi}(L_*(T^n))/Aut(\zz/l)=\LL \Lambda^l\tilde{S}^{\Psi}(L_*(T^n))/Aut(\zz/l)$$
where $\Psi=(\zz/l, \zz/l\sr S_l)$ is the permutation group responsible for the cyclic product. Since $T^n$ is a solid sheaf we further have by Theorem \ref{mainpreserve} a $Nis$-equivalence 
$$\LL \Lambda^l\tilde{S}^{\Psi}(L_*(T^n))/Aut(\zz/l)=\LL\Lambda^l \tilde{S}^{\Psi}(T^n)/Aut(\zz/l)$$
By Proposition \ref{forcompl} we have (after taking the associated sheaves)
$$\tilde{S}^{\Psi}(T^n)=\tilde{S}^{\Psi}({\bf A}^n_+)/(\tilde{S}^{\Psi}({\bf A}^n_+)-\tilde{S}^{\Psi}(\{0\}_+))=Quot_{\zz/l}Th(\rho)$$
where $\rho$ is the direct sum of $n$ copies of the regular representation of $\zz/l$. It remains to apply Proposition \ref{onm} taking into account that $((m_l^*)^{\oo i})/(Aut(\zz/l))=0$ if $i\ne 0(mod (l-1))$ and $((m_l^*)^{\oo i})/(Aut(\zz/l))=\zz/l$ otherwise and Theorem \ref{2009femb}. 
\end{proof}

\begin{example}\rm
\llabel{2009.09.15.ex1}
Let $l_n=\LL\Lambda_F^l(S^n_t)$ be the unstable version of $F(n)[n]$. Then $\Sigma^nl_n=L_n$ and $l_n$ is an odd or an even object depending on the parity of $n$. Therefore we may use Theorem \ref{tonextnew} together with Corollary \ref{homexact} to try to compute $\LL S^l_{tr}(l_n)$. One can see immediately that there exists $N$ such that $\Sigma^N\LL S^l_{tr}(l_n)$ can be obtained from elementary Tate object $l_i$ with $i\le l\cdot n$ by taking cones of morphisms i.e. that $\Sigma^N\LL S^l_{tr}(l_n)$ is $s$-stably a mixed Tate object. Let us compute $\Sigma^4\LL S^l_{tr}(l_2)$ over an algebraically closed field of characteristic zero more explicitly. 

 Let $X=\Sigma^1l_2$. Then $X$ is odd and applying to it Corollary \ref{homexact} we have a cofiber sequence
$$\Sigma^{l-1}X^{\oo l}\sr\Sigma^1\LL S^l_{tr}(X)\sr \LL S^l_{tr}(\Sigma^1 X)\sr \Sigma^{l}X^{\oo l}.$$
Since $\Sigma^1X=L_2$ we can rewrite it as
\begin{eq}
\llabel{minus1eq}
\Sigma^{2l-1}l_{2l}\sr \Sigma^1\LL S^l_{tr}(X)\sr  \LL S^l_{tr}(L_2)\sr \Sigma^{2l}l_{2l}
\end{eq}
Theorem \ref{tonextnew} shows that in the motivic notation the sequence (\ref{minus1eq})   is of the form
$${\bf F}_l(2l)[4l-1]\sr \Sigma^1\LL S^l_{tr}(X)\sr {\bf F}_l(2l)[4l]\oplus {\bf F}_l(l+1)[2l+2]\oplus {\bf F}_l(l+1)[2l+3]\stackrel{\partial_1}{\sr} {\bf F}_l(2l)[4l]$$
for some morphism $\partial_1$. 

Since the topological realization of $\Sigma^1S^2_t$ is the 3-sphere we know that the ordinary homology of $\LL S^l_{tr}(X)$ are of the form $H_{2l+1}=H_{2l+2}={\bf F}_l$ and the rest of the homology groups are zero (\cite{Nakaoka}).  Together with the properties of the topological realization functor considered in Section \ref{top.real} this implies that for $l>2$ the morphism $\partial_1$ is an isomorphism on the first summand and zero on two other summands and therefore 
\begin{eq}
\llabel{minus12}
\Sigma^2\LL S^l_{tr}(X)\cong  {\bf F}_l(l+1)[2l+3]\oplus {\bf F}_l(l+1)[2l+4].
\end{eq}
Applying Corollary \ref{homexact} to $l_2$ we get a cofiber sequence 
$$\Sigma^1 l_{2l}\sr \Sigma^1 \LL S^l_{tr}(l_2)\sr \LL S^l_{tr}(X)\sr \Sigma^2 l_{2l}$$
or in the motivic notation
$${\bf F}_l(2l)[2l+1]\sr \Sigma^1 \LL S^l_{tr}(l_2)\sr S^l_{tr}(X)\sr {\bf F}_l(2l)[2l+2].$$
Suspending it twice and using (\ref{minus12}) we get
$${\bf F}_l(2l)[2l+3]\sr \Sigma^3 \LL S^l_{tr}(l_2)\sr {\bf F}_l(l+1)[2l+3]\oplus {\bf F}_l(l+1)[2l+4]\stackrel{\partial_2}{\sr} {\bf F}_l(2l)[2l+4].$$
for some morphism $\partial_2$.   The restriction of $\partial_2$ to the first summand belongs to the group which is isomorphic to $H^{1,l-1}(k,{\bf F}_l)$ and which is  zero since the base field is algebraically closed (see \cite{}). It remains to compute the restriction of $\partial_2$ to the second summand. Topologically, $S^2_t$ is a 2-sphere and applying again the topological realization functor we conclude that $\LL S^l_{tr}(l_2)$ has only one non-trivial ordinary homology group in dimension $2l$.  Therefore, the restriction of $\partial_2$ to the second summand is a morphism of the form
$$\tau':{\bf F}_l(l+1)[2l+4]\sr {\bf F}_l(2l)[2l+4]$$
which defines an isomorphism on ordinary homology. Up to the multiplication by an element of $(\zz/l)^*$ there is a unique such morphism which corresponds to the generator $\tau_{l-1}$ of $H^{0,l-1}(Spec(k),\zz/l)$ and we conclude that there is an isomorphism of the form
$$\Sigma^4\LL S^l_{tr}(l_2)\cong {\bf F}_l(l+1)[2l+4]\oplus cone(\tau')$$
The topological realization of $cone(\tau')$ is trivial and does not affect the ordinary homology of $\LL S^l_{tr}(l_2)$ but the object itself is non-trivial. In particular, neither $\LL S^l_{tr}(l_2)$ nor any of its suspensions is a direct sum of elementary Tate objects. 
\end{example}

\subsection{Split proper Tate objects}
\llabel{PTO}
In this section we assume that the coefficient ring used to define finite correspondences is a field $F$ and that the base field $k$ is perfect.  
\begin{definition}
\llabel{longname}
An object $X$ in $H_{Nis,\af}(Cor(C,F))$ is called a split proper Tate object if it is isomorphic to a coproduct (direct sum) of objects of the form $\Sigma^iL_j$ for $i\ge 0$.
\end{definition}
We denote the full subcategory of split proper Tate objects by $\overline{SPT}$. Let further $\overline{SPT}_n$ (resp. $\overline{SPT}_{\le n}$, $\overline{SPT}_{\ge n}$) be the full subcategory in $\overline{SPT}$ which consists of direct sums of objects of the form $\Sigma^k L_n$ (resp. $\Sigma^k L_m$ for $m\le n$, $m\ge n$). All these subcategories are clearly closed under direct sums and tensor products.  

Since objects $L_j$ belong to $H_{Nis,\af}(Cor(Sm/k,F))$, Theorem \ref{2009femb} is applicable and we may consider $\overline{SPT}$ as a subcategory of $DM_{-}^{eff}$. In the standard $DM$ notation this subcategory consists of direct sums of objects of the form $F(n)[m]$ with $m\ge 2n$. 
We will also work with the subcategory $\overline{SDT}$ (resp. $\overline{SDT}_{\le n}$ etc.) in $DM_{-}^{eff}$ which consists of all direct sums of elementary Tate objects $F(n)[m]$ for $n\ge 0$ and $m\in\zz$ which exist in $DM_{-}$. 

Let us recall the following key result.
\begin{theorem}[see \cite{cancellationsub}]
\llabel{2009.08.13.th1}
Let $k$ be a perfect field. Then for any $Y\in DM_{-}^{eff}(k,F)$ we have:
$$Hom_{DM}(Y(i)[j],F(i')[j'])=\left\{
\begin{array}{ll}
Hom_{DM}(Y,F(i'-i)[j'-j])&\mbox{\rm for $i'\ge i$}\\
0&\mbox{\rm for $i'<i$}
\end{array}
\right.
$$
\end{theorem}
\begin{lemma}
\llabel{2009.08.13.l1}
For any $n\ge 0$ the subcategory $\overline{SDT}_n$ is abelian, semi-simple and closed under cones i.e. if in a distinguished triangle $X\sr Y\sr Z\sr X[1]$ one has $X,Y\in \overline{SDT}_n$ then $Z\in \overline{SDT}_n$ and the same holds for $\overline{SPT}_n$. 
\end{lemma}
\begin{proof}
By Theorem \ref{2009.08.13.th1} together with the fact that
$$
Hom_{DM}(F,F[i])=H^{p,0}(Spec(k),F)=\left\{
\begin{array}{ll}
F&\mbox{\rm for $i=0$}\\
0&\mbox{\rm for $i\ne 0$}
\end{array}
\right.
$$
we conclude that for any $n$
$$
Hom_{DM}(F(n)[i],F(n)[j])=\left\{
\begin{array}{ll}
F&\mbox{\rm for $i=j$}\\
0&\mbox{\rm for $i\ne j$}
\end{array}
\right.
$$
Therefore, the functor
$$M\mapsto \oplus_{i}Hom(F(n)[i],M)$$
defines an equivalence between $\overline{SDT}_n$ and the category of graded vector spaces $(V_i)_{i\in\zz}$ over $F$ such that $V_i=0$ for $i\ll 0$ which maps $\overline{SPT}_n$ to the subcategory of spaces such that $V_i=0$ for $i<2n$.  This shows that both categories are abelian and semi-simple.  Let $X\stackrel{f}{\sr} Y\sr Z\sr X[1]$ be a distinguished triangle with $X,Y\in \overline{SDT}_n$. Since $\overline{SDT}_n$ is abelian and semi-simple,  the morphism $f$ is isomorphic to a morphism of the form 
$$ker(f)\oplus Im(f)\stackrel{0\oplus Id}{\sr} coker(f)\oplus Im(f)$$
Since the cone of a direct sum of two morphisms is the direct sum of cones it remains to verify that $Z\in \overline{SDT}_n$ (resp. $Z\in \overline{SDT}_n$) if $X,Y\in \overline{SDT}_n$ (resp. $X,Y\in \overline{SPT}_n$) and $f=0$ or $f=Id$, which is obvious. 
\end{proof}
The following construction which we give here in the context of $\overline{SDT}$ is a particular case of  the general slice filtration constructions and can be extended to a much wider class of motives (see e.g. \cite{oversub}).  

For $X\in DM_{-}^{eff}$ and $n\in \zz$ consider the sub-functor of the functor representable by $X$ on $DM_{-}^{eff}$ which consists of morphisms $Y\sr X$ which admit a factorization of the form $Y\sr Z(n+1)\sr X$. Theorem \ref{2009.08.13.th1} implies immediately that  for $X=\oplus_{i\ge 0,j\in\zz}F(i)[j]^{\oplus n(i,j)}\in \overline{SDT}$ this sub-functor  is represented by a direct summand $s_{>n}(X)$ of $X$ which is identified with $\oplus_{i>n,j\in\zz}F(i)[j]^{\oplus n(i,j)}$. Set $s_{\le n}(X)=X/s_{>n}(X)$ and 
$$s_n(X)=s_{>n-1}(X)/s_{>n}(X)=ker(s_{\le n}X\sr s_{\le n-1}(X)).$$
This construction provides for each $X\in \overline{SDT}$ a collection of split distinguished triangles of the form
\begin{eq}
\llabel{gr1}
s_{>n}X\sr X\stackrel{p_n}{\sr} s_{\le n}X\sr s_{>n}X[1]
\end{eq}
and
\begin{eq}
\llabel{gr2}
s_nX\sr s_{\le n}X\stackrel{q_n}{\sr} s_{\le n-1}X\sr s_nX[1]
\end{eq}
whose terms are in $\overline{SDT}$ and which are natural in $X$ and commute in the obvious sense with the shift functor. 
\begin{remark}\rm
Note that for $X=\oplus_{i,j}F(i)[j]^{\oplus m(i,j)}$ we have 
\begin{enumerate}
\item $s_{>n}X=\oplus_{i>n,j}F(i)[j]^{\oplus m(i,j)}$ as a subobject of $X$,
\item $s_{\le n}X=\oplus_{i\le n,j}F(i)[j]^{\oplus m(i,j)}$ as a quotient object of $X$. 
\end{enumerate}
We also have isomorphisms $s_nX\cong \oplus_{j}F(m)[j]^{\oplus m(n,j)}$ and  in particular there is an isomorphism $X\cong \oplus_n s_nX$. However, this last isomorphism is not natural in $X$.
\end{remark}

\begin{proposition}
\llabel{2009.09.18.pr1}
Let $X=\oplus_{\alpha\in A} F(q_{\alpha})[p_{\alpha}]$ be an object of $\overline{SPT}$ such that for each $q$ the set $A_{*,\le q}=\{\alpha\in A\,|\,q_{\alpha}\le q\}$ is finite. Then $X$ is the direct product of the family $(F(q_{\alpha})[p_{\alpha}])_{\alpha\in A}$ in $DM^{eff}_{-}$.
\end{proposition}
\begin{proof}
Since direct product is an exact functor from families of $F$-vector spaces to $F$-vector spaces the functor $\prod_{\alpha\in A}Hom(-,F(q_{\alpha})[p_{\alpha}])$ is a cohomological functor which takes direct sums to direct products.  Therefore in order to show that the natural transformation 
$$Hom(-,X)\sr \prod_{\alpha\in A}Hom(-,F(q_{\alpha})[p_{\alpha}])$$
is an isomorphism it is sufficient to verify that it defines isomorphisms on a set of generators of $DM^{eff}_{-}$ for which we can take the set of objects of the form $M(U)[i]$ for $U\in Sm/k$ and $i\in\zz$. Since these objects are compact we need to show that for such $U$ and $i$ the natural map
\begin{eq}
\llabel{2009.09.18.eq1}
\oplus_{\alpha}Hom(M(U)[i],F(q_{\alpha})[p_{\alpha}])\sr \prod_{\alpha}Hom(M(U)[i],F(q_{\alpha})[p_{\alpha}])
\end{eq}
Since $X\in \overline{SPT}$ we have $p_{\alpha}\ge 2q_{\alpha}$ and therefore 
$$Hom(M(U)[i],F(q_{\alpha})[p_{\alpha}])=0$$
for $q_{\alpha}\ge dim(U)+i+1$. Together with our assumption on the finiteness of the sets $A_{*,\le q}$ this shows that the family $Hom(M(U)[i],F(q_{\alpha})[p_{\alpha}])$ has only finitely many non-zero members and therefore the map (\ref{2009.09.18.eq1}) is a bijection.
\end{proof}
\begin{lemma}
\llabel{fromsq}
Consider a distinguished triangle in $DM_{-}^{eff}(k,F)$ of the form
\begin{eq}
\llabel{sqfor}
X\sr Y\sr Z\sr X[1]
\end{eq}
Such that $X\in \overline{SDT}_{\le n}$ and $Z\in \overline{SDT}_n$. Then $Y\in \overline{SDT}_{\le n}$. Similarly if $X\in \overline{SPT}_{\le n}$ and $Z\in \overline{SPT}_n$  then $Y\in \overline{SPT}_{\le n}$.
\end{lemma}
\begin{proof}
Since $Z\in \overline{SDT}_n$, the morphism $Z\sr X[1]$ factors through a morphism $Z\sr s_{>(n-1)}X[1]$ and we get a diagram
$$
\begin{CD}
Z @>>> s_{>(n-1)}X[1] @.\\
@V=VV @VVV @.\\
Z @>>> X[1]  @>>> Y[1] @>>> Z[1]
\end{CD}
$$
where lower row is the shift of our original triangle by one step to the right.  By the usual properties of triangulated categories (see e.g. \cite[Prop. 1.4.6, p.58]{Nee}) this diagram can be extended 
to a commutative diagram of the form
$$
\begin{CD}
Z @>>> s_{>(n-1)}X[1] @>>> W @>>> Z[1]\\
@V=VV @VVV @VVV @VVV\\
Z @>>> X[1]  @>>> Y[1] @>>> Z[1]\\
@VVV @VVV @VVV @VVV\\
0 @>>> s_{\le (n-1)}X[1] @>f_2=Id>>  s_{\le (n-1)}X[1]  @>>> 0\\
@VVV @Vf_3VV @VVf_1V @VVV\\
Z[1] @>>> s_{>(n-1)}X[2] @>>> W[1] @>>> Z[2]
\end{CD}
$$
whose rows and columns are distinguished triangles. Since $X\in \overline{SDT}_{\le n}$ we have $s_{>(n-1)}(X)=s_n(X)\in \overline{SDT}_n$ and by Lemma \ref{2009.08.13.l1} we conclude that $W\in \overline{SDT}_n$. Moreover, it is easy to see from the upper distinguished triangle of our diagram that if $X,Z\in \overline{SPT}$ then $W[-1]\in \overline{SPT}_n$.  To prove the lemma it it is now sufficient to show that $f_1=0$ i.e. $Y[1]=W\oplus s_{\le (n-1)}(X)[1]$. This follows immediately from the fact that $f_2$ is an isomorphism and $f_3=0$ since triangles of the form (\ref{gr1}) for objects of $\overline{SDT}$ are split. 
\end{proof}
\begin{lemma}
\llabel{2009.09.04.l1}
Let $X\stackrel{f}{\sr} Y\stackrel{g}{\sr} Z$ be a sequence of morphisms in $\overline{SDT}$ such that for each $U\in Sm/k$, $i\in\zz$ the sequence of abelian groups
\begin{eq}
\llabel{2009.09.04.eq1}
0\sr Hom(M(U)[i],X)\sr Hom(M(U)[i],Y)\sr Hom(M(U)[i],Z)\sr 0
\end{eq}
is exact. Then the sequence $0\sr X\stackrel{f}{\sr} Y\stackrel{g}{\sr} Z\sr 0$ is split exact in $DM_{-}^{eff}(k,R)$.
\end{lemma}
\begin{proof}
Consider a distinguished triangle 
\begin{eq}
\llabel{2009.09.04.eq2}
X\stackrel{f}{\sr}Y\stackrel{g''}{\sr} cone(f)\sr X[1]
\end{eq}
Since $X$ is a direct sum of Tate objects and Tate objects are direct summands of shifts of objects of the form $M({\bf P}^n)$ the exactness of (\ref{2009.09.04.eq1}) implies that $g\circ f=0$ and therefore there exists a morphism $cone(f)\stackrel{g''}{\sr} Z$ such that $g=g''g'$. Since the maps $Hom(M(U)[i],X)\sr Hom(M(U)[i],Y)$ are injective the long exact sequence defined by the triangle (\ref{2009.09.04.eq2}) splits into short exact sequences of the form
$$0\sr Hom(M(U)[i],X)\sr Hom(M(U)[i],Y)\sr Hom(M(U)[i],cone(f))\sr 0$$
and we conclude that the map $g''$ defines isomorphisms on $Hom(M(U)[i],-)$ for all $U\in Sm/k$ and $i\in\zz$. Since $DM_{-}^{eff}(k,R)$ is generated as a triangulated category by objects of the form $M(U)$ for $U\in Sm/k$ we conclude that $g''$ is an isomorphism and $X\stackrel{f}{\sr} Y\stackrel{g}{\sr} Z$ extends to a distinguished triangle $X\stackrel{f}{\sr} Y\stackrel{g}{\sr} Z\stackrel{\partial}{\sr}X[1]$. Since sequences (\ref{2009.09.04.eq1}) are exact the morphism $\partial$ is zero on $Hom(M(U)[i],-)$ for all $U$ and $i$. Since $Z$ is a direct sum of Tate objects we conclude that $\partial=0$. This is equivalent to the assertion of the lemma. 
\end{proof}
\begin{lemma}
\llabel{2009.09.04.l2}
Let $X\stackrel{f}{\sr} Y\stackrel{g}{\sr} Z$ be a sequence of morphisms in $\overline{SPT}$ such that for all $n\ge 0$ the sequence
$$0\sr s_{\le n}X\sr s_{\le n}Y\sr s_{\le n}Z\sr 0$$
is split exact. Then the sequence $0\sr X\stackrel{f}{\sr} Y\stackrel{g}{\sr} Z\sr 0$ is split exact.
\end{lemma}
\begin{proof}
For any $U\in Sm/k$ we have $Hom(M(U),F(p)[q])=0$ for $q>p+dim(U)$. Since for a proper split Tate object $W$, $s_{>n}W$ is a direct sum of copies of object $F(p)[q]$ for $p>n$ and $q\ge 2p$ we have $Hom(M(U)[i],s_{>n}W)=0$ for all $n>i+dim(U)-1$. Therefore for each $U$ and $i$ there exists $n$ such that $Hom(M(U)[i],W)=Hom(M(U)[i],s_{\le n}W)$ for $W=X,Y,Z$. We conclude that the sequence  $X\stackrel{f}{\sr} Y\stackrel{g}{\sr} Z$ satisfies the condition of Lemma \ref{2009.09.04.l1} and therefore it is split exact. 
\end{proof}
\begin{proposition}
\llabel{splitmono}
For a morphism $f:X\sr Y$ in $DM_{-}^{eff}(k,F)$ one has:
\begin{enumerate}
\item if $X,Y\in \overline{SDT}_{\le n}$ and $s_i(f)$ is a split monomorphism for each $i\le n$ then $f$ a split monomorphism and $cone(f)\in \overline{SDT}_{\le n}$,
\item if $X,Y\in \overline{SPT}$ and $s_i(f)$ is a split monomorphism for each $i\ge 0$ then $f$ is a split monomorphism and $cone(f)\in \overline{SPT}$.
\end{enumerate}
\end{proposition}
\begin{proof}
To prove the first assertion let us shows that for each $m\le n$ there is a split exact sequence of the form 
\begin{eq}
\llabel{2009.09.04.eq4}
0\sr s_{\le m}X\stackrel{s_{\le m}f}{\sr}s_{\le m}Y\stackrel{\phi_{\le m}}{\sr} \oplus_{i\le m}(s_iY/s_iX)\sr 0
\end{eq}
where $s_iY/s_iX$ is the cokernel of $s_i(f)$ which is well defined since $s_i(f)$ is a split monomorphism. 

We proceed by induction on $m$. For $m=0$ the statement is obvious. To make an inductive step observe that by simple diagram search there exists a morphism $\phi_{\le m}$ which fits into the commutative diagram of the form
$$
\begin{CD}
s_mX @>s_m(f)>> s_mY @>>> s_mY/s_mX\\
@VVV @VVV @VVV\\
s_{\le m}X @>s_{\le m}(f)>> s_{\le m}Y @>\phi_{\le m}>> \oplus_{i\le m}(s_iY/s_iX)\\
@VVV @VVV @VVV\\
s_{\le (m-1)}X @>s_{\le m}(f)>> s_{\le (m-1)}Y @>\phi_{\le (m-1)}>> \oplus_{i\le (m-1)}(s_iY/s_iX)
\end{CD}
$$
and that for any such morphism the middle row is split exact. 

Consider now the second assertion. By the same argument as above we see that there are split exact sequences of the form (\ref{2009.09.04.eq4}). Since $Y$ is a direct sum of Tate objects and $s_iY/s_iX$ are proper split Tate objects  a simple argument shows that, for each choice of morphisms $\phi_{\le m}$ as above, there exists a unique morphism $\phi:Y\sr \oplus_{i\ge 0} (s_iY/s_iX)$ such that the composition of $\phi$ with the projection to $\oplus_{i\le m} (s_iY/s_iX)$ equals $\phi_{\le m}$. The sequence $X\stackrel{f}{\sr}Y\stackrel{\phi}{\sr}\oplus_{i\ge 0} (s_iY/s_iX)$ clearly satisfies the condition of Lemma \ref{2009.09.04.l2} and therefore is split exact. Proposition is proved.
\end{proof}
\begin{cor}
\llabel{2009.09.17.cor2}
For a morphism $f:X\sr Y$ in $DM_{-}^{eff}(k,F)$ one has:
\begin{enumerate}
\item if $X,Y\in \overline{SDT}_{\le n}$ and $s_i(f)$ is an isomorphism for each $i\le n$ then $f$ an isomorphism,
\item if $X,Y\in \overline{SPT}$ and $s_i(f)$ is an isomorphism for each $i\ge 0$ then $f$ is an isomorphism.
\end{enumerate}
\end{cor}
\begin{cor}
\llabel{2009.09.17.cor3}
Let $p:Spec(K)\sr Spec(k)$ be a field extension and let $f:X\sr Y$ be a morphism in $DM_{-}^{eff}(k,F)$. Then one has:
\begin{enumerate}
\item if $X,Y\in \overline{SDT}_{\le n}$ and $p^*(f)$ is an isomorphism then $f$ an isomorphism,
\item if $X,Y\in \overline{SPT}$ and $p^*(f)$ is an isomorphism then $f$ is an isomorphism.
\end{enumerate}
\end{cor}
\begin{proof}
The functors $s_i:\overline{SDT}\sr \overline{SDT}_i$ commute with the functor $p^*$. On the other hand the restriction of $p^*$ to $\overline{SDT}_i$ is an equivalence in view of the proof of Lemma \ref{2009.08.13.l1}. Therefore, for a morphism $f$ satisfying the conditions of the corollary, the morphisms $s_i(f)$ are isomorphisms and we conclude that $f$ is an isomorphism by Corollary \ref{2009.09.17.cor2}.
\end{proof}

Let $X_1\stackrel{f_1}{\sr}\dots\stackrel{f_{i-1}}{\sr} X_i\stackrel{f_i}{\sr}\dots$  be a sequence of morphisms in a triangulated category such that $\oplus_{i\ge 1}X_i$ exists. Recall that the homotopy colimit of this sequence is an object $hocolim_i X_i$ which is defined up to a non-canonical isomorphism by a distinguished triangle of the form
\begin{eq}
\llabel{2009.09.04.eq3}
\oplus_{i\ge 1}X_i\stackrel{\psi(f_1,\dots)}{\sr}\oplus_{i\ge 1}X_i\stackrel{p}{\sr}hocolim_i X_i\sr 
\oplus_{i\ge 1}X_i[1]
\end{eq}
where the morphism $\psi(f_1,\dots)$ is defined by the condition that its restriction to $X_i$ is $\iota_i-\iota_{i+1}f_i$ where $\iota_i:X_i\sr \oplus_{i\ge 1}X_i$ is the canonical embedding. 
\begin{cor}
\llabel{hocolim}
Let $X_1\stackrel{f_1}{\sr}\dots\stackrel{f_{i-1}}{\sr} X_n\stackrel{f_i}{\sr}\dots$ be a sequence of morphisms in $\overline{SPT}$. Then the direct sum $\oplus_{i\ge 1}X_i$ exists, the distinguished triangle (\ref{2009.09.04.eq3}) which defines $hocolim_i X_i$ splits and $hocolim_i X_i\in \overline{SPT}$.
\end{cor}
\begin{proof}
The direct sum $\oplus_{i\ge 1}X_i$ exists because objects $X_i$ belong to the image of the category $H_{Nis,\af}(Cor(Sm/k,F))$ which has all direct sums under the functor which respects direct sums. We obviously have $s_n(\psi(f_1,\dots))=\psi(s_n(f_1),\dots)$ for each $n\ge 0$ and since $\overline{SDT}_n$ is equivalent to the category of graded $F$-vector spaces one observes easily that the morphisms $s_n(\psi(f_1,\dots))$ are split monomorphisms. Our result follows now from Proposition \ref{splitmono}(2).
\end{proof}
\begin{cor}
\llabel{hocolim2}
Let $X_1\stackrel{f_1}{\sr}\dots\stackrel{f_{i-1}}{\sr} X_n\stackrel{f_i}{\sr}\dots$ be a sequence of morphisms in $\overline{SDT}_{\le n}$ such that $\oplus_{i\ge 1}X_i$ exists. Then the distinguished triangle (\ref{2009.09.04.eq3}) which defines $hocolim_i X_i$ splits and $hocolim_i X_i\in \overline{SDT}_{\le n}$.
\end{cor}
\begin{proof}
Same argument as in the proof of Corollary \ref{hocolim}.
\end{proof}
\begin{cor}
\llabel{allcl}
The subcategory $\overline{SPT}$ (resp. $\overline{SDT}_{\le n}$) is closed under direct summands i.e. contains images of projectors. 
\end{cor}
\begin{proof}
The image of a projector $p$ can be identified with the homotopy colimit of the sequence $X\stackrel{p}{\sr}\dots\stackrel{p}{\sr}X\stackrel{p}{\sr}\dots$. Since for any $X$ in $DM_{-}^{eff}$ the direct sum $\oplus_{i\ge 1} X$ of countably many copies of $X$ exists the result follows from Corollary \ref{hocolim} in the case of $\overline{SPT}$ and Corollary \ref{hocolim2} in the case of $\overline{SDT}_{\le n}$. 
\end{proof}
\begin{remark}\rm
The analog of Corollary \ref{hocolim2} (and therefore of Proposition \ref{splitmono}) for the whole category $\overline{SDT}$ is false. It is easy to see on the example of the sequence $X_i\stackrel{f_i}{\sr}X_{i+1}$ where $X_i=F(i-1)[i-1]$, $f_i=Id\oo \rho$ where $\rho\in H^{1,1}(k,F)$ is the class of $-1$, $k={\bf R}$ and $F={\bf F}_2$. The corresponding homotopy colimit $Y=hocolim_i X_i$ has the property that $s_n(Y)=0$ for all $n\ge 0$ but $Y\ne 0$.

I do not know whether or not the analog of Corollary \ref{allcl} holds for $\overline{SDT}$. 
\end{remark}

Our next goal is to prove  Theorem \ref{main2} which shows that $\overline{SPT}$ is closed under (standard) symmetric powers. Note that symmetric powers have only been defined on $H_{Nis,\af}(Cor(C,R))$ where $C$ is an $f$-admissible category  and $c(k)$ is invertible in $R$. Therefore, for the purposes of  Theorem \ref{main2} we must consider $\overline{SPT}$ as a full subcategory in $H_{Nis,\af}(Cor(C,F))$ where $C$ is the category of quasi-projective schemes, $F$ is a field and if $char(k)>0$ then $char(F)\ne char(k)$.
\begin{lemma}
\llabel{lcase}
Let $F$ be a field of characteristic $l>0$. Then for any $q\ge 1$ one has $\LL S^l_{tr}(\overline{SPT}_{\ge q})\subset \overline{SPT}_{\ge q+l-1}$.
\end{lemma}
\begin{proof}
For $i<l$ and $X\in \overline{SPT}_{\ge q}$ one has $\LL S^i(X)\in \overline{SPT}_{\ge iq}$. Together with Proposition \ref{directsums} this implies that the class of $X\in \overline{SPT}_{\ge q}$ for which the lemma holds is closed under direct sums.

It remains to show that for any $k\ge 0$ and $q\ge 1$ one has $\LL S^l_{tr}(\Sigma^kL^q)\in \overline{SPT}_{\ge q+l-1}$.  The proof is by induction on $k$. For $k=0$ the result follows from Theorem \ref{tonextnew}. By \cite[Corollary 15.8]{MVW}, $L$ is an even object and therefore $\Sigma^kL^q$ is odd or even depending on whether $k$ is odd or even. The inductive step follows now from Corollary \ref{homexact} and Lemma \ref{fromsq}.
\end{proof}
The present formulation of the following theorem is partly based on the considerations of \cite{patching}.
\begin{theorem}
\llabel{main2}
Let $k$ be our base field and $F$ the field of coefficients. Then for any integers $q\ge 1$, $n\ge 0$ one has:
\begin{enumerate}
\item if $char(F)=0$ then 
$$\LL S^n_{tr}(\overline{SPT}_{\ge q})\subset \overline{SPT}_{\ge nq}$$
\item if $char(F)=l>0$, $char(F)\ne char(k)$ and 
$$n=\sum_{i\ge 0} n_il^i$$
is the $l$-primary decomposition of $n$ then 
$$\LL S^n_{tr}(\overline{SPT}_{\ge q})\subset \overline{SPT}_{\ge w_l(q,n)}$$
where $w_l(q,n)=(\sum_i n_i)q+(\sum_i in_i)(l-1)$.
\end{enumerate}
\end{theorem}
\begin{proof}
The case of $char(F)=0$ is obvious. Suppose that $char(F)=l>0$. Let $X\in \overline{SPT}_{\ge q}$. By Proposition \ref{directsummand}, $\LL S^n_{tr}(X)$ is a direct summand of $\oo_i ((\LL S^l_{tr})^{\circ i}(X))^{\oo n_i}$ where as before $n=\sum_i n_i l^i$ is the $l$-primary decomposition of $n$. By Corollary \ref{allcl} it is sufficient to show that 
$$\oo_i ((\LL S^l_{tr})^{\circ i}(X))^{\oo n_i}\in \overline{SPT}_{\ge w_l(q,n)}$$
Since $w_l(q,n)=\sum_i n_i(q+(l-1)i)$ it is further sufficient to show that
$$(\LL S^l_{tr})^{\circ i}(X)\in \overline{SPT}_{\ge q+i(l-1)}$$
This follows by obvious induction from Lemma \ref{lcase}.
\end{proof}
\begin{cor}
\llabel{sinftyin}
Under the assumptions of the theorem one has $\LL S^{\infty}_{tr}(\overline{SPT}_{\ge q})\subset \overline{SPT}_{\ge q}$.
\end{cor}
\begin{remark}\rm{}\rm
Note that out computation of the $l$-th cyclic power of $L_n$ in Section \ref{TMO} shows that $\overline{SPT}$ is not closed under the generalized symmetric products $\LL S^{\Phi}_{tr}$.
\end{remark}\rm{}
\begin{remark}\rm{}\rm
\llabel{rmref2}
We can make the computations done in the proof of Lemma  \ref{lcase} more precise as follows. First observe that Corollary \ref{homexact} implies that there are cofiber sequences in $H_{Nis,\af}(Cor(C,R))$ of the form
\begin{eq}
\llabel{keven}
\Sigma^1\LL S^l_{tr}(\Sigma^{2i+1}L^m)\sr \LL S_{tr}^l(\Sigma^{2i+2}L^m)\sr \Sigma^{2(i+1)l}L^{lm}
\end{eq}
and
\begin{eq}
\llabel{kodd}
\Sigma^1\LL S^l_{tr}(\Sigma^{2i}L^m)\sr \LL S_{tr}^l(\Sigma^{2i+1}L^m)\sr \Sigma^{2il+2}L^{lm}
\end{eq}
for all $k\ge 0$. Using the topological realization functor (see Section \ref{top.real} below) together with the fact that we know the topological  homology of $\tilde{S}^{l}(S^n)$ with coefficients in ${\bf F}_l$ (see \cite{Nakaoka}) one can show that these sequences split so that we have isomorphisms
$$\LL S_{tr}^l(\Sigma^{2i}L^m)=\Sigma^1\LL S^l_{tr}(\Sigma^{2i-1}L^m)\oplus  \Sigma^{2il}L^{lm}$$
and
$$\LL S_{tr}^l(\Sigma^{2i+1}L^m)=\Sigma^1\LL S^l_{tr}(\Sigma^{2i}L^m)\oplus \Sigma^{2il+2}L^{lm}$$
Using obvious induction on $i$  one can now get explicit formulas for the isomorphism classes of $\LL S_{tr}^l(\Sigma^iL^m)$.
\end{remark}\rm{}

\newpage \section{Eilenberg-MacLane spaces and their motives}\llabel{motEM}
Recall that $c(k)=1$ if $char(k)=0$ and $c(k)=char(k)$ if $char(k)>0$.  In this section we let $S$ denote the ring $\zz[1/c(k)]$. If $M$ is an abelian monoid we let $M^+$ denote its group completion and $M[1/d]$ the colimit of the sequence
$$M\stackrel{x\mapsto d\cdot x}{\longrightarrow}M\stackrel{x\mapsto d\cdot x}{\longrightarrow}\dots$$
Clearly, for any $M$ one has $(M[1/d])^+=M^+[1/d]$ and for an abelian group $A$ one has $A[1/c(k)]=A\oo S$. 

\subsection{Motivic Dold-Thom Theorem}
The goal of this section is to prove Theorems \ref{susvoe} which is the motivic analog of the topological Dold-Thom theorem and then to give some sufficient conditions for the equivalence between $S^{\infty}[1/c(k)]$ and its group completion $S^{\infty}[1/c(k)]^+$. The proof of Theorem \ref{susvoe} in the context of normal schemes goes back to \cite{SusVoe}.  

Let us recall the following definition (see \cite{Swan}, \cite{traverso}).
\begin{definition}
\llabel{seminormal}
A scheme $U$ is called semi-normal if it is reduced and any finite morphism  $U'\sr U$ such that $U'$ is reduced and for any field $K$ the map $U'(K)\sr U(K)$ is bijective is an isomorphism. 
\end{definition}
Let $SN/k$ be the category of semi-normal schemes over $k$. By Lemma \ref{semifad} it is $f$-admissible. 
\begin{lemma}
\llabel{univhomep}
Let $U$ be a semi-normal affine scheme of finite type over $k$ and $f:U'\sr U$ be a universal homeomorphism of finite type. Then there exists $n\ge 0$ such that ${\cal O}(U')^{c(k)^n}\subset {\cal O}(U)$. In particular, if $char(k)=0$ then $f$ is an isomorphism. 
\end{lemma}
\begin{proof}
Note first that a semi-normal scheme is necessarily reduced and therefore ${\cal O}(U)\sr {\cal O}(U')$ is a monomorphism. Since $f$ is a universal homeomorphism it is separated, universally closed and quasi-finite. Therefore, by Zariski theorem (see e.g. \cite[Th. 1.8]{Miln}), $f$ is a finite morphism. Then $U'$ is affine and ${\cal O}(U')$ is a finitely generated module over ${\cal O}(U)$. 

Let $K$ be a field and $x:Spec(K)\sr U$ a $K$-point of $U$. Since $f$ is a universal homeomorphism there exists a purely inseparable field extension $K\subset K'$ and a $K'$ point $x'$ of $U'$ lying over $x$. Moreover since ${\cal O}(U')$ is a finitely generated module over ${\cal O}(U)$ there exists $n$ such that for any $x$ as above one may choose $K'$ of degree dividing $c(k)^n$ over $K$. 

Let $R={\cal O}(U){\cal O}(U')^{c(k)^n}$. We claim that $R={\cal O}(U)$. Indeed the morphism $Spec(R)\sr U$ is finite and since for a purely inseparable $K\subset K'$ one has $(K')^{deg K'/K}\subset K$, our choice of $n$ implies that for any $K$ the map $Spec(R)(K)\sr U(K)$ is a bijection.  Since $U$ is assumed to be semi-normal we conclude that $U=Spec(R)$
\end{proof}
\begin{lemma}
\llabel{symprod}
Let $U$ be an affine scheme of finite type over $k$ and $d_n:{\cal O}(S^{c(k)^n}U)\sr {\cal O}(U)$ be the map defined by the diagonal. Then one has
$$Im(d_n)={\cal O}(U)^{c(k)^n}\cdot k.$$
\end{lemma}
\begin{proof}
The inclusion "$\supset$" is obvious. To prove the opposite inclusion consider ${\cal O}(U)$ as a vector space over $k$ and choose a basis $e_i$, $i\in I$ for it. The basis for ${\cal O}(U^N)={\cal O}(U)^{\oo N}$ is formed by the tensor products $e_{\uu{i}}=e_{i_1}\oo\dots\oo e_{i_N}$, $\uu{i}=(i_1,\dots,i_N)\in I^N$. The natural action of $S_N$ on ${\cal O}(U^N)$
is permutational relative  to this basis. Therefore, the subspace of invariants ${\cal O}(S^NU)$ coincide  with the subspace generated by expressions of the form $e_{A}=\sum_{\uu{i}\in A}e_{\uu{i}}$ where $A$ runs through the orbits of the action of $S_N$ on $I^N$. For convenience let us choose a linear ordering on $I$. Then we may describe these orbits by pairs of sequences $(i_1,\dots, i_m;j_1,\dots,j_m)$ where $i_1<\dots<i_m$ are in $I$, $j_1,\dots,j_m\in \zz_{>0}$ and $j_1+\dots+j_m=N$. Such a sequence defines the orbit $A(\uu{i},\uu{j})$ which contains the element $e_{i_1}^{\oo j_1}\oo\dots\oo e_{i_m}^{\oo j_m}$. 

The image of $e_{A(\uu{i};\uu{j})}$ under the diagonal map ${\cal O}(S^NU)\sr {\cal O}(U)$ is the element 
$$x_{(\uu{i};\uu{j})}={{N!}\over{j_1!\dots j_m!}}e_{i_1}^{j_1}\dots e_{i_m}^{j_m}$$
In our case $N=c(k)^n$. Since $(x_1+\dots+x_m)^{c(k)}=x_1^{c(k)}+\dots+x_m^{c(k)}$ in $k[x_1,\dots,x_m]$ we know that ${{N!}\over{j_1!\dots j_m!}}=0$ in $k$ unless $m=1$ when it equals $1$. Therefore for $N=c(k)^n$ the image of the diagonal map is generated by elements of the form $e_i^{c(k)^n}$. Lemma is proved. 
\end{proof}
\begin{proposition}
\llabel{isaqfh} 
For any $X$ the functor $U\mapsto Hom(U,\amalg_{n\ge 0}S^nX)[1/c(k)]$ is a sheaf in the qfh-topology on $SN/k$.
\end{proposition}
\begin{proof}
Our functor is a filtered colimit of representable functors and \cite[Proof of Th. 3.2.9]{H1} implies that the associated $qfh$-sheaf is of the form 
$$U\mapsto  colim_{U'\sr U} Hom(U',\amalg_{n\ge 0}S^nX)[1/c(k)]$$
over all universal homeomorphisms $U'\sr U$. To prove the proposition it remains to show that for a universal homeomorphism $U'\sr U$ such that both $U$ and $U'$ are semi-normal, the map 
$$Hom(U,\amalg_{n\ge 0}S^nX)[1/c(k)]\sr Hom(U',\amalg_{n\ge 0}S^nX)[1/c(k)]$$
is bijective. Since semi-normal schemes are reduced it is injective. It remains to show that for a map $U'\sr S^iX$ there exists $n$ and a map $U\sr S^{ic(k)^n}X$ such that the diagram
$$
\begin{CD}
U' @>>> S^iX\\
@VVV @VVV\\
U @>>> S^{ic(k)^n}X
\end{CD}
$$
commutes. It is sufficient to consider the case of affine $U,U'$ and $X$. Then the claim follows from Lemmas \ref{symprod} and \ref{univhomep}.
\end{proof}

\begin{proposition}
\llabel{susvoeplus}
Let $X$ be such that $S^nX$ exist and $Cor^{eff}(U,X)$ be the monoid of effective finite correspondences from $U$ to $X$. Then for any semi-normal $U$ one has:
$$Hom(U, \amalg_{n\ge 0} S^{n}(X))[1/c(k)]=Cor^{eff}(U,X)[1/c(k)].$$
\end{proposition}
\begin{proof}
By \cite[Proposition 4.2.7]{SusVoe2new} $Cor(-,X)^{eff}[1/c(k)]$ is a $qfh$-sheaf and by Proposition \ref{isaqfh} the same holds for $S^{\infty}[1/c(k)](X_+)$. On the other hand by \cite[Theorem 6.8]{SusVoe} we have 
$$Cor^{eff}(U,X)[1/c(k)]=Hom(U, S^{\infty}(X))[1/c(k)].$$
for any normal $U$.  Since any scheme has a qfh-covering by normal schemes we conclude that this equality holds for all semi-normal $U$. 
\end{proof}
\begin{remark}\rm\llabel{whysemi1}
It is not hard to see that there are natural maps  
$$Hom(U, \amalg_{n\ge 0} S^{n}(X))[1/c(k)]\sr Cor^{eff}(U,X)[1/c(k)]$$
for all $U$ and it might be the case that for $char(k)>0$ these maps are bijective for all $U$. For $char(k)=0$ i.e. $c(k)=1$ they are not necessarily surjective if $U$ is not semi-normal. For example let $U$ be the cuspidal cubic and $X=\af$. Then the graph of the normalization map $X\sr U$ is in $Cor^{eff}(U,X)$ but clearly not in the image of  $Hom(U, \amalg_{n\ge 0} S^{n}(X))$.
\end{remark}
\begin{theorem}
\llabel{susvoe}
Let $k$ be a perfect field and $C$ an $f$-admissible category which is contained in the category of semi-normal schemes. Consider $\Lambda^r_{S}\Lambda^l_{S}$ as a functor from $\Delta^{op}C_+^{\#}$ to $\Delta^{op}Rad(C_+)$. Then one has
\begin{eq}
\llabel{newplus}
\Lambda^r_S\Lambda^l_S=a_{Nis}(S^{\infty}[1/c(k)])^+
\end{eq}
\end{theorem}
\begin{proof}
Since both sides of (\ref{newplus}) commute with filtered colimits it is enough to very that for $X\in C$ one has:
$$\Lambda^r_S\Lambda^l_S(X_+)=a_{Nis}(S^{\infty}[1/c(k)](X_+))^+$$
By definition,  $\Lambda^r_S\Lambda^l_S$ takes $X_+$  to the functor 
\begin{eq}
\llabel{toinf1}
U\mapsto Cor(U,X)\oo S=Cor(U,X)[1/c(k)]
\end{eq}
where $Cor(U,X)$ is the group of finite correspondences from $U$ to $X$. By Proposition \ref{susvoeplus} we have 
$$S^{\infty}[1/c(k)](X_+)=Cor^{eff}(-,X)[1/c(k)]$$
On the other hand Lemma \ref{hens} below together with the fact that $Cor(-,X)$ is a sheaf in the Nisnevich topology, implies that $Cor(-,X)=a_{Nis}(Cor^{eff}(-,X)^+)$. Theorem is proved.
\end{proof}
\begin{lemma}
\llabel{hens}
Let $U$ be a henselian local scheme. Then
$$Cor(U,X)=Cor^{eff}(U,X)^+$$
\end{lemma}
\begin{proof}
Since $Cor(U,X)=Cor(U_{red},X)$ and the same holds for $Cor^{eff}$ we may assume that $U$ is reduced. Let ${\cal Z}=\sum n_i z_i$ be a finite correspondence from $U$ to $X$ i.e. a relative finite cycle on $X_U=X\times U$ over $U$.  The individual points $z_i$ need not be relative cycles over $U$ but by the definition of a relative finite cycle  (cf. \cite{SusVoe2new}) we know that the points $z_i$ lie over the generic points of $U$ and that the closure $[z_i]$ of each $z_i$ is finite over  $U$. Therefore, in order to show that  $\cal Z$ is a difference of two effective relative cycles it is sufficient to show that for any $z$ which lies over a generic point of $U$ and such that its closure $[z]$ is finite over $U$ there exists an effective relative finite cycle $\cal Z$ such that $Supp({\cal Z})$ contains $z$. Since $U$ is assumed to be local henselian any effective relative cycle of relative dimension zero over $U$  is a sum of an effective relative finite cycle and a cycle whose support lies over the complement to the closed point  of $U$. Therefore it is sufficient to find an effective relative cycle $\cal Z$ of relative dimension  $0$ such that $Supp({\cal Z})$ contains $z$. 

The closure $[z]$ is a henselian local scheme. The image of the closed point of $[z]$ in $X$ lies on one or more of the irreducible components $X$. Replacing $X$ by one of such components we may assume that $X$ is irreducible and in particular equidimensional. We may further replace $X_U$ by an affine open neighborhood of the closed point of $[z]$.  This reduces the problem to the situation when we have a flat equidimensional morphism $p:X_U\sr U$ of some dimension $d$ of affine schemes such that $U$ is reduced and a point $z$ over a generic point of $U$ such that $[z]\sr U$ is finite and we need to find a effective relative cycle ${\cal Z}$ of relative dimension $0$ such that $z\in Supp({\cal Z})$. 

Since $X_U$ is flat and equidimensional the fundamental cycle of $X_U$ is a relative cycle of dimension $d$ over $U$ and clearly $z$ belongs to its support. Assume by induction that there exists an effective equidimensional cycle ${\cal Z}=\sum_i n_i z_i$ of relative dimension $r$ over $U$ such that $z\in Supp({\cal Z})$. If $r=0$ then we are done. Suppose that $r>0$. Since $X_U$ is affine there is a regular function $f$ on $X$ which is zero in the closed point of $[z]$ and non zero in the generic points of the closed fibers of $[z_i]\sr U$ for all $i$. By \cite{cancellationsub} the cycle $({\cal Z},D(f))=\sum_i n_i D(f_i)$ where $f_i$ is the restriction of $f$ to $[z_i]$ is a relative equidimensional cycle of relative dimension $r-1$ over $U$ and one obviously has $z\in Supp(({\cal Z},D(f)))$.
\end{proof}
\begin{remark}\rm
The previous remark (\ref{whysemi1}) shows that if we consider all schemes instead of semi-normal ones then (\ref{newplus}) stops being an isomorphism at least in characteristic zero. It is possible that there is still an $(Nis,\af)$-equivalence of the form (\ref{whysemi1}) for all schemes. It is also clear that there is an equivalence of the same form with respect to a modification of the Nisnevich topology which allows for semi-normalizations as coverings.
\end{remark}

In may cases, including the main case of Moore spaces considered in the next section it is important to know when one can further simplify the description of $\Lambda^r_S\Lambda^l_S(X)$ given in Theorem \ref{susvoe} by replacing $S^{\infty}[1/c(k)]^+$ with $S^{\infty}[1/c(k)]$. The goal of the rest of this section is to prove some partial results in this direction.
\begin{definition}
\llabel{2009.09.07.def1}
Let $C$ be an $f$-admissible subcategory. An object $X$ of $\Delta^{op}Rad(C_+)$ is said to satisfy condition (D1) if for any henselian local scheme $U$ of a pointed scheme from $C$, the abelian monoid $\pi_0(\Delta C_*(S^{\infty}[1/c(k)](X))(U))$ is a group.
\end{definition}
\begin{proposition}
\llabel{2009.09.07.pr1}
If $X$ satisfies condition (D1) then the natural morphism 
$$S^{\infty}[1/c(k)](L_*(X))\sr S^{\infty}[1/c(k)]^+(L_*(X))$$
is a $(Nis,\af)$-equivalence.
\end{proposition}
\begin{proof}
Consider the diagram of morphisms of monoids in $\Delta^{op}Rad(C_+)$ of the form
$$
\begin{CD}
S^{\infty}[1/c(k)](L_*(X))@>>> \Delta C_*(S^{\infty}[1/c(k)](L_*(X)))\\
@VVV @VVV\\
S^{\infty}[1/c(k)](L_*(X))^+ @>>> (\Delta C_*(S^{\infty}[1/c(k)](L_*(X))))^+\cong \Delta C_*(S^{\infty}[1/c(k)](L_*(X))^+)
\end{CD}
$$
By Proposition \ref{2009.07.24.1} the horizontal morphisms are $\af$-equivalences. We have $\pi_0(X)=\pi_0(L_*(X))$ since $L_*(X)\sr X$ is a projective equivalence. We also have 
$$\pi_0(\Delta C_*(Y))=\pi_0(\Delta C_*(\pi_0(Y))).$$
Indeed, for a bisimplicial set $B$ one has $\pi_0(\Delta B)=\pi_0(\Delta \pi_0^{vert}(B))$ where $\pi_0^{vert}(B)$ is obtained by replacing each column of $B$ by its $\pi_0$. Therefore, $\pi_0(\Delta C_*(Y))=\pi_0(\Delta \pi_0^{vert}(C_*(Y)))$. Consider $Y$ as going in the vertical direction. Then $\pi_0^{vert}(C_*(Y))=C_*(\pi_0(Y))$ since $C_*$ commutes with reflexive coequalizers (actually with all colimits).  

Finally,
$$\pi_0(S^{\infty}[1/c(k)](Y))=S^{\infty}[1/c(k)](\pi_0(Y))$$
since $S^{\infty}[1/c(k)]$ being a radditive extension commutes with reflexive coequalizers. We conclude that for any $Y$, we have 
$$\pi_0(\Delta C_*(S^{\infty}[1/c(k)](L_*(X))))=\pi_0(\Delta C_*(\pi_0(S^{\infty}[1/c(k)](L_*(X)))))=$$
$$=\pi_0(\Delta C_*(S^{\infty}[1/c(k)](\pi_0(L_*(X)))))=\pi_0(\Delta C_*(S^{\infty}[1/c(k)](\pi_0(X))))=$$
$$=\pi_0(\Delta C_*(\pi_0(S^{\infty}[1/c(k)](X))))=\pi_0(\Delta C_*(S^{\infty}[1/c(k)](X)))$$
Therefore, the right hand side vertical morphism  is a $Nis$-equivalence by our assumption and Lemma \ref{ge64} below.  We conclude that the left hand side vertical morphism is a $(Nis,\af)$-equivalence. 
\end{proof}
\begin{lemma}
\llabel{ge64}
Let $X$ be a  simplicial abelian monoid such that $\pi_0(X)$ is a group. Then the natural map $X\sr X^+$ is a weak equivalence.
\end{lemma}
\begin{proof}
See \cite[p. 381]{puppemonoid}.
\end{proof}
\begin{lemma}
\llabel{2009.09.07.l2}
For any $f$-admissible $C$ and any $X\in \Delta^{op}Rad(C_+)$ the object $S^1_s\wedge X$ satisfies condition (D1). 
\end{lemma}
\begin{proof}
One can easily see that  $\pi_0(\Delta C_*(S^{\infty}[1/c(k)](S^1_s\wedge X)))=0$.
\end{proof}
\begin{lemma}
\llabel{2009.09.07.l3}
For any $f$-admissible $C$ one has:
\begin{enumerate}
\item if $X, Y$ satisfy condition (D1) then $X\vee Y$ satisfies condition (D1),
\item if $f:X\sr Y$ is a morphism in $\Delta^{op}Rad(C_+)$ such that $X$ and $Y$ satisfy condition (D1) then $cone(f)$ satisfies condition (D1),
\end{enumerate}
\end{lemma}
\begin{proof}
The first assertion follows from the formula
$$\pi_0(\Delta C_*(S^{\infty}[1/c(k)](X\vee Y)))=\pi_0(\Delta C_*(S^{\infty}[1/c(k)](X)\times S^{\infty}[1/c(k)](X)))=$$
$$=\pi_0(\Delta C_*(S^{\infty}[1/c(k)](X)))\times \pi_0(\Delta C_*(S^{\infty}[1/c(k)](Y)))$$
where the first equality holds by Lemma \ref{wasnoth} and the second since by $C_*$ and $\pi_0$ commute with products. 

To prove the second assertion we may assume that $f$ is a term-wise coprojection. Then $cone(f)$ is projectively equivalent to the coequalizer of the reflexive pair $X\vee Y\dsr Y$ in which one map equals $f$ on $X$ and $Id_Y$ on $Y$ and the second one maps $X$ to the distinguished point and again equals $Id_Y$ on $Y$. Applying the functor  $ \pi_0(\Delta C_*(S^{\infty}[1/c(k)](-)))$ and using the fact that it commutes with reflexive coequalizers and takes $\vee$ to the direct product we conclude that $\pi_0(\Delta C_*(S^{\infty}[1/c(k)](cone(f))))$ is a reflexive coequalizer of a diagram of groups and therefore a group. 
\end{proof}
Recall that we let $S^1_t$ denote the pointed scheme $(\af-\{0\};1)$ which we identify with the radditive functor which it represents on $C_+$. 
\begin{proposition}
\llabel{sqgr}
For any $f$-admissible $C$ which is contained in $SN/k$ and $X\in \Delta^{op}Rad(C_+)$ the objects $S^1_t\wedge X$ satisfies condition (D1). 
\end{proposition}
\begin{proof}
We may clearly assume that $X\in \Delta^{op}C_+^{\#}$. In view of Proposition \ref{susvoeplus}, for any $U\in C$ and $Y_+\in C_+$ we have
$$Hom(U,S^{\infty}[1/c(k)](S^1_t\wedge Y_+))=Cor^{eff}[1/c(k)](U,(\af-\{0\})\times Y)/Cor^{eff}[1/c(k)](U,\{1\}\times Y)$$
Therefore, the $\pi_0$ set which we consider is the set of $\af$-homotopy classes of maps from $[U]$ to $[\af-\{0\}]\oo [Y]$, in the subcategory of $Cor(C,S)$ where morphisms are effective finite correspondences, modulo those classes which contain correspondences landing in $[\{1\}]\oo [Y]$. Let $\phi:\af-\{0\}\sr \af-\{0\}$ be the morphism $z\mapsto z^{-1}$. For a morphism $f:[U]\sr [\af-\{0\}]\oo [Y]$ in our category set $f^{-}=(\phi\oo Id_{[Y]})\circ f$. Let us show that $f+f^{-1}$ is $\af$-homotopic to a correspondence which goes to zero. Note first that we have
$$f+(\phi\oo Id_{[Y]})\circ f=(Id\oo Id_{[Y]})\circ f + (\phi\oo Id_{[Y]})\circ f=((Id+\phi)\oo Id_{[Y]})\circ f$$
Therefore it is enough to show that the finite correspondence $\phi+Id$ from $\af-\{0\}$ to itself is $\af$-homotopic  in the category of effective finite correspondences to a finite correspondences which lands in $1$. The set of finite correspondences of degree $d$ from any smooth $X$ to $\af-0$ is in a natural bijection with $Hom(X,S^d(\af-\{0\}))$. The standard theory of symmetric polynomials  shows that $S^2(\af-\{0\})=\af\times (\af-\{0\})$ where the projection to $\af$ corresponds to the polynomial $-X-Y$ and the projection to $\af-\{0\}$ to the polynomial $XY$. The correspondence $\phi+Id$ is represented under this bijection by a morphism $\af-\{0\}\sr  \af\times (\af-\{0\})$ whose image lies in $\af\times \{1\}$ and which is, therefore, $\af$-homotopic to the morphism which sends $\af-\{0\}$ to $\{1\}+\{1\}=(-2,1)$. 
\end{proof}
\begin{remark}
\llabel{strange}\rm
Let $c(k)=1$. In view of Proposition \ref{susvoeplus} the effective analog $\Lambda^r_{eff}$ of $\Lambda^r$ maps $(Cor^{eff})^{\#}$ to $C_+^{\#}$ and we have a pair of adjoint functors $\Lambda^l_{eff}$, $\Lambda^r_{eff}$ between these two categories which is analogous to the pair $\Lambda^l$, $\Lambda^r$. 

Consider the composition $\Lambda^l_{eff}\Lambda^r_{eff}:(Cor^{eff})^{\#}\sr (Cor^{eff})^{\#}$. By Proposition \ref{susvoeplus} this functor coincides {\em on objects} with the functors $S^{\infty}_{tr,eff}$ which is the obvious effective variant of the functor $S^{\infty}_{tr}$. Note however that these two functors {\em do not agree on morphisms}. For example, if $d\,Id:[X]\sr [X]$ is the $d$-th multiple of the identity map on $[X]$ then $\Lambda^l_{eff}\Lambda^r_{eff}(d\, Id)=\Lambda^l_{eff}(\times d)$ where $\times d$ is the multiplication by $d$ map $S^{\infty}(X)\sr S^{\infty}(X)$ while $S^{n}_{tr,eff}(d\, Id)=d^n\, Id$ and 
$$S^{\infty}_{tr,eff}(d\, Id)=\oplus_{n\ge 0} S^{n}_{tr,eff}(d\, Id)$$
Since $\Lambda^l_{eff}\Lambda^r_{eff}$ and $S^{\infty}_{tr,eff}$ do not agree on morphisms their extensions to the corresponding categories of simplicial objects do not agree on objects. This effect remains after we pass to group completions. In particular, the isomorphisms 
$$\Lambda^l\Lambda^r([X])\cong S^{\infty}_{tr}([X])^+$$
do not extend to isomorphisms between the functors $\Lambda^l\Lambda^r$ and $S^{\infty}_{tr}(-)^+$. 
\end{remark}

\subsection{Motivic Moore pairs and the main structure theorem}
\llabel{moore.section}
Let us recall the definition of motivic cohomology outlined in the introduction. Let  $C$ be an admissible category. Consider the motivic spheres $S^1_t=(\af-\{0\},1)$ and $S^q_t=(S^1_t)^{\wedge q}$ as radditive functors on $C_+$.  Set
$$l_{q,R}=\LL\Lambda^l_{R}(S^q_t)=\Lambda^l_{R}(L_*(S^q_t)).$$
We will write $l_q$ for $l_{q,\zz}$. The (reduced) {\em unstable} motivic cohomology of $X\in \Delta^{op}C_+^{\#}$ with coefficients in an abelian group $A$ is defined by the formula
\begin{eq}
\llabel{descr0}
\tilde{H}^{p,q}_{un}(X,A)=\left\{\begin{array}{lll}
Hom_{H_{Nis,\af}(Cor(C,\zz))}(\Lambda^lX,\,\,\Sigma^{p-q}(A\oo_{\LL}l_q))&\rm for&p\ge q\\\\
Hom_{H_{Nis,\af}(Cor(C,\zz))}(\Sigma^{p-q}\Lambda^lX,\,\,A\oo_{\LL}l_q)&\rm for&p\le q
\end{array}
\right.
\end{eq}
where $\oo_{\LL}$ is the derived tensor product given by
$$A\oo_{\LL}l_q=L_*(A)\oo_{\zz} l_q.$$
\begin{definition}
\llabel{2009.09.06.def1}
Let $A$ be an abelian group, $p,q\in\zz$ and $C$ be an admissible subcategory of $Sch/k$. A motivic Eilenberg-MacLane pair for $(A,p,q)$ on $C$ is a pair $(X,\iota)$ where $X\in \Delta^{op}Rad(C_+)$ and $\iota\in \tilde{H}^{p,q}_{un}(X,A)$, which represents $\tilde{H}^{p,q}_{un}(-,A)$ on $H_{Nis,\af}(C_+)$.
\end{definition}
By the adjunction between $\Lambda^r$ and $\LL\Lambda^l$ the functors $\tilde{H}^{p,q}_{un}(-,A)$ can be represented on $H_{Nis,\af}(C_+)$ by the spaces
$\Lambda^r\Sigma^{p-q}(A\oo_{\LL}l_q)$ for $p\ge q$ and 
$\Lambda^r\Omega^{p-q}_{Nis,\af}(A\oo_{\LL}l_q)$ for $p<q$. In particular, motivic Eilenberg-MacLane pairs exist for all $A$, $p$ and $q$. By the uniqueness of representing objects the isomorphism class of a $(X,\iota)$ in $H_{Nis,\af}(C_+)$ is well defined up to a canonical isomorphism and we will denote it by $(K(A,p,q)_C,\iota_{p,q})$. 
\begin{lemma}
\llabel{2009.09.17.l2}
For any $k$, any $A$, any $p,q\in\zz$ and any inclusion $i:C\sr D$ of admissible subcategories there is a canonical isomorphism
$$K(A,p,q)_C=i_{rad,+}(K(A,p,q)_D)$$
\end{lemma}
\begin{proof}
It follows immediately from the definition of $K(A,p,q)$ as a representing object and the adjunction between $i_{rad,+}$ and $\LL i^{rad}_+$.
\end{proof}

In what follows we will only consider the case of $p\ge q$. The case of $p<q$ is much more complicated. 

\begin{definition}
\llabel{2009.09.06.def2}
Let $A$ be an abelian group, $p\ge q\ge 0$ and $C$ be an admissible subcategory of $Sch/k$. A motivic Moore pair for $(A,p,q)$ on $C$ is a pair $(X,\phi)$ where $X\in \Delta^{op}Rad(C_+)$ and $\phi$ is an isomorphism $\LL\Lambda^l_{\zz}(X)\cong \Sigma^{p-q}(A\oo_{\LL} l_q)$ in $H_{Nis,\af}(Cor(C,\zz))$.
\end{definition}
\begin{proposition}
\llabel{2009.09.11.pr1}
If $(X,\phi)$ is a Moore pair for $(A,p,q)$ on $C$ then there is a natural class $\iota$ in $\tilde{H}^{p,q}_{un}(\Lambda^r_S\LL\Lambda^l_S(X), A\oo_{\zz} S)$ such that $(\Lambda^r_S\LL\Lambda^l_S(X),\iota)$ is an Eilenberg-MacLane pair for $(A\oo_{\zz} S,p,q)$ on $C$. 
\end{proposition}
\begin{proof}
It follows immediately from the standard adjunctions that for any $R$ and $A$ the image of the space $\Lambda^r_R\LL\Lambda^l_R(M(A,p,q))$ in $H_{Nis,\af}(C_+)$ represents the functor 
$$\tilde{H}^{p,q}_{un}(-,A\oo_{\LL,\zz} R)=Hom_{H_{Nis,\af}(Cor(C,\zz))}(-,\Sigma^{p-q}((A\oo_{\LL,\zz} R)\oo l_q)).$$
The claim of the proposition follows from the fact that $A\oo_{\LL,\zz} S=A\oo_{\zz} S$.
\end{proof}

\begin{proposition}
\llabel{2009.09.06.th1}
For any $f$-admissible $C$ which is contained in $SN/k$, any $p\ge q\ge 0$ such that $p>0$ and any finitely generated abelian group $A$ there exists a motivic Moore pair $(X(A,p,q),\phi)$ for $(A,p,q)$ on $C$ such that
\begin{enumerate}
\item $X(A,p,q)$ satisfies condition (D1),
\item for any $d>0$ there exists a morphism $m_d:X(A,p,q)\sr X(A,p,q)$ such that  
\begin{enumerate}
\item $\LL S^{\infty}(m_d)=\times d$ in $H_{Nis,\af}(C_+)$,
\item $\LL \Lambda_{\zz}^l(m_d)=d\cdot Id$.
\end{enumerate}
\end{enumerate}
\end{proposition}
\begin{proof}
It follows from Lemma \ref{2009.09.07.l3}(1) and Lemma \ref{wasnoth} that if $X(A,p,q)$ and $X(B,p,q)$ satisfy the conditions of the proposition for $(A,p,q)$ and $(B,p,q)$ respectively then $X(A,p,q)\vee X(B,p,q)$ satisfies its conditions for $(A\oplus B,p,q)$. Since any finitely generated abelian group is a direct sum of finite number of cyclic groups it means that it is sufficient to construct $X(\zz/n,p,q)$ for $n>0$ and $X(\zz,p,q)$. We consider the following two cases.
{\em $(q=0, p>0)$} Let $S^1_{s,gr}$ be a model of $S^1_s$ which is an abelian group e.g. $S^1_{s,gr}=K(\zz,1)$. Let $\cdot n:S^1_{s,gr}\sr S^1_{s,gr}$ be the multiplication by $n$ map with respect to the abelian group structure. Set
$$X(\zz,p,0)=S^1_{s,gr}\wedge S^{p-1}_s$$
$$X(\zz/n,p,0)=cone(\cdot n)\wedge S^{p-1}_s$$
These spaces satisfy the first condition of the proposition by Lemmas \ref{2009.09.07.l2} and \ref{2009.09.07.l3}(2). The second condition is easily seen to hold relative to the maps $m_d$ which are defined by the map $\cdot d$ on $S^1_{s,gr}$. 

{\em $(q\ge 1, p\ge q)$} Let $(-)^{n}:S^1_t\sr S^1_t$ be the map $z\mapsto z^n$. Set
$$X(\zz,p,q)=S^1_{t}\wedge S^{q-1}_t\wedge S^{p-q}_s$$
$$X(\zz/n,p,q)=cone((-)^n)\wedge  S^{q-1}_t\wedge S^{p-q}_s$$
These spaces satisfy the first condition of the proposition by Lemmas \ref{sqgr} and \ref{2009.09.07.l3}(2). Let us shows that they satisfy the second condition relative to the maps $m_d$ defined by the maps $(-)^d$ on $S^1_t$. The second half of the second condition is easy. To prove the first part observe that by  Lemma \ref{muls} below it is sufficient to verify that  the morphism $(-)^d$ defines on $S^{\infty}(L_*(S^1_t))$ and $S^{\infty}(L_*(cone((-)^n)))$ the morphisms which coincide with $\times d$ in $H_{Nis,\af}(C_+)$.

By Proposition \ref{2009.08.05.pr1} the map
$$S^{\infty}(L_*(S^1_t))\sr S^{\infty}(S^1_t)$$
is a projective equivalence. On the other hand the standard theory of symmetric polynomials implies that there is a homomorphism of monoids 
\begin{eq}
\llabel{2009.09.11.eq1}
S^{\infty}(S^1_t)\sr S^1_t
\end{eq}
which is an $\af$-equivalence. This proves the case of $S^{\infty}(L_*(S^1_t))$. Using (\ref{2009.09.11.eq1}) it is not hard to show further that there is an $\af$-equivalence of monoids
$$S^{\infty}(L_*(cone((-)^n)))\sr K(S^1_t\stackrel{(-)^n}{\sr} S^1_t)$$
where $K$ is the standard functor from complexes of abelian groups to simplicial objects, which proves the case of $S^{\infty}(L_*(cone((-)^n)))$.
\end{proof}
\begin{remark}\rm
Note that if $char(k)=0$ then conditions (1) and (2b) of Proposition \ref{2009.09.06.th1} imply condition  (2b). For $char(k)>0$ this is not necessarily the case because we do not know how to reconstruct  $S^{\infty}$ from $\Lambda_{\zz}^l$ even for objects satisfying condition (D1).
\end{remark}

\begin{lemma}
\llabel{muls}
Let $X$ be an object of $\Delta^{op}C_+^{\#}$ and $m_d:X\sr X$ a morphism such that $S^{\infty}(m_d)=\times d$ in $H_{Nis,\af}(C_+)$. Then for any $Y\in \Delta^{op}C_+^{\#}$ one has
$$S^{\infty}(m_d\wedge Id_Y)=\times d$$
in $H_{Nis,\af}(C_+)$. 
\end{lemma}
\begin{proof}
Let $f_d$ be the composition $X\sr S^{\infty}(X)\stackrel{\times d}{\sr} S^{\infty}(X)$ and $g_d$ be the composition of $X\stackrel{m_d}{\sr}X\sr S^{\infty}(X)$ or equivalently the composition of the natural map $X\sr S^{\infty}(X)$ with $S^{\infty}(m_d)$. Let
$$h:S^{\infty}(X)\wedge Y\sr S^{\infty}(X\wedge Y)$$
be the map which takes $(x_1+\dots+x_n;y)$ to $(x_1;y)+\dots + (x_n;y)$. Then there are commutative diagrams (in $\Delta^{op}C_+^{\#}$)
$$
\begin{CD}
S^{\infty}(X\wedge Y) @>\times d>> S^{\infty}(X\wedge Y) \\
@VS^{\infty}(f_d\wedge Id_Y)VV @AAA\\
S^{\infty}(S^{\infty}(X)\wedge Y) @>S^{\infty}(h)>> S^{\infty}(S^{\infty}(X\wedge Y))
\end{CD}
$$
and 
$$
\begin{CD}
S^{\infty}(X\wedge Y) @>S^{\infty}(m_d\wedge Id_Y)>> S^{\infty}(X\wedge Y) \\
@VS^{\infty}(g_d\wedge Id_Y)VV @AAA\\
S^{\infty}(S^{\infty}(X)\wedge Y) @>S^{\infty}(h)>> S^{\infty}(S^{\infty}(X\wedge Y))
\end{CD}
$$
where the right hand side arrows are from the standard triple structure on $S^{\infty}$. Since $f_d=g_d$ in $H_{Nis,\af}(C_+)$ we conclude that $\times d=S^{\infty}(m_d\wedge Id_Y)$ in this category as well.
\end{proof}

For any $X$ the isomorphism of Theorem \ref{susvoe} defines a morphism
$$dt_X:S^{\infty}(X)\sr \Lambda^r_S\Lambda^l_S(X)$$
\begin{proposition}
\llabel{2009.09.11.pr2}
Let $(X(A,p,q),\phi)$ be a Moore pair satisfying the conditions of Proposition \ref{2009.09.06.th1}. Then for any commutative algebra $R$ over $S=\zz[1/c(k)]$ the morphism 
$$\LL\Lambda_R(dt_{L_*(X(A,p,q))}):\LL\Lambda^l_R(\LL S^{\infty}(X(A,p,q)))\sr \LL\Lambda^l_R(\Lambda^r_S\LL\Lambda^l_S(X(A,p,q)))$$
is an isomorphism in $H_{Nis,\af}(Cor(C,R))$.
\end{proposition}
\begin{proof}
By Theorem \ref{susvoe} we have
$$S^{\infty}[1/c(k)]^+(L_*(X(A,p,q)))=\Lambda^r_S\LL\Lambda^l_S(X(A,p,q))$$
since $X(A,p,q)$ satisfies condition (D1) so does $L_*(X(A,p,q))$ and we may replace $S^{\infty}[1/c(k)]^+$ by $S^{\infty}[1/c(k)]$. On the other hand 
$$\LL\Lambda_R^l(\times c(k))=\LL\Lambda^l_R(S^{\infty}(m_{c(k)}))=S^{\infty}_{tr}(\LL\Lambda^l_R(m_{c(k)}))=S^{\infty}_{tr}(\Lambda^l_R(m_{c(k)}))=S^{\infty}_{tr}(c(k)\cdot Id)$$
is an isomorphism since $c(k)$ is invertible in $R$ (note that condition (2b) of Proposition \ref{2009.09.06.th1} implies a similar condition for $\Lambda_R^l$). Therefore the morphism
$$\LL\Lambda_R^l(\LL S^{\infty}(X(A,p,q)))\sr \LL\lambda_R^l(\LL S^{\infty}[1/c(k)](X(A,p,q)))$$
is an isomorphism in $H_{Nis,\af}(Cor(C,R))$. 
\end{proof}
Let us introduce the following notation:
$$M(A,p,q;R)_C=\LL\Lambda_R^l(K(A,p,q)_C)$$
The following result is the main theorem of this section. 
\begin{theorem}
\llabel{2009.09.12.th1}
Let $k$ be a perfect field, $C$ an $f$-admissible category which is contained in the category of semi-normal schemes, $A$ a finitely generated module over $S=\zz[1/c(k)]$ and $p,q$ two integers such that $p\ge q\ge 0$ and $p>0$. Then for any $S$-algebra $R$ there is an isomorphism in $H_{Nis,\af}(Cor(C,R))$ of the form
\begin{eq}
\llabel{2009.09.12.eq1}
M(A,p,q;R)_C\cong \oplus_{n\ge 0}S^{n}_{tr}(\Sigma^{p-q}((A\oo_{\LL,S} R)\oo l_{q,R}))
\end{eq}
such that for any $d>0$ the map $M(A,p,q;R)_C\sr M(A,p,q;R)_C$ defined by multiplication by $d$ in $A$ is of the form $\oplus_{n\ge 0} d^n\cdot Id$. A choice of such an isomorphism is determined  by a choice of a finitely generated abelian group $A'$, an isomorphism $A=A'\oo S$ and a motivic Moore pair $(X(A',p,q),\phi)$ satisfying the conditions of Proposition \ref{2009.09.06.th1}. 
\end{theorem}
\begin{proof}
Any finitely generated $S$-module $A$ is of the form $A=A'\oo S$ for a finite generated abelian group $A'$.  Let us choose a Moore pair $X(A',p,q)$ which satisfies the conditions of Proposition \ref{2009.09.06.th1}. By Proposition \ref{2009.09.11.pr1} and Proposition \ref{2009.09.11.pr2} we have an isomorphism
$$\LL\Lambda^l_R(K(A,p,q))\cong \LL \Lambda^l_R(\LL S^{\infty}(X(A',p,q)))=\LL S^{\infty}_{tr}(\LL \Lambda^l_R(X(A',p,q)))$$
By an obvious "universal coefficients formula" the isomorphism
$$\LL\Lambda^l_{\zz}(X(A',p,q))=\Sigma^{p-q}(A'\oo_{\LL,\zz} l_q)$$
defines an isomorphism
$$\LL \Lambda^l_R(X(A',p,q))=\Sigma^{p-q}((A'\oo_{\LL,\zz} R)\oo_{\LL,R} l_{q,R})=\Sigma^{p-q}((A\oo_{\LL,S} R)\oo_{\LL,R} l_{q,R})$$
and therefore an isomorphism of the form (\ref{2009.09.12.eq1}).

The map $\LL\Lambda^l_R(K(A,p,q))\sr \LL\Lambda^l_R(K(A,p,q))$ defined by multiplication by $d$ in $A$ corresponds under the isomorphism $\LL\Lambda^l_R(K(A,p,q))\cong \LL \Lambda^l_R(\LL S^{\infty}(X(A',p,q)))$ to the map $\LL\Lambda^l_R(\times d)$ and the properties of the map $m_d$ show that 
$$\LL\Lambda^l_R(\times d)=\LL\Lambda^l_R(\LL S^{\infty}(m_d))=\LL S^{\infty}_{tr}(\LL \Lambda^l_R(m_d))=\LL S^{\infty}_{tr}(d\cdot Id)=\oplus_{n\ge 0} d^n\cdot Id.$$
\end{proof}
\begin{cor}
\llabel{fcase}
Under the assumptions of the theorem assume in addition that $R=F$ is a field. Then a choice of $(X(A,p,q),\phi)$ defines isomorphisms 
$$M(A,p,q;F)_C=(M_0\oo M_1)\oplus M_0\oplus M_1$$
where
$$M_0=\LL S^{\infty}_{tr}(\Sigma^{p-q}((A\oo_{\zz} F)\oo_F l_{q,F}))$$
and
$$M_1=\LL S^{\infty}_{tr}(\Sigma^{p-q+1}(Tor_1^{\zz}(A,F)\oo_F l_{q,F})).$$
\end{cor}
\begin{proof}
It follows from the general case by Corollary \ref{directsuminfty} since we have a canonical decomposition 
$$A\oo_{\LL,\zz} F=(A\oo_{\zz} F)\oplus \Sigma^1Tor_1^{\zz}(A,F).$$ 
\end{proof}
\begin{cor}
\llabel{flz}
Under the assumptions of the theorem let  $R=F$ be a field of characteristic $l>0$ such that $l\ne char(k)$. Then a choice of $(X(\zz/l,p,q),\phi)$ defines isomorphisms:
$$M(\zz/l,p,q;F)_C=(M_0\oo M_1)\oplus M_0\oplus M_1$$
where $M_0=M(\zz,p,q;F)_C$ and $M_1=M(\zz,p+1,q;F)_C$.
\end{cor}
\begin{cor}
\llabel{2009.09.15.cor1}
Under the assumptions of the theorem assume in addition that $R=F$ is a field and $p\ge 2q$. Then one has
$$M(A,p,q;F)_C\in \overline{SPT}_{\ge q}$$
\end{cor}
\begin{proof}
The case when $p=q=0$ is obvious. For $p>0$ the statement  follows easily from Corollary \ref{fcase} and Theorem \ref{main2}.
\end{proof}
\begin{remark}
\llabel{2009.09.15.rem1}\rm
It is not clear whether or not there exist motivic Moore pairs $(X(A,p,q),\phi)$ satisfying the conditions of Proposition \ref{2009.09.06.th1} which define isomorphisms (\ref{2009.09.12.eq1}) which are different from the ones which are defined by the Moore pairs which are constructed in the proof of this proposition. One can look at this problem from the following angle. The composition $F=\LL\Lambda^l_{\zz}\Lambda^r_{\zz}$ is a cocomplete-triple on $H_{Nis,\af}(Cor(C,\zz))$. A choice of a motivic Moore pair for $(A,p,q)$ defines a structure of an $F$-coalgebra on $\Sigma^{p-q}(A\oo_{\LL,\zz} l_q)$ and a choice of such a structure defines an isomorphism of the form (\ref{2009.09.12.eq1}). Therefore in order to obtain an "exotic" isomorphism of this form we need a Moore space $X(A,p,q)$ which defines an exotic $F$-coalgebra structure on $\Sigma^{p-q}(A\oo_{\LL,\zz} l_q)$. On the other hand, one can easily see that such a structure is determined by the action of cohomological operations on the motivic cohomology of $X(A,p,q)$. In the topological context these observations show that there are no exotic coalgebra structures since there can be no non-trivial actions of cohomological operations on the cohomology of Moore spaces. In the motivic context they show that in order to construct an "exotic" Moore space for $\zz$ we need to find a pointed space $X$ whose motive is a Tate motive $\zz(p)[q]$ but the action of the motivic cohomological operations on $H^{*,*}(X,\zz)=H^{*-p,*-q}(Spec(k),\zz)$ is different from their action on the motivic cohomology of the point. Since we know almost nothing about unstable motivic cohomology operations of bi-degree $(i,j)$ where $2i\ge j\ge i$ the question of exotic Moore spaces remains open. 
\end{remark}
\begin{remark}\rm
\llabel{2009.09.15.rem2}\rm
The construction of the proof of Proposition \ref{2009.09.06.th1} can be easily adjusted to provide 
motivic Moore pairs for $(A,p,q)$ where $A$ is any $S$-module and $p>q$. However, it is not clear how to construct motivic Moore pairs for $(A,p,p)$ when $p>0$ and $A$ is an indecomposable (infinitely generated) $S$-module of rank greater than $1$. Multiple examples of such modules can be found in \cite{Fuchs}.
\end{remark}
\begin{remark}\rm
Theorem \ref{2009.09.12.th1} together with Example \ref{2009.09.15.ex1} shows that for $2q>p\ge q$ the image of $M(A,p,q;R)$ in $DM^{eff}_{-}$ is a "mixed Tate object" but not necessarily a pure Tate object. For $p<q$ this image may not even be a mixed Tate object as can be seen on the example of $M(\zz/l,0,1;\qq)$. The Eilenberg-MacLane object $K(\zz/l,0,1)$ is the scheme of $l$-roots of unity and its motive is not a Tate motive unless the $l$-root of unity is in $k$. 
\end{remark}\rm{}
\begin{theorem}
\llabel{mth1}
Under the assumptions of Theorem \ref{2009.09.12.th1} assume in addition that $R=F$ is a field, $p\ge 2q$ and that $k$ admits resolution of singularities in the sense of \cite{collection}. Then for any admissible subcategory $i:D\sr C$ of $C$ such that $D\subset Sm/k$ one has
$$\LL i^{rad}_{tr}(M(A,p,q;F)_D)=M(A,p,q;F)_C$$
\end{theorem}
\begin{proof}
Set $K_C=K(A,p,q)_C$ and $K_D=K(A,p,q)_D$ and let $M_C=\LL\Lambda^l_F(K_C)$ and $M_D=\LL\Lambda^l_F(K_C)$ be  the corresponding objects of $H_{Nis,\af}(Cor(-,F))$.  By Section \ref{changeofC} we have two pairs of adjoint functors connecting the categories $H_{Nis,\af}(Cor(-))$ and $H_{Nis,\af}((-)_+)$ over $C$ and $D$ respectively.  By Lemma \ref{2009.09.17.l2}, for any $i:D\sr C$, any $A$ and any $p,q$ one has 
\begin{eq}
\llabel{2009.09.15.eq2}
K_D=i_{rad,+}K_C
\end{eq}
Since $k$ be is a field with resolution of singularities and $D\subset Sm/k$ we further have
$$i_{tr,rad}M_C=i_{tr,rad}\LL\Lambda^l_F K_C=\LL\Lambda^l_F i_{rad,+} K_C=\LL \Lambda^l_F K_D=M_D$$
where the second equality holds by Theorem \ref{ressing} and therefore
$$\LL i^{rad}_{tr}(M_D)=\LL i^{rad}_{tr}i_{tr,rad}M_C$$
Since $p\ge 2q$, Corollary \ref{2009.09.15.cor1} implies that $M_C$ is in the image of the functor $\LL i^{rad}_{tr}$ and by Corollary \ref{2009.07.24.4} this implies that the adjunction $\LL i^{rad}_{tr}i_{tr,rad}M_C\sr M_C$ is an isomorphism.
\end{proof}
\begin{cor}
\llabel{2009.09.15.cor2}
Under the assumptions of Theorem \ref{2009.09.12.th1} assume in addition that $R=F$ is a field, $p\ge 2q$ and that $k$ admits resolution of singularities in the sense of \cite{collection}. Then for any admissible subcategory $D\subset Sm/k$ one has
$$M(A,p,q;F)_D\subset \overline{SPT}_{\ge q}$$
\end{cor}

\subsection{Topological realization functors}
\llabel{top.real}
In this section we assume that the base field is $\bf C$ and $R$ is the ring of coefficients for homology and correspondences. We set $C=QP/{\bf C}$ to be the category of quasi-projective schemes over $\bf C$.
 
Let $Top$ be the category of (all) topological spaces. Sending $X\in C$ to the topological space of its $\bf C$-points we get a functor
$$\pi:C\sr Top$$
We will need the following classical properties of this functor.
\begin{theorem}
\llabel{prop}\begin{enumerate}
\item $\pi$ commutes with disjoint unions,
\item $\pi$ commutes with fiber products,
\item $\pi$ commutes with finite group quotients,
\item $\pi$ takes universal homeomorphisms to homeomorphisms. 
\end{enumerate}
\end{theorem}
\begin{proof}
See \cite{GAGA}.
\end{proof}
\begin{proposition}
\llabel{2009.09.16.pr1}
For any closed embedding $i:Z\sr X$ in $QP/{\bf C}$, $\pi(i)$ is a closed embedding and $(X({\bf C}),Z({\bf C})))$ is a $CW$-pair.
\end{proposition}
\begin{proof}
It follows easily from the considerations of \cite{Hironaka2}.
\end{proof}

Consider the functor 
$$\pi^{rad}_{+}:\Delta^{op}Rad(C_+)\sr \Delta^{op}Rad(Top_+)$$
\begin{lemma}
\llabel{2009.09.16.l1}
The functor $\pi_{+}^{rad}$ takes projective equivalences to projective equivalences. In particular, for any $X\in \Delta^{op}rad(C_+)$ the morphism $\LL \pi^{rad}_{+}(X)\sr \pi^{rad}_+(X)$ is a projective equivalence.
\end{lemma}
\begin{proof}
A morphism of radditive functors is a projective equivalence if and only if it is a simplicial weak equivalence as a morphism of presheaves of sets. Such equivalences are preserved by inverse image functors defined by functors which commute with fiber products. Therefore the lemma follows from Theorem \ref{prop}(2).
\end{proof}

Let $\Lambda^r_{R,mod}$ be the "forgetting of transfers" functor from  $Rad(Cor(C,R))$ to presheaves of $R$-modules on $C$. Then  $\pi^{*}_{R-mod}\Lambda^r_{R,mod}(X)$ is a functor 
$$Rad(Cor(C,R))\sr PreShv_{R-mod}(Top)$$
Evaluating a pointed presheaf (resp. a presheaf of $R$-modules) $F$ on the standard cosimplicial object $\Delta^{\BB}_{top}$ in $Top$ we get a pointed simplicial set (resp. a simplicial $R$-module) $Sing_*(F)$. Composing $\pi^{rad}_+$ (resp. $\pi^*_{R-mod}\Lambda^r_{R,mod}$) with $Sing_*$ followed by the diagonal we get two functors
$$T_{\bf C}=\Delta Sing_*\pi^{rad}_+:\Delta^{op}Rad(C_+)\sr \Delta^{op}Sets_{\BB}$$
and
$$T_{\bf C}^{tr}=\Delta Sing_*\pi^*_{R-mod}\Lambda^r_{R,mod}:\Delta^{op}Rad(Cor(C,R))\sr \Delta^{op}R-mod$$
\begin{proposition}
\llabel{tcbb}
The functor $T_{\bf C}$ takes $(Nis,\af)$-equivalences to weak equivalences of simplicial sets and therefore defines a functor
$$t_{\bf C}:H_{Nis,\af}(C_+)\sr H_{\BB}^{top}$$
\end{proposition}
\begin{proof}
By Theorem \ref{maindelta} the class 
$$W_{Nis,\af}^{+}=cl_l((G_{Nis})_+\cup (G_{\af})_+)$$
of $(Nis,\af)$-equivalences in $\Delta^{op}rad(C_+)$ coincides with the class $cl_{\bdl}((G_{Nis}\amalg Id_C)_+\cup (G_{\af}\amalg Id_C)_+\cup W_{proj})$. Since $T_{\bf C}$ clearly takes $\bdl$-closures to $\bdl$-closures and commutes with $\amalg$ and since the class of weak equivalences in $\Delta^{op}Sets$ is $\bdl$-closed it is sufficient to verify that $T_{\bf C}$ takes $W_{proj}$, $(G_{Nis})_+$ and $(G_{\af})+$ to weak equivalences. The case of $W_{proj}$ follows from Lemma \ref{2009.09.16.l1}. Let us consider the case of $G_{Nis}$. Let $Q$ be a Cartesian square in $C$ of the form
$$
\begin{CD}
B @>>> Y\\
@VVV @VVpV\\
A @>j>> X
\end{CD}
$$
such that $p$ is etale, $j$ an open embedding and $Y\backslash B\sr X\backslash A$ is an isomorphism. Consider the morphism $f_+:(K_Q)_+\sr X_+$. For $U\in QP/{\bf C}$ we have $\pi^*(U)=U({\bf C})$ and $T_{\bf C}(U)=Sing_*(U({\bf C}))$ where $Sing_*$ is a singular simplicial set of a topological space. Therefore,  there is a push-out square of simplicial sets of the form
$$
\begin{CD}
Sing_*(B({\bf C}))_+\vee Sing_*(B({\bf C}))_+ @>>> Sing_*(Y({\bf C}))_+\times \Delta^1\\
@VVV @VVV\\
Sing_*(A({\bf C}))_+\vee Sing_*(Y({\bf C}))_+ @>>> T_{\bf C}((K_Q)_+)
\end{CD}
$$
and we need to verify that the obvious morphism $T_{\bf C}((K_Q)_+)\sr Sing_*(X({\bf C}))_+$ is a weak equivalence or equivalently that the square of singular simplicial sets 
\begin{eq}
\llabel{2009.12.19.oldeq1}
\begin{CD}
Sing_*(B({\bf C}))_+ @>>> Sing_*(Y({\bf C}))_+\\
@VVV @VVV\\
Sing_*(A({\bf C}))_+ @>>> Sing_*(X({\bf C}))_+
\end{CD}
\end{eq}
is a homotopy push-out square. This follows from Lemma \ref{2009.12.19.l2}.

It remains to consider the case of $(G_{\af})_+$ i.e. to show that for any $X\in QP/{\bf C}$ the map 
$$Sing_*((X\times \af)({\bf C}))_+\sr Sing_*(X({\bf C})_+)$$
is a weak equivalence. This follows from the fact that $(X\times \af)({\bf C})=X({\bf C})\times {\bf C}$ and ${\bf C}$ is contractible.
\end{proof}
\begin{lemma}
\llabel{2009.12.19.l2}
Let $Q$ 
$$
\begin{CD}
B @>>> Y\\
@VVV @VVV\\
A @>>> X
\end{CD}
$$
be an upper distinguished square in $Sch/{\bf C}$. Then the associated square of the form (\ref{2009.12.19.oldeq1}) is a homotopy push-out square.
\end{lemma}
\begin{proof}
Since the spaces of $\bf C$-points of $A$, $B$, $X$ and $Y$ admit triangulations the condition that the square (\ref{2009.12.19.oldeq1}) is homotopy push-out square is equivalent to the condition that the obvious map $q:hp(Q({\bf C}))\sr X({\bf C})$ where $hp(Q({\bf C}))$ is the push-out  of the diagram
$$
\begin{CD}
B({\bf C})\times \partial\Delta^1 @>>> A({\bf C})\amalg Y({\bf C})\\
@VVV\\
B({\bf C})\times\Delta^1
\end{CD}
$$
is a weak equivalence. By \cite[Cor. 1.4, p.93]{Mayqf} (see also \cite{Dugger}), it is sufficient to show that for any point $x\in X({\bf C})$ there exists an open neighborhood $U'$ of $p$ such that for any open subset $U$ of $U'$ the map $q^{-1}(U)\sr U$ is a weak equivalence. This map is clearly isomorphic to the map $q_U:hp(Q({\bf C})_U)\sr U$ where $Q({\bf C})_U$ is the pull-back of the square $Q({\bf C})$ to a square over $U$. 

Using the fact that etale morphisms define local homeomorphisms on the spaces of ${\bf C}$-points one can show that for any $x\in X({\bf C})$ there exists $U'$ such that $Q({\bf C})_{U'}$ is isomorphic to a square of the form
$$
\begin{CD}
(\amalg_i (U'-Z({\bf C})\cap U'))\amalg (U'-Z({\bf C})\cap U') @>>> (\amalg_i (U'-Z({\bf C})\cap U'))\amalg U'\\
@VVV @VVV\\
U'-Z({\bf C})\cap U' @>>> U'
\end{CD}
$$
Then the same is true for $Q({\bf C})_U$ for any $U\subset U'$ which easily implies that the maps $q_U$ are homotopy equivalences.
\end{proof}
\begin{lemma}
\llabel{tcsmash}
Functor $t_{\bf C}$  commutes with the smash products.
\end{lemma}
\begin{proof}
One can easily see that for any $X,Y\in \Delta^{op}Rad(C_+)$ there is a natural map $T_{\bf C}(X)\wedge T_{\bf C}(Y)\sr T_{\bf C}(X\wedge Y)$. Since the class of weak equivalences of pointed simplicial sets is $\bdl$-closed it is sufficient to check that it is a weak equivalence for $X=X'_+$, $Y=Y'_+$ where $X',Y'\in C$. In this case our map is an isomorphism because both $\pi^*$ and $Sing_*$ commute with direct products.
\end{proof}
One can easily see that there are canonical isomorphisms:
\begin{eq}
\llabel{isos}
t_{\bf C}(S^1_s)=S^1
\end{eq}
\begin{eq}\llabel{isot}
t_{\bf C}(S^1_t)=S^1
\end{eq}
Let $H(R-mod)$ be the homotopy category of simplicial $R$-modules.
\begin{proposition}
\llabel{tctr}
The functor $T^{tr}_{\bf C}$ takes $(Nis,\af)$-equivalences to weak equivalences of simplicial $R$-modules and therefore defines a functor
\begin{eq}
\llabel{tctrh}
t^{tr}_{\bf C}:H_{Nis,\af}(Cor(C,R))\sr H(R-mod)
\end{eq}
\end{proposition}
\begin{proof}
Let $\phi$ be the forgetting functor from $R$-modules to sets. Then by construction we have
$$\phi\circ T^{tr}_{\bf C}=T_{\bf C}\Lambda^r_R$$
The functor $\phi$ reflects weak equivalences. The functor $\Lambda^r_R$ respects $(Nis,\af)$-equivalences by Theorem \ref{2009th1}. The functor $T_{\bf C}$ respects equivalences by Proposition \ref{tcbb}. We conclude that $T^{tr}_{\bf C}$ respects equivalences.
\end{proof}
Let 
$$H_R:\Delta^{op}Sets_{\BB}\sr \Delta^{op}R-mod$$ 
be the functor which takes a pointed simplicial set $X$ to the free $R$-module $H_R(X)$ generated $X$.
\begin{proposition}
\llabel{stcom}
The square
\begin{eq}
\begin{CD}
H_{Nis,\af}(C_+) @>t_{\bf C}>> H_{\BB}^{top}\\
@V\LL\Lambda^l_RVV @VH_RVV\\
H_{Nis,\af}(Cor(C,R)) @>t_{\bf C}^{tr}>> H(R-mod)
\end{CD}
\end{eq}
commutes up to a natural isomorphism.
\end{proposition}
\begin{proof}
We may interpret $H_{Nis,\af}(C_+)$ as a localization of $\Delta^{op}C_+^{\#}$ and consider $\Lambda^l$ instead of $\LL\Lambda^l$.  By definition, $H_R T_{\bf C}(X)$ is the free $R$-module generated by the pointed simplicial set $\Delta Sing_{*}(X({\bf C}))$ and $T_{\bf C}^{tr} \Lambda_R^l(X)$ is the simplicial $R$-module $\Delta Sing_*(\pi^*\Lambda^r_{R,ab}\Lambda^l_R(X))$. The natural transformation $X({\bf C})\sr \pi^*\Lambda^r_{R,ab}\Lambda^l_R(X)$ together with the universal property of free $R$-modules provide us with a natural transformation of the form
\begin{eq}
\llabel{neq3}
H_R T_{\bf C}\sr T_{\bf C}^{tr} \LL\Lambda_R^l
\end{eq}
Since the forgetting functor $\psi$ from $R$-modules to abelian groups reflects equivalences and one has
$$\psi H_R T_{\bf C}(X)=(\psi H_{\zz} T_{\bf C}(X))\oo R$$
and 
$$\psi T_{\bf C}^{tr} \Lambda_R^l(X)=(\psi T_{\bf C}^{tr} \Lambda_{\zz}^l(X))\oo R$$
its is enough to consider the case $R=\zz$. Since both sides of (\ref{neq3}) commute with the simplicial suspension (because of Theorem \ref{prop}(1)) and the suspension on $H(\zz-mod)$ reflects isomorphisms, it is enough to verify that for $X\in \Delta^{op}C_{+}^{\#}$ the morphism 
$$\phi H_{\zz} \Delta Sing_*\pi^{rad}_+(\Sigma^1_sX)\sr \Delta Sing^*\pi^{rad}_+\Lambda^r_{\zz}\Lambda^l_{\zz}(\Sigma^1_sX)$$
is a weak equivalence. Consider the commutative square
\begin{eq}
\llabel{2009.09.17.eq1}
\begin{CD}
\phi H_{\nn} \Delta Sing_*\pi^{rad}_+(\Sigma^1_sX)@>>> \Delta Sing^*\pi^{rad}_+S^{\infty}(\Sigma^1_sX)\\
@VVV @VVV\\
\phi H_{\zz} \Delta Sing_*\pi^{rad}_+(\Sigma^1_sX)@>>> \Delta Sing^*\pi^{rad}_+\Lambda^r_{\zz}\Lambda^l_{\zz}(\Sigma^1_sX)
\end{CD}
\end{eq}
where $H_{\nn}(-)$ is the free abelian monoid functor. The left hand side vertical arrow in this square is a weak equivalence by Lemma \ref{ge64}. 

Let $i:SN/{\bf C}\sr QP/{\bf C}$ be the embedding of the subcategory of semi-normal schemes. By Lemma \ref{2009.09.17.l1} we have $\pi^{rad}_+=\pi^{rad}_+\LL i^{rad}_+i_{rad,+}$. 
The morphism $i_{rad,+}S^{\infty}(\Sigma^1_sX)\sr i_{rad,+}\Lambda^r_{\zz}\Lambda^l_{\zz}(\Sigma^1_sX)$ is a $Nis$-equivalence by Proposition \ref{susvoeplus} and Lemma \ref{hens}. 
Therefore $\LL i^{rad,+}i_{rad,+}S^{\infty}(\Sigma^1_sX)\sr \LL i^{rad,+}i_{rad,+}\Lambda^r_{\zz}\Lambda^l_{\zz}(\Sigma^1_sX)$ is a $Nis$-equivalence and we conclude by Proposition \ref{tcbb} that the right hand side vertical arrow of (\ref{2009.09.17.eq1}) is a weak equivalence. 

It remains to verify that the morphism
\begin{eq}
\llabel{2009.09.16.eq1}
\phi H_{\nn} \Delta Sing_*\pi^{rad}_+(\Sigma^1sX)\sr \Delta Sing_*\pi^{rad}_+S^{\infty}(\Sigma^1_sX)
\end{eq}
is a weak equivalence for any $X\in \Delta^{op}(QP/{\bf C})_+^{\#}$. Since the class of weak equivalences of simplicial sets is $\bdl$-closed it is sufficient to consider the case of $X\in (QP/{\bf C})_+$. By Theorem \ref{prop}(1,2,3) we conclude that for such an $X$ one has 
$$\pi^{rad}_+(S^{\infty}(\Sigma^1_sX))=S^{\infty}(\Sigma^1_s(X({\bf C})))$$
where on the right $\Sigma^1_s(X({\bf C}))$ is considered as a simplicial topological space. 
The fact that  the morphism (\ref{2009.09.16.eq1}) is a weak equivalence follows now from Proposition \ref{2009.09.16.pr1} and  the Dold-Thom theorem in the form which asserts that for a  simplicial topological space $X$ whose terms are $CW$-complexes, the maps
$$H_{\nn} \Delta Sing_*(\Sigma^1_s X)\sr \Delta Sing_*(S^{\infty}(\Sigma^1_s X))$$
is a weak equivalences.
\end{proof}
\begin{lemma}
\llabel{2009.09.17.l1}
Let $i:SN/{\bf C}\sr QP/{\bf C}$ be the embedding of semi-normal schemes. Then the morphism $\pi^{rad}_+\LL i^{rad}_+ i_{rad,+}\sr \pi^{rad}_+$ is a projective equivalence.
\end{lemma}
\begin{proof}
It is clearly sufficient to consider the non-pointed case. Let $sn:QP/{\bf C}\sr SN/{\bf C}$ be the semi-normalization functor which is the left adjoint to $i$. Then the adjunction between $sn$ and $i$ shows that there is a natural isomorphism $i_{rad}=sn^{rad}$. Since $i_{rad}$ respects projective equivalence so does $sn^{rad}$ and therefore $i_{rad}=\LL sn^{rad}$. By Lemma \ref{2009.09.16.l1} we have $\pi^{rad}=\LL\pi^{rad}$. Combing these observations together we get:
$$\pi^{rad} \LL i^{rad} i_{rad}=\LL\pi^{rad}\LL i^{rad}\LL sn^{rad}=\LL(\pi\circ i\circ sn)^{rad}=\LL\pi^{rad}$$
The last equality holds by Theorem \ref{prop}(4) since the morphism $i\circ sn\sr Id$ is a universal homeomorphism.
\end{proof}

\begin{proposition}
\llabel{tctroo}
The functor $t_{\bf C}^{tr}$ commutes with (derived) tensor products.
\end{proposition}
\begin{proof}
Note first that there is a natural transformation on the level of simplicial objects of the form
$$T_{\bf C}^{tr}(X)\oo_{R} T_{\bf C}^{tr}(Y)\sr T_{\bf C}^{tr}(X\oo Y)$$
Since the class of weak equivalences of simplicial $R$-modules is $\bdl$-closed 
it is enough to verify that it is a weak equivalence for $X,Y\in Cor(C,R)$ i.e. for $X=\Lambda^l(X'_+)$, $Y=\Lambda^l(Y'_+)$ where $X',Y'\in C$. In this case the claim follows from Proposition  \ref{stcom} and the K\"unneth isomorphism theorem for topological homology.
\end{proof}
Proposition \ref{stcom}, isomorphisms (\ref{isos}), (\ref{isot}) and Proposition \ref{tctroo} imply that there are a canonical isomorphism 
\begin{eq}
\llabel{isoln}
t_{\bf C}^{tr}(A\oo_{\LL,\zz}l_n)=(A\oo_{\LL,\zz} R)[n]
\end{eq}
and in particular $t_{\bf C}^{tr}(L_n)=R[2n]$. Using the definition of motivic cohomology given at the beginning of Section \ref{moore.section} in combination with Proposition \ref{stcom}, isomorphisms (\ref{isoln}) for $R=\zz$ and the suspension isomorphisms in the topological cohomology we get canonical maps:
$$\tilde{H}^{p,q}_{un}(X,A)\sr \tilde{H}^p(t_{\bf C}(X),A)$$
For an $R$-module $A$ let $K(A,p)$ and $K(A,p,q)_{QP/{\bf C}}$
be the topological and the motivic Eilenberg-MacLane spaces representing the functors $\tilde{H}^{p}(-,A)$ and $\tilde{H}^{p,q}_{un}(-,A)$ respectively. 
\begin{lemma}
\llabel{twoks}
For $p\ge q$ there is a canonical isomorphism 
$$t_{\bf C}(K(A,p,q)_{QP/{\bf C}})\sr K(A,p).$$ 
\end{lemma}
\begin{proof}
We have
$$t_{\bf C}(K(A,p,q))=t_{\bf C}\Lambda^r_{\zz}(A\oo_{\LL,\zz} \Sigma^{p-q}l_q)=\phi t_{\bf C}^{tr}(A\oo_{\LL,\zz} \Sigma^{p-q}l_q)=\phi(A[p])=K(A,p).$$
\end{proof}
\begin{remark}\rm
\llabel{plessq}
The statement of Lemma \ref{twoks} is false for $p<q$ at least when $A$ is not a torsion abelian group. For example, $K(\zz,0,1)=pt$ while $K(\zz,0)=\zz$.
\end{remark}
\begin{remark}\rm
Note that in the case of $R=\zz/l$ we get
$t_{\bf C}^{tr}(L_n)=\zz/l[2n]$ while it would be more natural to have $t_{\bf C}^{tr}(L_n)=\mu_l[2n]$ where $\mu_l$ is the group of $l$-roots of unity in $\bf C$. The reason that we get $\zz/l$ instead of $\mu_l$ is that  the isomorphism (\ref{isot}) defines an identification of $\mu_l$, which is the fiber of the $l$-th power map on $S^1_t$, with $\zz/l$ which is the fiber of the $l$-th power map on the circle. If the circle is oriented counter clock-wise then this identification corresponds to the choice of the $l$-th root of unity with the smallest argument. 
\end{remark}
In the next section we will need to work with $K(A,p,q)_{SmQP/{\bf C}}$ instead of $K(A,p,q)_{QP/{\bf C}}$. Let $j:SmQP/{\bf C}\sr QP/{\bf C}$ be the embedding of smooth quasi-projective schemes to all quasi-projective schemes. By Lemma \ref{2009.09.17.l2} we have an equivalence
$$j_{rad,+} K(A,p,q)_{QP/{\bf C}}=K(A,p,q)_{SmQP/{\bf C}}$$
which, by adjunction, defines a morphism
$$\psi(A,p,q):\LL i^{rad}_+K(A,p,q)_{SmQP/{\bf C}}\sr K(A,p,q)_{QP/{\bf C}}$$
\begin{proposition}
\llabel{2009.09.17.pr1}
For $p\ge q$ the morphism $t_{\bf C}(\psi(A,p,q))$ is a weak equivalence.
\end{proposition}
\begin{proof}
Note first that the category $SmQP/{\bf C}$ has direct products (fiber products over the point) and that the functor $j$ respects these products. This implies that $\LL j^{rad}_+$ commutes with direct products up to a canonical equivalence. Therefore both  $t_{\bf C}(\LL j^{rad}_+(K(A,p,q)_{SmQP/{\bf C}}))$ and $t_{\bf C}(K(A,p,q)_{QP/{\bf C}})$ are $H$-spaces with an inverse map. A morphism of such $H$-spaces is a weak equivalence if it defines a weak equivalence on homology i.e. it is sufficient to show that 
$$H_{\zz}t_{\bf C}(\LL j^{rad}_+(K(A,p,q)_{SmQP/{\bf C}}))\sr H_{\zz}t_{\bf C}(K(A,p,q)_{QP/{\bf C}})$$
is a weak equivalence of simplicial abelian groups.  By Proposition \ref{stcom} this is equivalent to the showing that 
\begin{eq}
\llabel{2009.09.17.eq3}
t^{tr}_{\bf C}(\LL\Lambda^l_{\zz}\LL j^{rad}_+(K(A,p,q)_{SmQP/{\bf C}}))\sr t^{tr}_{\bf C}(\LL\Lambda^l_{\zz}(K(A,p,q)_{QP/{\bf C}}))
\end{eq}
is an isomorphism. For $q\le 1$ there are models of $K(A,p,q)_{QP/{\bf C}}$ which lie in $\Delta^{op}(SmQP/{\bf C})_+$ and therefore $\psi(A,p,q)$ itself is a weak equivalence. Therefore we may assume that $p>0$. Then (\ref{2009.09.17.eq3}) is an isomorphism by Theorem \ref{mth1} since $\LL\Lambda^l_{\zz}\LL j^{rad}_+=\LL j^{rad}_{tr}\LL\Lambda^l_{\zz}$. 
\end{proof}
\begin{cor}
\llabel{2009.09.17.cor1}
Let 
$$t_{Sm,{\bf C}}=\Delta Sing_* \LL(\pi\circ j)^{rad}_+$$
be the topological realization functor on $H_{Nis,\af}((SmQP/{\bf C})_+)$. Then for $p\ge q$, there 
are canonical weak equivalences
$$t_{Sm,{\bf C}}(K(A,p,q)_{SmQP/{\bf C}})=t_{\bf C}(K(A,p,q)_{QP/{\bf C}})=K(A,p).$$
\end{cor}

\subsection{Application to stable operations.}
In this section we will show that over a field $k$ of characteristic zero the algebra of all stable operations in motivic cohomology with coefficients in ${\bf F}_l$ coincides with the motivic Steenrod algebra. Since motivic Steenrod operations have been defined only in cohomology of smooth schemes we will work with smooth schemes over $k$.  Let $H^{*,*}=H^{*,*}(Spec(k),{\bf F}_l)$.  Let $K_n$ be the motivic Eilenberg-MacLane space $K(\zz/l,2n,n)_{Sm/k}$ and $K_{2n}^{top}$ be its topological counterpart $K(\zz/l,2n)$. The abelian group of all stable operations is given by
$$M^{*,*}={\lim_{\leftarrow}}{}_{n} \,\, H^{*+2n,*+n}(K_n,{\bf F}_l)$$
where the homomorphisms of the system are defined by the morphisms
\begin{eq}
\llabel{mor1}\Sigma^1_TK_n\sr K_{n+1}
\end{eq}
Let ${\cal A}^{*,*}$ be the motivic Steenrod algebra. Since operations from ${\cal A}^{*,*}$ are stable with respect to $\Sigma_T$ they act on $M^{*,*}$. Denote by $\iota$ the element in $M^{*,*}$ whose restriction to $K_n$ is the canonical class $\iota_n$. Acting by elements of ${\cal A}^{*,*}$ on $\iota$ we get a map
$$u:{\cal A}^{*,*}\sr M^{*,*}$$
which is a homomorphism of left $H^{*,*}$-modules. 
\begin{theorem}
\llabel{mainst}
Let $k$ be a field of characteristic zero. Then $u$ is an isomorphism.
\end{theorem}
\begin{proof}
The proof of this theorem occupies the rest of this section and ends right before Remark \ref{false}.
\end{proof}
Let $M_n=\LL \Lambda^l_{\zz}(K_n)$ be the class of $K_n$ in $H_{Nis,\af}(Cor(Sm/k,{\bf F}_l))$. In view of Corollary \ref{2009.09.15.cor2} there exist objects $M'_n$ in $\overline{SPT}$ such that $M_n=\Sigma_T^n M'_n$. Morphisms  (\ref{mor1}) define a sequence
$$M'_0\sr M'_1\sr M'_2\sr\dots$$
Let $M$ be its homotopy colimit. By Corollary \ref{hocolim} the distinguished triangle which defines $M$ splits. Therefore the motivic cohomology of $M$ coincide with $M^{*,*}$ and by the same corollary $M\in \overline{SPT}$.

Let $(e_{\alpha})_{\alpha\in A}$ be the basis of ${\cal A}^{*,*}(k)$ over $H^{*,*}$ which consists of admissible monomials (see \cite[pp.40-41]{Redpub}).  Let $p_{\alpha},q_{\alpha}$ be defined by the condition that $e_{\alpha}\in {\cal A}^{p_{\alpha},q_{\alpha}}$. Set ${\cal A}=\oplus_{\alpha\in A}{\bf F}_l(q_{\alpha})[p_{\alpha}]$. The bi-degrees of admissible monomials are such that $p_{\alpha}\ge 2q_{\alpha}$ i.e. ${\cal A}\in \overline{SPT}$.  
\begin{lemma}
\llabel{2009.09.18.l1}
The motivic cohomology of ${\cal A}$ is canonically isomorphic to ${\cal A}^{*,*}$ and there is a unique morphism 
$$\tilde{u}:M\sr {\cal A}$$
which defines on the motivic cohomology homomorphism $u$. 
\end{lemma}
\begin{proof}
The definition of an admissible monomial given in \cite[p.40]{Redpub} implies immediately that for any $q\ge 0$ the set $A_{*,\le q}$ of $\alpha\in A$ such that $q_{\alpha}\le q$ is finite. Let $H^{*,*}(q)[p]$ be the free bi-graded module over $H^{*,*}$ with the generator in bi-degree $(p,q)$. Then we have 
$$H^{*,*}({\cal A})=\prod_{\alpha\in A}^{gr} H^{*,*}(q_{\alpha})[p_{\alpha}]=\oplus_{\alpha\in A}^{gr} H^{*,*}(q_{\alpha})[p_{\alpha}]={\cal A}^{*,*}$$
where $\prod^{gr}$ and $\oplus^{gr}$ are the direct product and the direct sum in the category of bi-graded modules and the middle equality holds because the sets $A_{*,\le q}$ are finite. This proves the first assertion of the lemma. 

The condition that $\tilde{u}$ defines homomorphism $u$ on motivic cohomology is equivalent to the condition that the composition of $u$ with the projection to the summand corresponding to $e_{\alpha}$ is $u(e_{\alpha})$. The uniqueness and existence of such $\tilde{u}$ follows from the finiteness of the sets $A_{*,\le q}$ and Proposition \ref{2009.09.18.pr1}.
\end{proof}
To prove Theorem \ref{mainst} we need to show that $\tilde{u}$ is an isomorphism. Since both $M$ and ${\cal A}$ belong to $\overline{SPT}$ we may apply Corollary \ref{2009.09.17.cor3}(2) and assume from now on that $k$ is an algebraically closed field.  

Let us choose a primitive $l$-th root of unity in $k$ and let $\tau$ be the corresponding element of $H^{0,1}=H^{0,1}(Spec(k))$. Then  by \cite[Cor. 4.3, p.254]{Suslin.collection} we have
\begin{eq}
\llabel{2009.09.18.eq3}
H^{*,*}={\bf F}_l[\tau].
\end{eq}
For a module $A$ over $F[\tau]$ we will write $A\oo_0F$ (resp. $A\oo_1F$) for the tensor product of $A$ with ${\bf F}_l$ with respect to the homomorphism $F[\tau]\sr F$ which takes $\tau$ to $0$ (resp. to $1$).

\begin{lemma}
\llabel{alg2}
$u:M^*\sr N^*$ a homomorphism of non-negatively graded free modules over $F[\tau]$ where $F$ is a field and $gr(\tau)=1$. Assume that $u_1=u\oo_1 Id$ is surjective and $u_0=u\oo_0 Id$ injective. Then $u$ is an isomorphism.
\end{lemma}
\begin{proof}
If $x\in M^n$ and $u(x)=0$ then the image of $x$ in $M^*\oo_{0} F$ is zero i.e. $x=\tau x'$. Then $0=u(\tau x')=\tau u(x')$ and since  $N^*$ is free we have $u(x')=0$. A simple induction on $n$ implies now that $u$ is injective.

Let us show that $u$ is surjective. Let $x\in N^{n}$. Since $u_1$ is surjective we have $x=u(z)+(\tau-1)y$ for some $z\in M^*$ and $y\in N^*$. Let $m$ be the smallest integer $\ge n$ such that $z\in M^{\le m}$ and $y\in N^{\le m}$.  Then $x_{m+1}=0$ and $z_{m+1}=0$ and therefore $\tau y_m =0$. Since $N^*$ is free we conclude that $y_m=0$. Suppose that $m>n$. Then $u(z_m)=-\tau y_{m-1}$. 

The condition that $u_0$ is a monomorphism shows that if $u(z_m)=-\tau y_{m-1}$ then there exists $z'_{m-1}$ such that $z_m=\tau z'_{m-1}$ and, since $N^*$ is a free module $u(z'_{m-1})=-y_{m-1}$. Set $z''=z-z_m+z'_{m-1}$ and $y''=y-y_{m-1}$. Then one has again $x=u(z'')+(\tau-1)y''$. By induction we conclude that we may assume that $m=n$. Then $(\tau-1)y=x-u(z)$ implies that $y_{n}=0$ and $x=u(z_n)+\tau y_{n-1}$. By obvious induction on degree we may assume that $y_{n-1}$ is in the image of $u$. Then $x\in Im(u)$. 
\end{proof}
\begin{lemma}
\llabel{caseofst}
Let $X\in \Delta^{op}Rad((QP/{\bf C})_+)$ be such that $\LL\Lambda^l_{{\bf F}_l}(X)$ is a direct sum of objects of the form $\Sigma^il_j$. Then the natural homomorphisms
\begin{eq}
\llabel{oo1iso}
\tilde{H}_{un}^{p,*}(X,\zz/l)\oo_1{\bf F}_l\sr \tilde{H}^p(t_{\bf C}(X),{\bf F}_l)
\end{eq}
are isomorphisms.
\end{lemma}
\begin{proof}
By Proposition \ref{stcom} we have
$$\tilde{H}^p(t_{\bf C}(X),{\bf F}_l)=Hom_{H({\bf F}_l-mod)}(H_{{\bf F}_l}t_{\bf C}(X),{\bf F}_l[p])=Hom(t_{\bf C}^{tr}\LL\Lambda^l_{{\bf F}_l}(X), {\bf F}_l[p])$$
and the homomorphism $\tilde{H}_{un}^{p,*}(X,\zz/l)\sr \tilde{H}^p(t_{\bf C}(X),{\bf F}_l)$ which defines (\ref{oo1iso})
is the homomorphism
\begin{eq}
\llabel{oo2iso}
\begin{array}{c}
(\oplus_{q>p}Hom_{H_{Nis,\af}}(\Sigma^{q-p}Y,{\bf F}_l(q)[q]))\oplus(\oplus_{p\ge q\ge 0}Hom_{H_{Nis,\af}}(Y,{\bf F}_l(q)[p]))\\
\downarrow\\
 Hom_{H({\bf F}_l-mod)}(t_{\bf C}^{tr}(Y),{\bf F}_l[p])
 \end{array}
\end{eq}
for $Y=\LL\Lambda^l_{{\bf F}_l}(X)$. Therefore, it is sufficient to verify that (\ref{oo2iso}) defines an isomorphism for $Y=\Sigma^il_j$ i.e. that the maps
$$\tilde{H}^{p-i-j,*-j}_{un}(Spec(k)_+,\zz/l)\oo_1{\bf F}_l\sr \tilde{H}^p_{un}(S^{i+j},{\bf F}_l)$$
are isomorphisms. This is equivalent to (\ref{2009.09.18.eq3}).
\end{proof}
The motivic Adem relations demonstrated in \cite{Redpub} imply that for any $k$ of characteristic zero there is a homomorphism 
\begin{eq}\llabel{assas}
{\cal A}^{*,*}(k)\sr {\cal A}^*
\end{eq}
which sends $P^i$ to $P^i$ and $\beta$ to $\beta$ and which is of the form
$$H^{*,*}(k)\sr H^{*,*}(\bar{k})={\bf F}_l[\tau]\stackrel{\tau\mapsto 1}{\sr} {\bf F}_l$$
on $H^{*,*}$. The same relations imply immediately the following result.
\begin{lemma}
\llabel{assasiso}
Let $k$ be an algebraically closed field of characteristic $0$. Then (\ref{assas}) defines an isomorphism
\begin{eq}\llabel{assasisoeq}
{\cal A}^{*,*}(k)\oo_1{\bf F}_l\sr {\cal A}^*
\end{eq}
\end{lemma}
\begin{lemma}
\llabel{useunique}
The square
\begin{eq}
\llabel{comsq1}
\begin{CD}
{\cal A}^{*,*} @>\iota_n>>  \tilde{H}^{*,*}(K_n)\\
@VVV  @VVV\\
{\cal A}^* @>\iota_{2n}>> \tilde{H}^{*}(K_{2n}^{top})
\end{CD}
\end{eq}
commutes. 
\end{lemma}
\begin{proof}
We will only give a sketch of the argument. From general functoriality it is sufficient to verify that the image $P^{i,n}$ (resp. $\beta P^{i,n}$) of the motivic class $P^i\iota_n$ (resp. $\beta P^i\iota_n$) in $\tilde{H}^*(K_{2n}^{top})$ is $P^i\iota_{2n}$ (resp. $\beta P^i\iota_{2n}$). Knowing the Cartan formula and Adem relations for the motivic reduced powers we can deduce that the family of operations defined by the classes $P^{i,n}$ and $\beta P^{i,n}$ satisfy the list of properties which uniquely characterize the reduced power operations (see e.g. \cite{SE}).
\end{proof}

\begin{proposition}
\llabel{topu}
Let $k$ be an algebraically closed field of characteristic zero.  Then the  homomorphism 
\begin{eq}
\llabel{topueq}
u\oo_1 Id:{\cal A}^{*,*}\oo_{1}{\bf F}_l\sr M^{*,*}\oo_1{\bf F}_l
\end{eq}
is an isomorphism.
\end{proposition}
\begin{proof}
Since both sides of (\ref{topueq}) remain unchanged when we pass from an algebraically closed field to its algebraically closed extension we may assume that $k={\bf C}$. 

The vertical arrows of (\ref{comsq1}) factor as
$${\cal A}^{*,*}\sr {\cal A}^{*,*}\oo_1 {\bf F}_l\sr {\cal A}^*$$
$$\tilde{H}^{*,*}(K_n)\sr \tilde{H}^{*,*}(K_n)\oo_1{\bf F}_l\sr \tilde{H}^*(K_{2n}^{top})$$
where the second arrows of both factorizations are isomorphisms - the first one by Lemma \ref{assasiso} and the second one by Lemma \ref{caseofst}. We conclude that the maps
$${\cal A}^{*,*}\oo_1{\bf F}_l\sr \tilde{H}^{*,*}(K_n)\oo_1{\bf F}_l$$
are isomorphic to the maps 
$${\cal A}^*\sr \tilde{H}^*(K_{2n}^{top})$$
defined by the action of the topological Steenrod algebra on $\iota_{2n}$ and one verifies easily that these isomorphisms identify the maps defined by (\ref{mor1}) with the similar maps defined by the topological suspension morphisms
$$\Sigma^2K_{2n}^{top}\sr K_{2n+2}^{top}.$$
We conclude that (\ref{topueq}) is isomorphic to the map  
\begin{eq}
\llabel{isoastar}
{\cal A}^*\sr {\lim_{\leftarrow}}{}_{n} \,\,\tilde{H}^{*+2n}(K_{2n}^{top})
\end{eq}
defined by the action of the topological Steenrod algebra on the canonical cohomology classes of $K_{2n}^{top}$ which is an isomorphism by \cite{SE}. 
\end{proof}
\begin{proposition}
\llabel{nou}
Let $k$ be an algebraically closed field of characteristic zero. Then the homomorphism 
$$u\oo_0 Id:{\cal A}^{*,*}\oo_{0}{\bf F}_l\sr M^{*,*}\oo_0{\bf F}_l$$ 
is a monomorphism.
\end{proposition}
\begin{proof}
Following the notations of \cite{MoVo} let $P^I$ denote the element of ${\cal A}^{*,*}$ corresponding to an admissible sequence $I=(\epsilon_0,s_1,\dots,s_k,\epsilon_k)$. We need to show that for any non-trivial linear combination $P=\sum a_I P^I$ there exists $n$ such that $P(\iota_n)\ne 0$ in $\tilde{H}^{*,*}(K_{n})\oo_0{\bf F}_l$. Bi-stability of operations together with the universal property of $\iota_n$ imply that it is sufficient to find any $X\in \Delta^{op}(Sm/k)^{\#}_+$ and a class $w\in \tilde{H}^{*,*}(X)$ such that $P(w)\ne 0$ in $\tilde{H}^{*,*}(X)\oo_0{\bf F}_l$. Using the computation of the action of motivic Steenrod algebra on the cohomology of $B\mu_l$ and the proof of \cite[Proposition VI.2.4]{SE} one can easily see that it can be done by taking $X=(B\mu_l)^N$ for some $N$ (cf. \cite[Proposition 11.4]{MoVo}).
\end{proof}

{\em End of the proof of Theorem \ref{mainst}:} To prove the theorem we need to show that the morphism $\tilde{u}$ of Lemma \ref{2009.09.18.l1} is an isomorphism. As was noted above, we may apply Corollary \ref{2009.09.17.cor3}(2) and assume that $k$ is an algebraically closed field.
A morphism between split Tate objects is an isomorphism if and only if it defines isomorphism on the motivic cohomology. Applying Lemma \ref{alg2} to the homomorphism
$$u:M^{*,*}\sr {\cal A}^{*,*}$$
considered as a homomorphism of graded modules over 
$$H^{*,*}(Spec(k),{\bf F}_l)={\bf F}_l[\tau]$$
with respect to the "weight" grading we see that it is sufficient to show that $u_1=u\oo_1 Id$ is  surjective and $u_0=u\oo_0 Id$ is injective. This is done in Propositions \ref{topu} and \ref{nou}.

\begin{remark}
\rm\llabel{false}
An unstable analog of Theorem \ref{mainst} appears to be false i.e. the motivic cohomology of individual spaces $K_n$ are not generated as algebras by elements of the form $P^I(\iota_n)$. In particular, \cite[Lemma 2.2]{zl} is probably false.

Let the base field be $\bf C$. In view of the unstable analog of Proposition \ref{topu} for any class $a\in \tilde{H}^{*,*}(K_n)$ there exist $m$ such that $\tau^ma$ can be represented as a polynomial of $P^I(\iota_n)$. We claim that there are classes $a$ for which the smallest $m$ satisfying this condition is $>0$.

If we consider all the motivic Eilenberg-MacLane spaces then it is obvious. Indeed, the motivic cohomology of weight zero of $K(\zz/l,n,0)$ are the same as the topological cohomology of $K(\zz/l,n)$. On the other hand all the motivic power operations shift the weight so applying $P^I$ to the canonical element in $\tilde{H}^{n,0}$ we get elements in $\tilde{H}^{p,q}$ with $q>0$. 

Existence of such classes for the spaces $K_n=K(\zz/l,2n,n)$ is less obvious. They do exist if the unstable motivic cohomology operations
$$P^i_{KM}:\tilde{H}^{p,q}\sr \tilde{H}^{p+2i(l-1),lq}$$
constructed in the context of the higher Chow groups by Kriz and May can be extended to the motivic cohomology of objects such as $K_n$. In that case the first example I know would be $w=Sq^{16}_{KM}Sq^7Sq^3Sq^1(\iota_3)$ for which one would expect 
$$\tau w =Sq^{16}Sq^7Sq^3Sq^1(\iota_3)$$
but which can not be obtained as a polynomial of $P^I(\iota_3)$ itself. The complexity of the example is due to the fact that one needs to find an admissible sequence with low excess and high weight shift. 

It is possible that the motivic cohomology classes of $K(\zz/l,p,q)$ for $p\ge 2q$ can be represented as polynomials of classes obtained from the canonical one using both stable and unstable power operations.
\end{remark}

\newpage \section{Appendices}
\subsection{Admissible categories}
\llabel{apadmissible}
\begin{definition}
\llabel{admissible}
A full subcategory $C$ of $Sch/k$ is called admissible if
\begin{enumerate}
\item $Spec(k)$ and $\af$ are in $C$
\item for $X$ and $Y$ in $C$ the product $X\times Y$ is in $C$
\item if $X$ is in $C$ and $U\sr X$ is etale then $U$ is in $C$
\item for $X$ and $Y$ in $C$ the coproduct $X\amalg Y$ is in $C$
\end{enumerate}
If in addition $C$ is closed under the formation of quotients with respect to actions of finite groups it will be called $f$-admissible. 
\end{definition}  
\begin{lemma}
\llabel{adm}
The categories of all schemes of finite type, of smooth quasi-projective schemes and of smooth quasi-affine schemes over any field are admissible.
\end{lemma}
\begin{lemma}
\llabel{qpfa}
The categories of quasi-projective and quasi-affine schemes over any field are $f$-admissible.
\end{lemma}
\begin{lemma}
\llabel{normok}
The categories of normal quasi-projective and normal quasi-affine schemes over a perfect field are $f$-admissible.
\end{lemma}
\begin{proof}
The only non-trivial point is to check that the product of two normal quasi-projective schemes over a perfect field is normal. This follows from \cite[6.8.5]{EGA4} and \cite[17.15.14.2]{EGA4}.
\end{proof}
\begin{lemma}
\llabel{semifad} 
The categories of quasi-projective and quasi-affine semi-normal schemes over a perfect field are $f$-admissible.
\end{lemma}
\begin{proof}
Let us consider for example the quasi-projective case. The  product of two  semi-normal schemes over a perfect field is semi-normal by \cite[Corollary 5.9]{traverso}. The fact that a scheme etale over a semi-normal one is semi-normal follows from the results of \cite{traverso} as well. Let us show that if $X$ is semi-normal then $X/G$ is semi-normal. Let $p:U\sr X/G$ be the semi-normalization of $X/G$. Then the projection $X\sr X/G$ factors through $p$ by the universal property of semi-normalizations and since $p$ is a universal homeomorphism and $X$ is reduced we conclude that $X\sr U$ is invariant under $G$-action. Hence we get an inverse $X/G\sr U$.
\end{proof}
\begin{remark}\rm{}
The categories of (semi-)normal quasi-projective and quasi-affine schemes over a non-perfect field are not admissible since the product of the spectra of two inseparable extensions need not be normal.
\end{remark}\rm{}
\begin{remark}\rm{}
Using the references provided above it is easy to see that the smallest admissible category over any field consists of disjoint unions of (smooth) schemes $X$ such that for some $n$ there exists an etale morphism $X\sr {\bf A}^n$. The smallest $f$-admissible category over a perfect field consists of finite group quotients of $X$ as above. I do not know whether the same category is $f$-admissible over any field.
\end{remark}\rm{}

\subsection{Finite group quotients in additive categories}

We will need some computations which apply to categorical quotients for finite group actions in any additive category $A$. In our case the category will be $Cor(C,R)$. Let $X$ be an object with an action of a finite group $G$. For an element $g\in G$ we let $[g]$ denote the corresponding automorphism $X\sr X$. For a subgroup $H$ of $G$ we let $p_H:X\sr X/H$ denote the projection. Note that $p_H[hg]=p_H[g]$ for any $h\in H$ and $g\in G$. 

If $L$ and $M$ are two subgroups of $G$ and $g\in G$ is such that $gLg^{-1}\subset M$ then there is a unique morphism 
$$p_{L,M,g}:X/L\sr X/M$$
such that 
\begin{eq}\llabel{udef1}
p_{L,M,g}p_L=p_M[g]
\end{eq}
We will write $p_{L,M}$ instead of $p_{L,M,1}$. Also for any $L\subset M$ in $G$ there exists a unique morphism
$$p^{M,L}:X/M\sr X/L$$
such that for any choice of a set of representatives $[L\backslash M]$ of $L\backslash M$ in $M$ one has 
\begin{eq}\llabel{udef2}
p^{M,L}p_M=\sum_{g\in [L\backslash M]}p_L[g].
\end{eq}
These observations are obvious from the definition of categorical quotients. We will need the following result.
\begin{proposition}
\llabel{start2}
For any $X$, $G$ and $H$ as above one has:
\begin{enumerate}
\item $p_{H,G}p^{G,H}=|G|/|H|\, Id_{X/G}$
\item Let $[H\backslash G/H]\subset G$ be a set of representatives for double coset classes in $G$ with respect to $H$. Then one has
\begin{eq}
\llabel{uff2}
p^{G,H}p_{H,G}=\sum_{x\in [H\backslash G/H]}p_{H\cap x^{-1}Hx,H,x}p^{H,H\cap x^{-1}Hx}
\end{eq}
\end{enumerate}
\end{proposition}
\begin{proof}
The first equality is between morphisms from $X/G$ to itself. Therefore it is sufficient to check that the compositions of these morphisms with $p_G$ coincide. We have:
$$p_{H,G,1}p^{G,H}p_G=p_{H,G,1}\sum_{g\in [H\backslash G]} p_{H}[g]=\sum_{g\in [H\backslash G]} p_G[g]=|G|/|H|p_G$$
where the first equality holds by (\ref{udef2}), the second by (\ref{udef1}) and the third because $p_G[g]=p_G$ for all $g$.

To prove the second statement let us choose a set of representatives of double coset classes $[H\backslash G/H]\subset G$. Let us further choose sets of representatives $[(H\cap x^{-1}Hx)\backslash H]\subset H$ for all $x\in [H\backslash G/H]$. Then an element $g$ of $G$ can be written in a unique way as the product $hxu$ where $h\in H$, $x\in [H\backslash G/H]$ and $u\in [(H\cap x^{-1}Hx)\backslash H]$. In particular, elements of the form $xu$ for $x\in [H\backslash G/H]$ and $u\in [(H\cap x^{-1}Hx)\backslash H]$ give us a set of representatives for $H\backslash G$ which we denote by $A$. 

Since (\ref{uff2}) is an equality between two morphisms from $X/H$ to itself it is sufficient to check that their compositions with $p_H$ coincide i.e. that
\begin{eq}
\llabel{uff3}
p^{G,H}p_{H,G,1}p_H=\sum_{x\in [H\backslash G/H]}p_{H\cap x^{-1}Hx,H,x}p^{H,H\cap x^{-1}Hx}p_H
\end{eq}
We have
$$\sum_{x\in [H\backslash G/H]}p_{H\cap x^{-1}Hx,H,x}p^{H,H\cap x^{-1}Hx}p_H=\sum_{x\in [H\backslash G/H]}p_{H\cap x^{-1}Hx,H,x}\sum_{u\in [(H\cap x^{-1}Hx)\backslash H]}p_{H\cap x^{-1}Hx}[u]=$$
$$=\sum_{x\in [H\backslash G/H]}\sum_{u\in [(H\cap x^{-1}Hx)\backslash H]}p_H[x][u]=\sum_{x\in [H\backslash G/H]}\sum_{u\in [(H\cap x^{-1}Hx)\backslash H]}p_H[xu]$$
where the first equality holds by (\ref{udef2}), the second by (\ref{udef1}) and the third because $[x][u]=[xu]$. 
On the other hand
$$p^{G,H}p_{H,G,1}p_H=p^{G,H}p_G=\sum_{g\in A}p_H[g]=\sum_{x\in [H\backslash G/H]}\sum_{u\in [(H\cap x^{-1}Hx)\backslash H]}p_H[xu].$$
\end{proof}
\begin{cor}
\llabel{uffnorm}
If $H$ is a normal subgroup in $G$ and $[H\backslash G]$ is a set of representatives for the left conjugacy classes of $H$ in $G$ then 
\begin{eq}
\llabel{uff32}
p^{G,H}p_{H,G}=\sum_{g\in [H\backslash G]}p_{H,H,g} 
\end{eq}
\end{cor}

Recall from \cite[p. 149]{McLane2} that a fork diagram $\partial_0,\partial_1:X\dsr Y\stackrel{e}{\sr} Z$ in a category is called a split fork diagram if $e\partial_0=e\partial_1$ and there exist morphisms $s:Z\sr Y$ and $t:Y\sr X$ such that $es=Id$ and $\partial_0t=Id$, $\partial_1t=se$. As was shown in \cite[p. 149]{McLane2}, every split fork diagram is an absolute coequalizer diagram. 
\begin{proposition}
\llabel{splitc}
Let $X$, $G$, $H$ be as in Proposition \ref{start2} and assume in addition that the index $|G|/|H|$ of $H$ in $G$ is invertible in the $Hom$-groups of our category. Then the fork diagram
$$\partial_0,\partial_1:\oplus_{g\in G} X/(H\cap g^{-1}Hg)\dsr X/H\stackrel{e}{\sr} X/G$$
where 
$$\partial_0=\oplus_{g\in G} p_{H\cap g^{-1}Hg,H}$$
$$\partial_1=\oplus_{g\in G} p_{H\cap g^{-1}Hg,H,g}$$
and $e=p_{H,G}$, is a split fork diagram and in particular an absolute coequalizer diagram. 
\end{proposition}
\begin{proof}
As in the proof of Proposition \ref{start2} let us choose a set of representatives for the double cosets $[H\backslash G/H]\subset G$. 
$$s=(|H|/|G|)p^{G,H}$$
$$t=(|H|/|G|)\oplus_{x\in [H\backslash G/H]}p^{H,H\cap x^{-1}Hx}$$
The relation $es=Id$ follows from Proposition \ref{start2}(1) and by Proposition \ref{start2}(2) we have
$$se=(|H|/|G|)p_{H,G}p^{G,H}=(|H|/|G|)\sum_{x\in [H\backslash G/H]}p_{H\cap x^{-1}Hx,H,x}p^{H, H\cap x^{-1}Hx}$$
We further have
$$\partial_0 t=(|H|/|G|)\oplus_{x\in [H\backslash G/H]}p_{H\cap x^{-1}Hx,H}p^{H,H\cap x^{-1}Hx}=$$
$$=(|H|/|G|)\sum_{x\in [H\backslash G/H]}(|H|/|H\cap x^{-1}Hx|)Id_{X/H}=(|H|/|G|)|H\backslash G|Id_{X/H}=Id_{X/H}$$
and
$$\partial_1 t=(|H|/|G|)\oplus_{x\in [H\backslash G/H]}p_{H\cap x^{-1}Hx,H,x}p^{H,H\cap x^{-1}Hx}=se.$$
\end{proof}
\begin{remark}\rm
The proof of Proposition \ref{splitc} shows that it remains valid if we replace the first term of the fork by $\oplus_{g\in [H\backslash G/H]} X/(H\cap g^{-1}Hg)$ for any set of representatives of double cosets $[H\backslash G/H]\subset G$.
\end{remark}

\subsection{Radditive functors}
\llabel{radsum}
In this appendix we reproduce for the convenience of the reader some definitions and results of \cite{SRFsub} which we use in this paper.

Let $C$ be a small category with finite coproducts. A radditive functor on $C$ is a functor $C^{op}\sr Sets$  (i.e. a presheaf of sets) such that $F(X\amalg Y)=F(X)\times F(Y)$. We denote the category of radditive functors by $Rad(C)$. Representable functors are radditive by definition of coproducts. The category $Rad(C)$ is complete and cocomplete with limits, filtered colimits and reflexive coequalizers (i.e. coequalizers of pairs of arrows $f,g$ which have a common section $s$, $gs=fs=Id$) being the same as for presheaves but with different coproducts which are compatible for representable functors with coproducts in $C$ (see \cite[Prop. {2009elprop}]{SRFsub}). 

The category of radditive functors on an additive category $A$ is naturally equivalent to the category of additive contravariant functors from $A$ to the category of abelian groups. If $C$ has a final object $pt$ and $C_+$ is the full subcategory of the category of pointed objects in $C$ which consists of objects of the form $(X\amalg pt, i_{pt})$ then the category of radditive functors on $C_+$ is naturally equivalent to the category of pointed radditive functors on $C$ (see \cite[Lemma {2007pointed}]{SRFsub}). If $C$ is a formal closure with respect of finite coproducts of a subcategory $C'$ then the category of radditive functors on $C$ is naturally equivalent to the category of presheaves on $C'$. 

A morphism of simplicial radditive functors $F:X\sr Y$ is called a projective equivalence (resp. projective fibration) if for any $U\in C$ the map of simplicial sets $F_U:X(U)\sr Y(U)$ is a weak equivalence (resp. Kan fibration). The classes $W_{proj}$ and $Fib_{proj}$ of projective equivalences and projective fibrations generate a closed model structure on $\Delta^{op}Rad(C)$ which is called the projective closed model structure. Its homotopy category is denoted by $H(C)$.

For any set of morphisms $E$ in $\Delta^{op}Rad(C)$ one defines in the usual way the class of $E$-local objects and the class of $E$-local equivalences $cl_l(E)$. The localization of $H(C)$ by $cl_l(E)$ always exists and we denote it by $H(C,E)$. If the projective closed model structure is left proper then there exists a left Bousfield localization (see \cite{Hirs}) of the projective closed model structure by $E$ and $cl_l(E)$ coincides with its class of weak equivalences. In the case when the class of projective equivalences is closed under finite coproducts, which holds in all of the examples mentioned above,  the projective closed model structure is left proper and in addition the class of $E$-local equivalences is closed under coproducts for any $E$.

For any functor $F:C\sr Rad(C')$ where $C$ and $C'$ are categories with finite coproducts the Kan extension $F^*$ of $F$ takes radditive functors to radditive functors and defines a $F^{rad}:Rad(C)\sr Rad(C')$ which we call the radditive extension of $F$.  The radditive extensions commute with filtered colimits and reflexive coequalizers. If $F$ commutes with finite coproducts then the right adjoint $F_*$ to $F^*$ takes radditive functors to radditive functors and defines a right adjoint $F_{rad}$ to $F^{rad}$ which also commutes with filtered colimits and reflexive coequalizers. 

Let $C^{\#}$ be the full subcategory of $Rad(C)$ which consists of filtered colimits of representable functors and $\bar{C}$ the full subcategory of $Rad(C)$ which consists of coproducts of representable functors (the category $C^{\#}$ is also known as the category of ind-objects over $C$. See e.g. \cite[Sec. 8.2.4, p.70]{SGA41}). Since a finite coproduct of representable radditive functors is representable we have $\bar{C}\subset C^{\#}$. 

\begin{proposition}[{\cite[Prop. 3.18 \comment{2009lres3}]{SRFsub}}]
\llabel{2009ref0}
There are a functor $L_*:Rad(C)\sr \Delta^{op}Rad(C)$ and a natural transformation $L_*\sr \iota$ where $\iota:Rad(C)\sr \Delta^{op}Rad(C)$ is the natural embedding,  such that
\begin{enumerate}
\item for any $X\in \Delta^{op}Rad(C)$ the object $L_*(X)$ belongs to $\Delta^{op}\bar{C}$ and the morphism $L_*(X)\sr X$ is projective equivalence,
\item for any $X\in \Delta^{op}C$ the morphism $L_*(X)\sr X$ is a simplicial homotopy equivalence.
\end{enumerate}
\end{proposition}
\begin{proposition}[{\cite[Prop. 4.10 \comment{preloc}, Th. 4.8 \comment{2009th2}]{SRFsub}}]
\llabel{2009ref1}
One has:
\begin{enumerate}
\item the functor $\Delta^{op}C^{\#}\sr H(C)$ is a localization,
\item for any $F:C\sr Rad(C')$, the functor $F^{rad}$ takes projective equivalences between objects of $\Delta^{op}C^{\#}$ to projective equivalences.
\end{enumerate}
\end{proposition}
Proposition \ref{2009ref1} implies that any $F:C\sr Rad(C')$ defines a functor $H(C)\sr H(C')$ which we denote by $\LL F^{rad}$. 
\begin{proposition}[{\cite[Lemma 4.12 \comment{2007fcc}]{SRFsub}}]
\llabel{2009ref2}
Let $F:C\sr C'$ be a functor which commutes with finite coproducts. Then the functor $F_{rad}$ takes projective equivalences to projective equivalences and the resulting functor ${\bf R} F_{rad}:H(C')\sr H(C)$ is right adjoint to $\LL F^{rad}$.
\end{proposition}
\begin{theorem}[{\cite[Th. 4.19 \comment{2007el1}]{SRFsub}}]
\llabel{2009ref3}
Let $F:C\sr Rad(C')$ be a functor and $E$, $E'$ be sets of morphisms in $\Delta^{op}C^{\#}$ and $\Delta^{op}Rad(C')$ respectively such that for any $f\in E$ and $U\in C$ one has
$$F^{rad}(f\amalg Id_U)\in cl_l(E')$$
Then for any $f\in cl_l(E)\cap \Delta^{op}C^{\#}$ one has $F^{rad}(f)\in cl_l(E')$ and in particular $\LL F^{rad}$ defines a functor $H(C,E)\sr H(C',E')$. 
\end{theorem}
\begin{theorem}[{\cite[Th. 4.20 \comment{2007eadj}]{SRFsub}}]
\llabel{2009ref4}
Let $F:C\sr C'$ be a functor which commutes with finite coproducts and $E$, $E'$ be sets of morphisms in $\Delta^{op}C^{\#}$ and $\Delta^{op}(C')^{\#}$ respectively such that for any $f\in E$ and $U\in C$ one has $F^{rad}(f\amalg Id_U)\in cl_l(E')$ and for any $f'\in E'$ and $U'\in C'$ one has $F_{rad}(f'\amalg Id_{U'})\in cl_l(E)$. 

Then for any $f\in cl_l(E)\cap \Delta^{op}C^{\#}$ one has $F^{rad}(f)\in cl_l(E')$ and for any $f'\in cl_l(E')$ one has $F_{rad}(f')\in cl_l(E)$ and the resulting functors $\LL_E F^{rad}$ and ${\bf R}_E F_{rad}$ between $H(C,E)$ and $H(C,E')$ are adjoint.
\end{theorem}
\begin{proposition}[{\cite[Cor. 4.21 \comment{2007anfe}]{SRFsub}}]
\llabel{2009ref5}
Under the assumptions of Theorem \ref{2009ref4} one has:
\begin{enumerate}
\item if $F$ is a full embedding then $\LL_E F^{rad}$ is a full embedding and ${\bf R}_E F_{rad}$ is a localization,
\item if $F$ is surjective on isomorphism classes of objects then ${\bf R}_E F_{rad}$ reflects isomorphisms.
\end{enumerate}
\end{proposition}
A morphism $X\sr Y$ in $Rad(C)$ is called a coprojection if it is isomorphic to the canonical morphism of the form $X\sr X\amalg A$. A morphism $f$ in $\Delta^{op}Rad(C)$ is called a term-wise coprojection if each term of $f_i$ of $f$ is a coprojection. 

A sequence $X\sr Y\sr Z$ is called a term-wise coprojection sequence if $X\sr Y$ is a coprojection and $Y\sr Z$ is isomorphic to $Y\sr Y/X$.  For the general notion of a cofiber sequence used below see \cite{Quillen} or \cite{Hovey}.
\begin{proposition}[{\cite[Th. 3.41 \comment{2009cfs}, Cor. 3.53 \comment{2009cs1}]{SRFsub}}]
\llabel{2009ref6}
If $C$ is a pointed category then:
\begin{enumerate}
\item any term-wise coprojection sequence $(X\sr Y\sr Z)$ in $\Delta^{op}C^{\#}$ extends in a natural way to a cofiber sequence $(X\sr Y\sr Z, \, Z\sr X\vee_{\LL}\Sigma(X))$ in $H(C)$, 
\item any cofiber sequence $(X\sr Y\sr Z, \, Z\sr X\vee_{\LL}\Sigma(X))$ in $H(C)$ is isomorphic to the cofiber sequence defined by a coprojection sequence $X'\sr Y'\sr Z'$ in $\Delta^{op}\bar{C}$ where $\bar{C}$ is the full subcategory of coproducts of representable functors in $Rad(C)$,
\item for any cofiber sequence $(X\stackrel{f}{\sr} Y\stackrel{g}{\sr} Z, \, Z\sr X\vee_{\LL}\Sigma(X))$ in $H(C)$ one has $g\in cl_{l}(\{X\sr pt\})$ and $(pt\sr Z)\in cl_{l}(\{f\})$.
\end{enumerate}
\end{proposition}

The main technical tool which is used in the proofs of the results cited above is the notion of a $\bdl$-closed class defined as follows. Let $D$ be a category which has coproducts and filtered  colimits. For a set $K$ and an object $X$ of $D$ we let $X\oo K=\amalg_K X$ denote the coproduct of $K$ copies of $X$. Similarly for a simplicial set $K$ and an object $X$ of $\Delta^{op}D$ we let   $X\oo K$ denote the simplicial object with terms $X_n\oo K_n$. 

Let $f,g:X\sr Y$ be two morphisms in $\Delta^{op}D$. An elementary simplicial homotopy from $f$ to $g$ is defined in the usual way as a morphism $h:X\oo\Delta^1\sr Y$ such that $f=h\circ (Id_X\oo i_0)$ and $g=h\circ (Id_X\oo i_1)$ where $i_0,i_1$ are the standard morphisms $\Delta^0\sr \Delta^1$. Let $\sim$ be the smallest equivalence relation on morphisms such that $f\sim g$ if there exists an elementary simplicial homotopy from $f$ to $g$.  A morphism $f:X\sr Y$ in $\Delta^{op}D$ is said to be a simplicial homotopy equivalence defined if there exists a morphism $g:Y\sr X$ such that $f\circ g\sim Id_Y$ and $g\circ f\sim Id_X$. 

A class of morphisms $E$ in a category $\Delta^{op}D$ is said to be closed under filtered colimits if for any pair of filtered systems $(X_i)_{i\in I}$, $(Y_i)_{i\in I}$ and any  morphism of systems $(f_i):(X_i)\sr (Y_i)$ such that $f_i\in E$ one has $f=colim_i f_i\in E$.
\begin{definition}
\llabel{bdlclosed}
A class of morphisms $E$ in $\Delta^{op}D$ is said to be $\bdl$-closed if it satisfies the following conditions:
\begin{enumerate}
\item simplicial homotopy equivalences are in $E$,
\item if $f$ and $g$ are morphisms such that the composition $gf$ is
defined and two out of three morphisms $f, g, gf$ are in $E$ then the
third is in $E$,
\item if $f:B\sr B'$ is a morphism of bisimplicial objects over $D$
such that the rows or columns of $f$ are in $E$ then the diagonal morphism
$\Delta(f)$ is in $E$,
\item $E$ is closed under filtered colimits.
\end{enumerate}
\end{definition}
For any class of morphisms $E$ there exists the smallest $\bdl$-closed class which contains $E$ and we denote it by $cl_{\bdl}(E)$.

\begin{theorem}[{\cite[Th. 3.51 \comment{2009main}, Cor. 3.52 \comment{2009main2}]{SRFsub}}]
\llabel{maindelta}
Let $C$ be a small category with finite coproducts and $E$ a set of morphisms in $\Delta^{op}C^{\#}$. Then one has
$$cl_l(E)\cap \Delta^{op}C^{\#}=cl_{\bdl}(E\amalg Id_C)$$
and
$$cl_l(E)=cl_{\bdl}((E\amalg Id_C)\cup W_{proj})$$
where $E\amalg Id_C$ is the set of morphisms of the form $f\amalg Id_U$ for $f\in E$ and $U\in C$.
\end{theorem}
The following somewhat technical result is used in the proof of Proposition \ref{deltaver} below. We denote by $CofEnds$ the class of morphisms between objects in $\Delta^{op}Rad(C)$ which are cofibrant in the projective closed model structure. 
\begin{proposition}[{\cite[Prop. 3.28 \comment{2009.07.23.1}]{SRFsub}}]
\llabel{technical1}
Let $C$ be a small category with finite coproducts and $E$ a class of morphisms in $\Delta^{op}Rad(C)$ which satisfies the following conditions:
\begin{enumerate}
\item $E$ contains $W_{proj}$,
\item $E$ satisfies the 2-out-of-3 property,
\item $E\cap \Delta^{op}C^{\#}\cap CofEnds$ is closed under coproducts,
\item for $f\in E\cap \Delta^{op}C^{\#}\cap CofEnds$ and $i\ge 0$ one has $f\oo Id_{\partial\Delta^i}\in E$,
\item for a morphism of push-out squares
$$
\left(
\begin{CD}
f_1\,\, @. f_2\\
f_3\,\, @. f_4
\end{CD}
\right)
:
\left(
\begin{CD}
X_1 @>>> X_2\\
@VgVV @VVV\\
X_3 @>>> X_4
\end{CD}
\right)
{\longrightarrow}
\left(
\begin{CD}
X_1' @>>>X_2'\\
@Vg'VV @VVV\\
X_3' @>>> X_4'
\end{CD}
\right)
$$
such that all the objects are in $\Delta^{op}C^{\#}\cap CofEnds$, the morphisms $g$, $g'$ are cofibrations and term-wise coprojections and $f_1,f_2,f_3$ are in $E$ one has $f_4\in E$,
\end{enumerate}
Then $E$ is $\bdl$-closed.  
\end{proposition}

Let $A$ be an additive category. The category $Rad(A)$ of radditive functors on $A$ is easily seen to be equivalent to the category of additive functors from $A$ to the category of abelian groups. In particular it is an abelian category. The objects of the smallest subcategory $\bar{A}$ of $Rad(A)$ which contains $A$ and is closed under direct sums are projective in $Rad(A)$ and any object from $Rad(A)$ can be epimorphically covered by an object from $\bar{A}$. In particular $\bar{A}$ has enough projectives. Let $Cmpl_{-}(Rad(A))$ be the category of complexes bounded from the above over $Rad(A)$. We have the usual Dold-Kan correspondence i.e. an adjoint pair of functors
$$N:\Delta^{op}Rad(A)\sr Cmpl_{-}(Rad(A))$$
$$K:Cmpl_{-}(Rad(A))\sr \Delta^{op}Rad(A)$$
where $N$ is the normalization functor and $K$ is its right adjoint. As for any abelian category, the functor $N$ is a full embedding and $K$ is a localization. The composition $N\circ K$ coincides with the canonical truncation functor $\tau^{\le 0}$ which "removes" the negative cohomology objects of a complex. 

Due to the fact that the representable functors are projective in $Rad(C)$, a morphism $f$ in $\Delta^{op}Rad(A)$ is a projective equivalence if and only if $N(f)$ is a quasi-isomorphism. Therefore $N$ takes projective equivalences to quasi-isomorphisms and $K$ takes quasi-isomorphisms to projective equivalences. After passing to the corresponding localizations we obtain a pair of adjoint functors
$$N_{proj}:H(A)\sr D_{-}(Rad(A))$$
$$K_{proj}:D_{-}(Rad(A))\sr H(A)$$
such that $N_{proj}$ is a full embedding. 

For a class $E$ of morphisms in a triangulated category let $cl_{vl}(E)$ denote the (left) Verdier closure of $E$ i.e. the class of morphisms whose cones belong to the localizing subcategory generated by cones of morphisms from $E$. 

\begin{proposition}
\llabel{deltaver}
Let $E$ be a set of morphisms in $\Delta^{op}Rad(A)$. Then one has
$$N_{proj}(cl_l(E))\subset cl_{vl}(N_{proj}(E))$$
(for simplicity we omit the natural projection $\Phi:\Delta^{op}Rad(A)\sr H(A)$ from our notation).
\end{proposition}
\begin{proof}
By Theorem \ref{maindelta} we have
\begin{eq}
\llabel{mainadd}
cl_l(E)=cl_{\bdl}((E\oplus Id_A)\cup W_{proj})
\end{eq}
where we use $\oplus$ to denote coproducts in the additive context. Consider the class $F=N_{proj}^{-1}(cl_{vl}(N_{proj}(E)))$. To prove that it contains $cl_l(E)$ it is sufficient in view of (\ref{mainadd}) to prove that it is $\bdl$-closed and contains $W_{proj}$ and $E\oplus{Id_A}$. The 
later is obvious. To prove that $F$ is $\bdl$-closed we use \cite[Prop. {bdlchar}]{SRFsub} which implies easily that it is sufficient to establish that $F$ satisfies the diagonal condition of Definition \ref{bdlclosed}. Let us apply Proposition \ref{technical1} to $F$. the first three conditions of the lemma are obvious. The fourth condition follows from the fact that $N(X\oo\partial\Delta^i)\cong N(X)[i-1]$. To prove the fifth condition observe that for a push-out square
$$\begin{CD}
X_1 @>>> X_2\\
@VgVV @VVV\\
X_3 @>>> X_4
\end{CD}$$
such that $g$ is a monomorphism, the sequence of complexes 
$$0\sr N(X_1)\sr N(X_2)\oplus N(X_3)\sr N(X_4)\sr 0$$
is exact. Since coprojections in an additive category are monomorphisms this implies that in the context of the fifth condition of Proposition \ref{technical1} we get a morphism of exact sequences of complexes 
$$\begin{CD}
0\sr N(X_1) @>>> N(X_2)\oplus N(X_3) @>>> N(X_4) \sr 0\\
@VN(f_1)VV @VN(f_2)\oplus N(f_3)VV @VN(f_4)VV \\
0\sr N(X_1') @>>> N(X_2')\oplus N(X_3') @>>> N(X_4') \sr 0
\end{CD}
$$
and therefore a distinguished triangle
$$cone(N(f_1))\sr cone(N(f_2))\oplus cone(N(f_3))\sr cone(N(f_4))\sr cone(N(f_1))[1]$$
The fact that $N(f_4)\in F$ follows now from the standard properties of Verdier closure. Proposition is proved. 
\end{proof}

\def\cprime{$'$}

\comment{
\bibliography{alggeom}
\bibliographystyle{plain}
}

\end{document}